\documentclass[11pt,a4paper]{amsart} 
\usepackage{amsfonts,amssymb,bbm} 
\usepackage[cmtip,all]{xy}                     
\usepackage{enumerate,mathrsfs}

\usepackage{color}
\usepackage[colorlinks,plainpages]{hyperref}

\makeatletter
\def\subsection{\@startsection{subsection}{2}%
  \z@{.5\linespacing\@plus.7\linespacing}{-.5em}%
  {\normalfont\bfseries}}
   \renewcommand{\@secnumfont}{\bfseries}
\def\subsubsection{\@startsection{subsubsection}{3}%
 \z@{.5\linespacing\@plus.7\linespacing}{-.5em}%
  {\normalfont\itshape}}
\makeatother

\def\C{\Bbb{C}}

\def\k{\mathbbm{k}}

\def\Q{\Bbb{Q}}
\def\R{\Bbb{R}}

\def\di{\partial}

\def\bl{\langle}
\def\br{\rangle}

\def\liml{\lim\limits}
\def\suml{\sum\limits}

\def\capl{\mathop\cap\limits}

\def\tg{\tilde{g}}

\def\tN{{\tilde{N}}}

\def\tI{\tilde{I}}

\def\tz{{\tilde{z}}}

\def\hf{\hat{f}}
\def\hG{{\widehat{G}}}
\def\hg{{\hat{g}}}

\def\hM{{\widehat{M}}}
\def\hR{{\widehat{R}}}

\def\hU{\widehat{U}}
\def\hV{\widehat{V}}

\def\tphi{\tilde{\phi}}

\def\hw{\hat{w}}
\def\hx{\hat{x}}
\def\txi{{\tilde{\xi}}}
\def\hy{\hat{y}}
\def\hz{\hat{z}}

\def\al{\alpha}\def\Ga{\Gamma}\def\be{\beta}\def\De{\Delta}

\def\La{\Lambda}

\def\cA{\mathcal A}\def\ca{\mathfrak a}

\def\cK{{\mathscr K}\!}
\def\cO{\mathcal O}

\def\cR{\mathscr{R}}

\def\cm{{\frak m}}
\def\cp{{\frak p}}

\def\um{{\underline{m}}}

\def\one{{1\hspace{-0.1cm}\rm I}}\def\zero{\mathbb{O}}

\newcommand{\ber}{\begin{array}{l}}\newcommand{\eer}{\end{array}}
\newcommand{\bpm}{\begin{pmatrix}}\newcommand{\epm}{\end{pmatrix}}
\newcommand{\bbm}{\begin{bmatrix}}\newcommand{\ebm}{\end{bmatrix}}
\newcommand{\bM}{\begin{matrix}}\newcommand{\eM}{\end{matrix}}
\newcommand{\bee}{\begin{enumerate}}\newcommand{\eee}{\end{enumerate}}
\newcommand{\bei}{\begin{itemize}}\newcommand{\eei}{\end{itemize}}

\def\wrt{with respect to }
\def\sset{\subset}\def\sseteq{\subseteq}\def\smin{\setminus}

\def\Mat{{Mat}_{m\times n}(R)}

\newtheorem{Lemma}{Lemma}[section]\newcommand{\bel}{\begin{Lemma}}\newcommand{\eel}{\end{Lemma}}
\newtheorem{Theorem}[Lemma]{Theorem}\newcommand{\bthe}{\begin{Theorem}}\newcommand{\ethe}{\end{Theorem}}
\newtheorem{Proposition}[Lemma]{Proposition}\newcommand{\bprop}{\begin{Proposition}}\newcommand{\eprop}{\end{Proposition}}
\newtheorem{Corollary}[Lemma]{Corollary}\newcommand{\bcor}{\begin{Corollary}}\newcommand{\ecor}{\end{Corollary}}
\newtheorem{Assumptions}[Lemma]{Assumptions}
\newtheorem{Definition}[Lemma]{Definition}
\newcommand{\bed}{\begin{Definition}}
\newcommand{\eed}{\end{Definition}}
\newtheorem{Definition-Proposition}[Lemma]{Definition-Proposition}

\def\bpr{~\\{\em Proof.\ }}
\newcommand{\epr}{{\hfill\ensuremath\blacksquare}}
\newtheorem{Remark}[Lemma]{Remark}
\newcommand{\beR}{\begin{Remark}\rm}
\newcommand{\eeR}{\end{Remark}}
\newtheorem{Example}[Lemma]{Example}
\newcommand{\bex}{\begin{Example}\rm}
\newcommand{\eex}{\end{Example}}
\newtheorem{Problem}[Lemma]{Problem}
\newcommand{\bprob}{\begin{Problem}\rm}
\newcommand{\eprob}{\end{Problem}}

\newcommand{\bet}{\begin{tabular}{cccccccc}}\newcommand{\eet}{\end{tabular}}
\newcommand{\beq}{\begin{equation}}\newcommand{\eeq}{\end{equation}}

\newcommand{\bin}[2]{\binom{#1}{#2}}

\title[]{P\MakeLowercase{airs of} L\MakeLowercase{ie-type and  large orbits of group actions on filtered modules.}\\
(A \MakeLowercase{characteristic-free approach to finite determinacy}.)}
\author{A\MakeLowercase{lberto} F. B\MakeLowercase{oix}, G\MakeLowercase{ert}--M\MakeLowercase{artin} G\MakeLowercase{reuel and}
D\MakeLowercase{mitry} K\MakeLowercase{erner}}
\address{Department of Applied Mathematics, Science, Materials Engineering and Electronic Technology, Universidad Rey Juan Carlos, C/Tulip\'an, s/n, 28933, M\'ostoles, Spain.}
\address{Fachbereich Mathematik, Universit\"at Kaiserslautern, Erwin--Schr\"odinger Str., 67663 Kaiserslautern, Germany.}
\address{Department of Mathematics, Ben Gurion University of the Negev, P.O.B. 653, Be'er Sheva 84105, Israel.}
\email{alberto.fboix@urjc.es}
\email{greuel@mathematik.uni-kl.de}
\email{dmitry.kerner@gmail.com}
\thanks{A.F. Boix was partially supported by Israel Science Foundation (grant No. 844/14) and Spanish Ministerio de Econom\'ia y Competitividad
 MTM2016-7881-P and PID2019-104844GB-I00}
\thanks{D.Kerner was partially supported by Israel Science Foundation (grant No. 1910/18)}

\begin{document}
\maketitle

\begin{abstract}
Finite determinacy for mappings has been classically thoroughly studied  in numerous scenarios in the real- and
complex-analytic category and in the differentiable case. It means that the map-germ is determined, up to a
given equivalence relation, by a finite part of its Taylor expansion.
The equivalence relation is usually given by a group action and the first step is always to reduce the determinacy
question to an ``infinitesimal determinacy", i.e., to the tangent spaces at the orbits of the group action.
In this work we formulate a universal, characteristic-free approach to finite determinacy,
not necessarily
over a field, and for a large class of group actions. We do not restrict to pro-algebraic or Lie groups, rather we
introduce the notion of ``pairs of (weak) Lie type", which are groups together with a substitute for the tangent space to the orbit
such that the orbit
is locally approximated by its tangent space, in a precise sense. This construction may be considered as a kind of replacement of the exponential resp. logarithmic maps. It is of independent interest as it provides a general method to pass from the tangent space to the orbit of a group action in any characteristic.
In this generality we establish the ``determinacy versus infinitesimal determinacy" criteria, a far reaching
generalization of numerous classical and recent results,
together with some new applications.
\end{abstract}
\bigskip

\noindent {\bf Key words:} {finite determinacy, pairs of Lie-type, sufficiency of jets, infinitesimal stability, group actions on modules.}

\noindent{\bf 2020 Mathematics Subject Classification: }{Primary: 58K40, 58K50. Secondary: 14B05, 14B07, 16W70. }

\setcounter{secnumdepth}{6} \setcounter{tocdepth}{1}
\tableofcontents

\section{Introduction}
\subsection {Prologue}

Let $f$ be the germ at the origin of a real- or complex-analytic
function or a $C^\infty$-function of several
variables $x =(x_1, \ldots, x_p)$. Famous results of
\cite{Mather1968}, \cite{Tougeron1968}  and many others on finite
determinacy of function-germs bound the order of determinacy in
terms of the jacobian ideal  of $f$: \bei
\item
if $\cm^2\cdot Jac(f)\supseteq \cm^{N+1}$ then $f$ is $N$-right-determined,
\item
if $\cm^2\cdot Jac(f)+\cm\bl f\br\supseteq \cm^{N+1}$ then $f$ is $N$-contact-determined.
\eei
Here $\cm=\bl x_1,\dots,x_p\br$ is the maximal ideal and $Jac(f)$ is the ideal generated by the partials of $f$.
 A function $f$ is $N$-right (resp. $N$-contact) determined if every $g$ whose Taylor expansion up to order $N$ coincides
  with that of $f$ lies in the same $\cR$- (resp. $\cK$-) orbit as $f$. Here
$\cR$ (resp. $\cK$) is the {\em right group} (resp. {\em contact group}) acting on the ring of germs by analytic or $C^\infty$-coordinate change (resp. additionally by multiplication with a unit).

\medskip

These results have been generalized to numerous group actions and
rings in the following way. Let $M$
be a  space of maps, usually a filtered module over a ring, together with a
fixed action $G\circlearrowright M$ of a nice subgroup $G$ of the
contact group $\cK$. The classical statements compare the tangent
space $T_{(Gf,f)}$ to the group orbit $Gf$ at $f$ with the filtration of $M$
given by the subspaces $\cm^{i} \cdot M$. Let  $f\in M$ be some element and consider the statement:
\beq\label{Eq.Mather.classical.theorem}
\ber
 \text{\em i.  \ Suppose the tangent space to the orbit of  $f$ satifies:
 \quad $\cm\cdot T_{(Gf,f)}\supseteq\cm^{N+1}\cdot M.$}
\\ \text{\em  ii. Then the orbit of $f$ is large in the sense: \quad $Gf\supseteq\{f\}+\cm^{d_N+1}\cdot M$,}
\eer
\eeq
where $d_N$ is some integer depending on $N$.
 Whenever the statement \eqref{Eq.Mather.classical.theorem}$ii.$ holds one says that $f$ is $d_N$-determined \wrt the $G$-action.

\medskip

A statement like \eqref{Eq.Mather.classical.theorem}  can be rephrased in saying ``a large tangent space implies a large orbit''.
 To prove such a statement basically two different methods have been used. Primarily
 the integration of vector fields and the use of the exponential map in characteristic zero with the space of maps $M$ involving formal or analytic power series or germs of $C^\infty$-maps and with $G$  an algebraic group or a Lie group (after reduction
  to a finite dimensional parameter space). Secondly, power series methods with $M$ involving formal power series over a field of positive characteristic.
   However, in different scenarios for different kinds of $M$ and groups, these methods had always to be adapted and modified.
   One of the aims of our paper is to give a unified approach.

   \medskip

\subsection{Goals and methods} \label{Sec.Intro.Goals}
The goals of our current work are three-fold:
\bee[ \ \ A.]
\item To extend the whole theory to arbitrary filtered modules $M$ over a base ring $\k$ of  any characteristic, not necessarily a field.
\item To broaden the class of admissible group actions in a characteristic-free way, by combining the characteristic zero approach and the use of the exponential map with the power series approach in positive characteristic.
\item To show how the general theory may be applied, not only to recover most
of the previously known results, but also to obtain some new ones.
\eee

\medskip

To deduce the inclusion as in  \eqref{Eq.Mather.classical.theorem}$(ii)$ from  the assumption \eqref{Eq.Mather.classical.theorem}$(i)$
we study the orbit $Gz$ of an element $z \in M$. For a given ``higher order" element $w\in M$ we want to prove $z+w\in Gz$, i.e.,
 to solve the equation $z+w=gz$, for the unknown $g\in G$.
This is done in two steps, as follows.

\medskip

\quad {\em Step 1.} First one establishes an ``order-by-order" solution, i.e., a sequence $\{g_n\}$ of elements of $G$ satisfying $g_nz\to z+w$.
 The convergence is taken in the filtration topology with the limit being an element in the closure of the orbit, $z+w\in \overline{Gz}$.
  This step involves the main new construction.
\medskip

\quad {\em Step 2.} To pass from an order-by-order solution $\{g_n\}$ to an ordinary solution $g\in G$, we use various approximation results.
For example, if all the equations involve power series we invoke first the Theorem of Popescu to ensure a formal solution over the completion
 (see Theorem \ref{Thm.Popescu}) and  then we use Artin-type approximations to ensure an ordinary solution. For $C^\infty$-equations we use
 an  approximation Theorem of Tougeron-type. 

\subsection{Main construction} \label{Sec.Intro.MainConstruction}

Fix a base ring $\k$ and a filtered $\k$-module $M=M_0\supset
M_1\supset M_2\cdots$  (usually not finitely generated over $\k$)
with the filtration topology, i.e., $\{M_i\}$  is a fundamental
system of neighbourhoods of $0 \in M$ (the $M_i$ are both open and
closed). Denote the group of $\k$-linear automorphisms of $M$ by
$GL_\k(M)$. The filtration of $M$ induces a natural descending
filtration of $GL_\k(M)$ by the normal subgroups $GL^{(i)}_\k(M)$,
consisting of elements $g$ that preserve the filtration and such
that $g$ and $g^{-1}$ are of the form $\one + \phi$  with $\phi$ an
element of \beq End_\k^{(i)}(M) := \{ \phi \in End_\k(M) \  |  \
\phi (M_j) \subseteq M_{j+i} \ \forall j \geq 0 \}, \eeq
with $End_\k^{(0)}(M)$ the endomorphisms of $M$ that respect the filtration.
\medskip

With the filtration topology we get as topological closure

 \bei
 \item for a submodule $\La\sseteq M$, $\overline{\La}=\capl_{i\ge1}(\La+M_i)$;
 \item for a  subgroup $G\sseteq GL_\k(M)$, $\overline{G}=\capl_{i\ge1}\left(G\cdot GL^{(i)}_\k(M)\right)$;
 \item for a  subgroup  $G\sseteq GL_\k(M)$ and an element $z\in M$, $\overline{Gz}=\capl_{i\ge1}\left( Gz+M_i\right)$;
 \item for a  submodule  $T\sseteq End_\k(M)$ and $z\in M$, $\overline{T(z)}=\capl_{i\ge1}\left( T(z)+M_i\right)$.
 \eei

Since the filtration topology is first-countable, the closure $\overline X$ of a subset $X \subset M$ consist
 of the points $x \in M$ for which there exists a sequence $x_n \in X$ converging to $x$.
 The same holds for subsets of  $ GL_\k(M)$.
If $\k$ is Noetherian and $M$ finitely generated (over $\k$) then any submodule of  $M$ is already closed.\\

Any subgroup $G$ of $GL_\k(M)$ carries the induced filtration $G^{(i)} := G \cap GL_\k^{(i)}(M)$.
 A special role here plays the following subgroup of $GL^{(1)}_\k(M)$:
\beq \label{eq3}
G^{(1)} := G \cap GL^{(1)}_\k(M).
\eeq

We call $G^{(1)}$ the  {\em (topologically) unipotent part of G}. Note that  if $\bigcap M_i =0$ (which will be the case in most of our applications) then $g \in G^{(1)}$ satisfies $\lim_{i\to \infty} (g-\one)^i =0$, but we use the term ``topologically unipotent''  also if $\bigcap M_i \neq 0$ (e.g. in the $C^\infty$-case  $\bigcap \cm^i$ contains the flat functions).

Since any element of
$G^{(1)}$ is of the form $\one + \phi$  with $\phi \in End_\k^{(1)}(M)$,
$G^{(1)}$ induces the identity on the ``linear" part $M_1/M_2$.\\

The key notion, which is introduced in this paper, is that of a ``pair of (weak) Lie-type'' of a subgroup $G \subset GL_\k(M)$. It is a pair
$$(T_{(G^{(1)},M)},G^{(1)}),$$
where $T_{(G^{(1)},M)}$ is a certain (not unique) submodule of $End_\k(M)$ that will be a substitute of the tangent space  to the orbit
 of the action of $G^{(1)}$ on $M$.
 The pair $(T_{(G^{(1)},M)},G^{(1)})$ is called of  {\em (weak) Lie type} if there exists
a (weak) substitution of the classical exponential and logarithmic
maps $T_{(G^{(1)},M)} \rightleftarrows G^{(1)}$. The choice of such substitutions is part of the data.
For a precise definition of (weak) Lie type and the even more general notions of pointwise (weak) Lie type we refer   to section  \ref{Sec.Pairs.of.Lie.type}.
  \medskip

If $\k$ contains the subring $\Q$ (e.g., if $\k$ is a field of
characteristic zero) then many groups admit the standard exponential and logarithmic
maps and they are trivially of Lie-type. In positive characteristic
however, the standard exponential map cannot be defined, but
nevertheless many groups are of weak Lie type. For example, for $R =
\k [[x]]$, the ring of power series over an arbitrary
field $\k$ in finitely many variables $x$,
the group of $\k$-algebra automorphisms $Aut_\k(R) \subset
GL_\k(M)$, acting on $M=R^n$ component-wise by coordinate change,
gives a pair of weak Lie-type (cf. Example  \ref{Ex.GL.times.GL.times.Aut.of.weak.Lie.type}). We mention that the group of
$R$-module automorphisms, $GL_R(M)\subset GL_\k(M)$, gives a pair of
Lie type for any ring $\k$, in any characteristic (see Example \ref{Ex.Lie.type.pair.GLM}).

\subsection{Main results}\label{Sec.Intro.Main.Results}
Let $\k$ be a ring, $M$ a filtered $\k$-module and $G\sseteq GL_\k(M)$ a subgroup
with induced filtration. We set
$$T_{(G^{(i)},M)} := T_{(G^{(1)},M)} \cap End_\k^{(i)}(M),$$
and prove the following general criterion for finite determinacy (under the slightly weaker assumptions of ``pointwise Lie type",
  see Theorems \ref{Thm.Finite.Determinacy.pointwise.Lie.type} and
\ref{Thm.Finite.Determinacy.weak.Lie.pairs}).
Fix some $z\in M$ and $i,N\ge 0$.
\beq\ber \label{Main.result.1}
\text{\em i. If $(T_{(G^{(1)},M)},G^{(1)})$ is a pair of Lie type and $M_{N+1}\sseteq \overline{T_{(G^{(i+1)},M)}(z)}$}
\\
\text{\em \quad\quad\quad\quad\quad\quad\quad\quad\quad then $\{z\}+M_{N+1}\sseteq\overline{G^{(i+1)}z}$.}
   \\\text{\em ii. If $(T_{(G^{(1)},M)},G^{(1)})$ is a pair of weak Lie type and $M_{N+k}\sseteq \overline{T_{(G^{(i+k)},M)}(z)}$ for any $k > 0$,}
   \\\text {\em   \quad\quad\quad\quad\quad\quad   then $\{z\}+M_{N+k}\sseteq\overline{G^{(i+k)}z}$ for any  $k>\max(0,N-2i-ord(z))$.}
      \eer\eeq
Here $ord(z) := sup \{ j \ | \ z \in M_j \}$ is the {\em order} of $z$.
This statement linearizes the determinacy question and reduces it to the level of the tangent space.
Part $i.$ implies that $z$ is $N$-determined and generalizes the classical results in characteristic zero as in \eqref{Eq.Mather.classical.theorem} (if we take $i=0$).
Part $ii.$ implies that $z$ is $(2N - ord(z))$-determined and
generalizes the known results in positive characteristic
as in \eqref{Eq.Greuel.Pham.thm.positive.char} below.

The assumption in (4.ii.) is however a condition
  for all $k > 0$, which is stronger than some known cases,
    but the conclusion is also stronger since it implies determinacy of $z$ + (terms of order $N+k$) already by the
     smaller group $G^{(k)} \subset G^{(1)}$ (for $i=0$).
      In many cases the assumption in (4.ii.) for all $k>0$ is already implied by the assumption for $k=1$ (cf. Corollary \ref{corollary of finite determinacy} and Example \ref{Ex.for.cond.50}).

The upper bound of {\em 4.ii} on the order of determinacy is  $(2N-2i- ord(z))$, which is weaker than $N$ of  {\em 4.i.}
 This $(2N-2i- ord(z))$   cannot be significantly improved, see section  \ref{Sec.Finite.Determinacy.Sharpness.of.bounds}.

\

Statement  \eqref{Main.result.1} may be rephrased as ``large tangent space implies determinacy".
We also prove the converse statement, ``determinacy implies large tangent space" (under slightly weaker assumptions, cf. Theorems
  \ref{Thm.Finite.Determinacy.pointwise.Lie.type} and
\ref{Thm.Finite.Determinacy.weak.Lie.pairs}), which reads:
\beq\ber \label{Main.result.2}
\text{\em i. Let $(T_{(G^{(1)},M)},G^{(1)})$ be a pair of  Lie type and
suppose $z\in M$}\\
\text{\em satisfies $\{z\}+M_{N+1}\sseteq \overline{G^{(1)}z}$. Then
$M_{N+1}\sseteq\overline{T_{(G^{(1)},M)}(z)}$.}
\\\text{\em ii.  Let $(T_{(G^{(1)},M)},G^{(1)})$ be a pair of weak Lie type and suppose $z\in M$ }
 \\\text{\em satisfies $\{z\}+M_{N+k}\sseteq \overline{G^{(k)}z}$ for any $k\ge1$.}\\
\text{\em Then $M_{N+k}\sseteq\overline{T_{(G^{(k)},M)}(z)}$ holds for any $k>N-ord(z)$.}
\eer\eeq

\medskip

In section  \ref{Sec.Applications.Examples} we couple these statements with the approximation theorems of Popescu, Artin
 and Tougeron to get from an order-by-order solution an algebraic, resp. analytic, resp. $C^\infty$ solution.
   The proposed generality allows to recover finite determinacy statements for many particular scenarios, e.g., for germs of functions, of maps on smooth and non-smooth spaces and of matrices.  In particular,

\bee[ i.]
\item when $\k$ is a field of characteristic zero this recovers numerous classical results e.g., by Mather, Gaffney,
Bruce - Du Plessis - Wall, Damon, and many others;
\item  when $\k$ is a field  of positive characteristic this gives other known results, e.g., those of
Boubakri - Greuel - Markwig and Greuel - Pham;
\item the notion of weak Lie-type might be potentially useful not only in prime characteristic, but also in mixed characteristic,
 as we do not impose any kind of restriction on the base $\k$;
\item beyond this we get new results, e.g on relative determinacy results for non-isolated singularities, section \ref{determinacy for maps that preserve a subscheme} and \ref{determinacy for non-isolated singularities}.
\eee

\subsection{Remarks}
In order to put our results into perspective, we finish this introduction by giving references to previous results  (by far not complete).

\subsubsection{}
Investigations on determinacy had been classically  restricted to the real resp. complex case,
with $M$ being formal or analytic power series, or to germs of $C^\infty$-maps. In order to apply methods and
 results from Lie groups or algebraic groups, the setting was immediately reduced to a finite dimensional parameter space (either a finite
  jet space or the parameter space of a semi-universal deformation) by assuming some kind of ``isolated singularity''.
The proofs used essentially complex or real analysis, integration of vector fields, and topology.
See \cite{Wall-1981}, \cite{Damon84} for a short introduction and \cite[section 2.7.2]{B.K.motor} for some more recent history.
It was observed in \cite{Bruce.du-Plessis.Wall} that in fact the essential ingredient for a statement like \eqref{Eq.Mather.classical.theorem}
is the unipotency of the group action. In \cite{B.K.motor}
this idea was used to extend  \eqref{Eq.Mather.classical.theorem}
to Henselian rings over a field of characteristic zero and to filtered groups possessing a (formal)
exponential and logarithmic map, or at least an  ``order-by-order" version of these maps.
These exponential and logarithmic maps were the basis of the construction, thus the method seemed to be inapplicable to the case when the
base ring does not contain the rational numbers.
Furthermore, to define the tangent space one had to restrict to some particular class of group actions, though broad enough to include most
of the known scenarios.

\subsubsection{}
Another direction of generalization was to  positive
characteristic. This study was initiated (to the best of our knowledge) in \cite{GrKr.90} and then continued in \cite{Bou09}, \cite{Boubakri.Gre.Mark},
\cite{Greuel-Pham.mather-yau.char.positive}, \cite{Pham}.
In \cite{Greuel-Pham.2018} the authors considered the case of matrices over the ring $R=\k[[x]]$,  $\k$ an arbitrary field of any characteristic, with $M=Mat_{m\times n}(R)$, and
 the group $G=GL(m,R)\times GL(n,R)\rtimes Aut_\k(R)$. The
  proved result was  ($\cm=\bl x_1,...,x_p\br$ the maximal ideal,
  $A \in Mat_{m\times n}(R)$):
\beq\label{Eq.Greuel.Pham.thm.positive.char}
\ber
\text{\em If \ $\cm\cdot \widetilde{T}_{G}(A)\supseteq \cm^{N+1} \cdot Mat_{m\times n}(R)$ then $GA\supseteq\{A\}+ \cm^{2N+1-ord(A)}\cdot Mat_{m\times n}(R)$.}
\\\text{\em In particular, $A$ is then $(2N-ord(A))$ determined.}
\eer\eeq
Here $\widetilde{T}_{G}(A)$ is the {\em tangent image}, i.e., the image of the tangent map of the orbit map $G \to Mat_{m\times n}(R)$, $g \mapsto gA$.

In characteristic zero  $\widetilde{T}_{G}(A)$ coincides with the tangent space $T_{G}(A)$ to the orbit $GA$ at $A$, but in
positive characteristic $\widetilde{T}_{G}(A) $ differs in general from $T_{G}(A)$. In our general framework the module $T_{(G^{(1)},M)}(z)$
of a pair of weak Lie-type is a generalization of the tangent image.
For $m=1=n$ we get the contact determinacy of function germs, recovering \cite[Theorem 3]{Boubakri.Gre.Mark}.

Note that these bounds coincide with those of equation \eqref{Main.result.1}.

\subsubsection{}
The approximation step ({\em Step 2 } in section  \ref{Sec.Intro.Goals}) is more or less standard,
 we repeat it briefly in section  \ref{Sec.Approximation.Theorems}.  The case $C^\infty$ is more involved, we treat it in \ref{Sec.Approximation.Cinfty.case}.

\medskip

In some cases no approximation theorems are possible. For example, analytic questions of dynamical systems or differential equations
are notoriously difficult when compared to the formal ones. Even in such cases, establishing that ``two objects are order-by-order equivalent"
is a significant result.

\subsubsection{}
In this paper we address only the topologically unipotent part  $G^{(1)}$ of the group $G$. However, this has no significant impact on the finite determinacy,
as for most ``reasonable" groups over a local ring $R$ we have $\cm\cdot T_{G}\sseteq T_{(G^{(1)},M)} \sseteq T_{G}$ for the maximal ideal $\cm\sset R$.
Accordingly, the orders of determinacy under $G$ and $G^{(1)}$ differ at most by one in the case of a Lie type
pair, and by two in the case of a weak Lie type pair.

\subsubsection{}
We do not consider the group of left-right or $\cA$-equivalences. This group is not a subgroup of $GL_R(R^n)\rtimes
Aut_\k(R)$, as the $\cA$-action is not $\k$-linear. Thus our method
does not apply directly. Certain modifications are needed to establish
a (weak) Lie-pair structure for this group and to obtain the
``determinacy vs. infinitesimal determinacy" statements. We hope to
report on this soon.

 \subsubsection{}
In many cases a result of type   \eqref{Eq.Mather.classical.theorem} is not yet a complete solution. The tangent module
can be rather complicated, and to check the condition  $\cm\cdot T_{(Gf,f)}\supseteq\cm^{N+1}\cdot M$ in particular cases
can be a difficult task (although, when $R$ is the ring of power series, standard basis methods
provide effective algorithms, cf. \cite{Greuel-Pham.2017b} or
\cite{Abzal-Kanwal-Pfister.2017}). For example,
for matrices over local rings and various groups acting on them, one gets non-trivial questions on the annihilators of quotient modules,
 see \cite{B.K.fin.det.1}, \cite{B.K.fin.det.2},  \cite{Kerner}.

\subsection{Contents of the paper}
\bei
\item
Section \ref{Sec.Preparations} is preparatory, we review the relevant facts about filtered rings,  the Implicit Function Theorem (with the ``unit main part")
 filtered modules, and  the associated filtration on $GL_{\k} (M)$ and its subgroups.\\
In section \ref{Sec.Approximation.Theorems} we recall some relevant facts on Artin approximation and their extensions to non-polynomial equations.
\medskip

 \item
 Sections \ref{Sec.Pairs.of.Lie.type} and \ref{Sec.Finite.Determinacy.Criteria} form the core of the paper.
 In section  \ref{Sec.Pairs.of.Lie.type} we introduce  the pairs of (pointwise/weak) Lie type. As was briefly mentioned in section \ref{Sec.Intro.MainConstruction},
  these are groups together with a substitute for the tangent space at the orbit of the action. The orbit is locally approximated by its tangent
space and there are some substitutions for the classical exponential/logarithmic maps.\\
We show that the class of such pairs is rich enough. It contains the main interesting subgroups of $GL_\k(M)$, in particular
 the groups $Aut_\k(R)$ (under certain assumptions on $R$), $GL_R(M)$, the groups of contact equivalences for matrices, and the (semi-)direct products of these groups (cf. Theorem \ref{Thm.Aut(R).pointwise.Lie}).
\medskip

\item
Section  \ref{Sec.Finite.Determinacy.Criteria} contains the main results of this
paper (as indicated in section \ref{Sec.Intro.Main.Results}), the Finite Determinacy Theorems \ref{Thm.Finite.Determinacy.pointwise.Lie.type} and
\ref{Thm.Finite.Determinacy.weak.Lie.pairs}.
 These results establish the ``determinacy versus infinitesimal determinacy" criteria in a characteristic free way.\\
 The determinacy bounds for pairs of weak Lie type are weaker than those for pairs of Lie type. We show in section \ref{Sec.Finite.Determinacy.Sharpness.of.bounds}
  that these weaker bounds are often sharp.\\
In section \ref{Sec.Finite.Determinacy.vs.Infinitesimal.Stability} we translate the finite determinacy
 into the infinitesimal stability in the traditional way: an element $z\in M$ is finitely determined iff its fibres are infinitesimally stable
  on the punctured neighborhood $Spec(R)^\times$.
  \\
  The ring $C^\infty(\R^p,0)$ does not have the Artin approximation property, thus we cannot pass from $\overline{Gz}$ to $Gz$ by using
   the general/standard theory of section \ref{Sec.Approximation.Theorems}. However, we prove the relevant approximation statement by using the
    Lie type notions, section \ref{Sec.Approximation.Cinfty.case}.
\medskip

\item Finally, in section  \ref{Sec.Applications.Examples}
we couple the Finite Determinacy Theorems
\ref{Thm.Finite.Determinacy.pointwise.Lie.type} and
\ref{Thm.Finite.Determinacy.weak.Lie.pairs}
with various approximation results and apply these to several particular scenarios.
 We recover and extend numerous classical results and obtain some new results (in arbitrary characteristic), e.g.
  the right determinacy of germs of functions, the right
 indeterminacy of germs of maps, the contact determinacy of germs of maps,
 determinacy of maps relative to a space-germ, relative determinacy of non-isolated singularities, relative algebraization of power series, determinacy of matrices.
\eei

\subsection{Acknowledgements} We thank D. Popescu for highly valuable explanations about the generality of Theorem \ref{Thm.Popescu},
and L. Narvaez Macarro
and M. Schulze for their explanations about logarithmic derivations.
We also thank G. Belitski and I. Tyomkin for valuable advices. Finally we thank the referee for the numerous helpful comments.

\section{Preparations}\label{Sec.Preparations}

\subsection{Notations and assumptions} \label{Sec.Notations}
In this paper all rings are supposed to be associative, commutative and unital,
with subrings having the same identity element as the ambient ring.
  We fix a base ring $\k$, not necessarily a field,
of any characteristic.
 Consider a  $\k$-module filtered by submodules, $M=M_0\supset M_1\supset
M_2\cdots$, usually not finitely generated over $\k$.

We consider also a (commutative, associative) $\k$-algebra $R$, filtered by a chain of ideals,
 $R=I_0\supset I_1\supset\cdots$ with $I_j\cdot I_k\sseteq I_{j+k}$,  e.g., $I_j =I^j$  for some ideal $I$.

If $R$ is local  (i.e., a ring with the unique maximal ideal $\cm$) then the classically
 considered filtration is by  the powers of  $\cm$.

If $M$ is also an $R$-module, not just a $\k$-module, then the filtrations of $R$, $M$ are supposed to be compatible, i.e.,  $I_j M_i\sseteq M_{i+j}$ for any $i,j$.

When a group $G$ acts on $M$ one fixes some action
$G\circlearrowright R$, possibly trivial, and then the action $G\circlearrowright M$
is assumed to be $R$-multiplicative, i.e., $g(f\cdot z)=g(f)\cdot g(z)$  for any
$g \in G$, $f\in R$, $z\in M$.

On any filtered object we use the filtration topology.
 Thus $\{M_i\}$ is a fundamental system of neighbourhoods of $0 \in M$,
 the $M_i$ are both open and closed. (See e.g. \cite[Chapter III, section 2, no.1, Corollary to Proposition 4]{General.Topology}.)
The completion of $M$ \wrt this topology is $\widehat M:=\liml_{\leftarrow}M/M_i.$
 The orbit closure of an element $z$ is $\overline{Gz}=\cap_i(Gz+M_i)$.

\subsubsection{Typical rings}\label{Sec.Examples.of.rings}

Let $R$ be a $\k$-algebra of any characteristic (e.g., $\k$ a field, the classically considered case).
We use the multi-variable notation (always finitely many variables), $x=(x_1,\dots,x_p)$.
Typical examples for $R$ are:
\bei
 \item  the formal power series, $\k[[x]]$;
 \item  the algebraic power series, $\k\langle x\rangle$;
 \item  the convergent power series,   $\k\{x\}$, when $\k$ is a complete normed ring;
 \item  the germs of $C^\infty$-functions, $C^\infty(\R^p,0)$;
 \item the quotients of these rings by some ideal, \ $ \k[[x]]/J$, $\k\{x\}/J$, $\k\langle x \rangle/J$,  $C^\infty(\R^p,0)/J.$
 \eei

More generally, we often assume the condition:  $R\sseteq \k[[x]]/J$ and $R$ contains the images of $\{x_i\}\sset \k[[x]]$.
In this case the completion map $R\hookrightarrow\hR$  is injective and $\hR$ is of Cohen-type.
 This is satisfied by many rings, not necessarily over a field, see e.g., page 214 in \cite{Lu}.

\subsubsection{Implicit function Theorem with unit linear part, $\bf IFT_\one$}\label{Sec.IFT.one}
We often need to solve an equation $z+w=gz$ for the unknown $g\in G$.
 In many cases this condition can be stated as  a system of implicit function equations.
  For example, let $R\sseteq\k[[x]]$, $M=R$, $G=Aut_\k(R)$, then given $f\in R$ and its perturbation $h\in R$ we want to resolve $f(y)=f(x)+h(x)$, $y=gx$ for the unknown $g\in G$.

Moreover (and this is one of the main results of section 3), in certain cases this implicit function equation can be transformed to a simpler form (with  $x \in R^n$, and $y =gx$  an $n$-tuple of unknowns):
\beq\label{Eq.IFT}
y+h(y,x)=x,  \quad h(y,x)\in \langle y^2,x\cdot y \rangle R^n:=\langle\{y_i \cdot y_j, x_i \cdot y_j\}\rangle R^n \sset R[[y]]^{n}.
\eeq
Here the second condition means that $h$ is ``of higher order".
Usually implicit function equations are studied with the assumption that $R$ contains  a field. This assumption is not needed for   equation \eqref{Eq.IFT}.

\bed (cf. \cite[section 3.4]{B.K.IFT})
We say that $IFT_\one$   holds over $R$
 if the equation  $y+h(y,x)=x$ has a solution  $y(x)\in R^{n}$ for any $h(y,x)$
satisfying $h(y,x)\in \langle y^2,x\cdot y \rangle R^{n}$,
 and moreover $y(x)-x\in \langle\{x_i \cdot x_j\}\rangle R^n\sset R^n$.
\eed

\bex\label{Ex.rings.with.IFT1}
Examples of rings for which $IFT_\one$ holds include:
\bee[i.]
\item $\k[[x]]/J, \ \k\langle x\rangle/J$, and $\k\{x\}/J$ for $\k$ a normed field, complete w.r.t. its norm.
\medskip
\item More generally, let $\k$ be
either a field or a discrete valuation ring, and $\{W_n\}_{n\geq 1}$ be a Weierstrass system of rings
over $\k$ \cite[page 2]{Denef-Lipshitz}. Then, for any $n\geq 1,$ the ring  $W_n/J$ satisfies $IFT_\one$ \cite[page 4]{Denef-Lipshitz}.
\medskip
\item $C^\infty(\R^p,0)/J$. For $J=(0)$ see \cite[page 79]{Raynaud}. 
Otherwise one can assume $J\sseteq (x)^2$. Then one lifts the equations to  $C^\infty(\R^p,0)$, resolves there and sends the solutions back to $C^\infty(\R^p,0)/J$.
\medskip
    \eee
\eex
\smallskip

\subsection{Filtered modules  and group actions}\label{Sec.Filtered.Modules.Group.Actions}
\subsubsection{Induced filtrations on $End_\k(M)$ and $GL_\k(M)$} \label{Sec.Preparations.Induced.Filtration} Let $\k$ be a commutative ring.
Fix a filtered $\k$-module,  $M=M_0\supset M_1\supset \cdots$ and consider the set of all $\k$-linear
 endomorphisms, $End_\k(M)$. The filtration of $M$
induces a filtration of $End_\k(M)$ by $\k$-submodules
 \begin{eqnarray}
End_\k(M)\supseteq   End^{(0)}_\k(M)\supset   End^{(1)}_\k(M)\supset\cdots,
\text{where} \quad \quad \\ End^{(i)}_\k(M):=\{\phi\in End_\k(M)|\ \phi(M_j)\sseteq M_{j+i},\ \forall j\ge0\}. \nonumber
\end{eqnarray}
Here $End^{(0)}_\k(M)$ is the module of endomorphisms that preserve the filtration, while we call
$End^{(1)}_\k(M)$ the module of  {\em topologically nilpotent endomorphisms} (see the remark after equation  (\ref{eq3})).

\medskip

The {\em order} of an element $z\in M$ is defined as
$ord(z):=sup\{j|\ z\in M_j\}$ In particular, $ord(z) =\infty$
if $z\in \cap_{j=1}^\infty M_j$.
Similarly, for any morphism $\phi\in End_\k(M)$ the order is
$ord(\phi):=sup\{j|\ \phi\in End^{(j)}_\k(M)\}$.

\medskip

Denote the group of all $\k$-linear automorphisms
of $M$ by $GL_\k(M)$.
Define the subgroup of automorphisms that preserve the filtration,
\beq
GL^{(0)}_\k(M):=\{g|\   g,g^{-1}\circlearrowright M_i,\ \forall i\}\sseteq GL_\k(M).
\eeq
It is filtered by its subgroups,
\beq\label{Eq.group.filtration}
GL^{(i)}_\k(M):=GL^{(0)}_\k(M)\cap\left(\{\one\}+End^{(i)}_\k(M)\right), \ i\geq 1,
\eeq
i.e., the subgroups of $GL^{(0)}_\k(M)$, such that the elements and their inverses are of the form $\one + \phi$ with $\phi \in End^{(i)}_\k(M)$.
 $GL^{(1)}_\k(M)$ is the subgroup of {\em topologically unipotent automorphisms}. One readily checks that $GL^{(i)}_\k(M)$
  is indeed a subgroup of $GL^{(0)}_\k(M)$.  For example, if $g\in GL^{(i)}_\k(M)$ then $g^{-1}\in GL^{(i)}_\k(M)$: since
 $g|_{{M_j}/{M_{i+j}}}=\one|_{{M_j}/{M_{i+j}}}$ for all $j$, $g^{-1}|_{{M_j}/{M_{i+j}}}=\one|_{{M_j}/{M_{i+j}}}$ and thus
 $g^{-1}-\one\in End^{(i)}_\k(M)$.

As the simplest case take
$M = \mathcal O_{(\C^n,o)}=\C\{x\}$, $M_i = \langle x\rangle^i$ and $G=\cR$, the right group  acting by $\phi \cdot f = f \circ \phi^{-1}$ for
$\phi \circlearrowright (\C^n,0)$. Then $GL^{(i)}_\k(M)$ is the subgroup $\cR^{(i)}$ leaving the $i + 1$-jet of every $f \in \mathcal O_{(\C^n,0)}$ unchanged.
\medskip

The  group $GL^{(0)}_\k(M)$ will be always taken as the ambient group. For any subgroup $G \subseteq GL^{(0)}_\k(M)$ we get the induced  filtration,
\beq
G^{(i)}:=G\cap GL^{(i)}_\k(M),\quad i=0,1,\dots \, .
\eeq
These are always normal subgroups, $G^{(0)}\vartriangleright G^{(1)}\vartriangleright\cdots$, as their action on the quotient $M/M_{i+1}$ is trivial.
\medskip

\subsubsection{The structure of $GL^{(i)}_\k(M)$ and $GL^{(i)}_R(M)$}
\label{Sec.Examples.of.groups}

 In (\ref{Eq.group.filtration}) we define $GL^{(i)}_\k(M)$ as a subset of $\{\one\}+End^{(i)}_\k(M)$, where $\k$ can be any commutative ring.
 We give now examples of group actions where both sets coincide.

 \bee[i.]

\item Suppose $M$ is complete \wrt the filtration $\{M_i\}$. Then we have for all $ i\ge1$:
\beq\label{Eq.Glk(M)}
 GL^{(i)}_\k(M)=\{\one\}+End^{(i)}_\k(M) =\{ \one+\phi, \ \phi\in End^{(i)}_\k(M)\}.
\eeq
The inclusion $GL^{(i)}_\k(M)\sseteq\{\one\}+End^{(i)}_\k(M)$ follows just from the definition, equation \eqref{Eq.group.filtration}.
To show $GL^{(i)}_\k(M)\supseteq\cdots$  it is enough to check  that any endomorphism of the form $\one+\phi$,   for $\phi\in End^{(1)}_\k(M)$,  is invertible.
Indeed, its inverse, $\sum_{j=0}^\infty (-\phi)^j$, is a well defined $\k$-linear operator on $M$ (by completeness of $M$  and the topological nilpotence of $\phi$).

If $M$ is not complete  then statement \eqref{Eq.Glk(M)} does not
necessarily hold, regardless of how nice is $\k$. For example, suppose $\k$ is a
field and $M=\k[x]$. Consider the operator $(1+x)\in End_\k(M)$, acting by $p(x)\mapsto (1+x)p(x)$. Though it is topologically unipotent, it is not invertible.
\medskip

\item   Now consider a ring extension $\k \subset R$ such that $M$ is also a filtered module over the larger ring $R$.
Suppose that $GL^{(1)}_\k(M)=\{\one\}+End^{(1)}_\k(M)$ holds for the $\k$-module  $M$, then we have also
$GL^{(1)}_R(M)=\{\one\}+End^{(1)}_R(M)$.
\smallskip

\item
Take a filtered ring $R$ and let $M$ be a finitely generated $R$-module with a filtration $M=M_0\supset M_1\supset\cdots$.
Let $J(R)$ be the Jacobson radical of $R$, if $R$ is local then $J(R)$ is the maximal ideal. Assume $M_j\sseteq J(R)\cdot M$, for $j\gg1$.
Then the group of $R$-linear automorphisms, $GL_R(M)$, is related to the module of $R$-linear endomorphisms $End_R(M)$:
\beq\label{Eq.GlR(M)}
\forall\ i\ge1: \quad GL^{(i)}_R(M)=\{ \one\}+End^{(i)}_R(M).
\eeq
As before, the inclusion ``$\sseteq$'' follows from the definition.
For ``$\supseteq$" it is enough to check that $\one-\phi$    is invertible  for $\phi\in End^{(1)}_R(M)$. This can be proved as follows.
By our assumption there is an integer $d$ with $\phi^d(M)\sseteq M_d\sseteq J(R)\cdot M$.
Now notice:
 \beq
\Big(\one-\phi\Big)\Big(\one+\phi+\cdots+\phi^{d-1}\Big)=\one-\phi^{d}\equiv \one\ \ mod\ \ J(R)\cdot End_R(M).
\eeq
 Setting $\psi:=\one- \phi^{d}$ we get $M\subset \psi(M) + J(R)\cdot M$.
By Nakayama's Lemma $M=\psi(M)$, i.e. $\psi$ is surjective and hence bijective (\cite[Theorem 2.4]{Matsumura}).
%
 Therefore   $\one-\phi$ is invertible, with $(\one-\phi)^{-1}=\Big(\one+\phi+\cdots+ \phi^{d-1}\Big)(\one-\phi^{d})^{-1}$ .
\medskip
\eee

\subsubsection{Actions involving ring automorphisms}\label{Sec.Examples.Actions.with.Aut}
Let $\k$ be a ring and $R$ be a (not necessarily Noetherian) filtered $\k$-algebra, e.g., one of the rings in section \ref{Sec.Examples.of.rings}.
 \bee[i.]

\item
Considering $R$ as a $\k$-module we get the group $GL_\k(R)$ and its filtration.
 Consider the subgroup $Aut_\k(R)\sset GL_\k(R)$ of  $\k$-linear
  automorphisms of this ring. This group is naturally filtered,
\beq
Aut^{(i)}_\k(R):=Aut_\k(R) \cap GL^{(i)}_\k(R),
\eeq
and every $g\in Aut^{(i)}_\k(R)$, $i\ge1$, is of the form  $\one+\phi$
 for some $\phi\in End^{(i)}_\k(R)$, by equation \eqref{Eq.group.filtration}.
The multiplicativity property of $g\in Aut_\k(R)$ imposes the ``almost Leibniz rule" for $\phi$,
\beq\label{Eq.16}
\phi(f_1f_2)=\phi(f_1)f_2+f_1\phi(f_2)+\phi(f_1)\phi(f_2).
\eeq
We claim: if $IFT_\one$  holds over $R$ (see  \S\ref{Sec.IFT.one}) then \eqref{Eq.16} is the only condition on $\phi$. Namely,
\beq
Aut^{(i)}_\k(R)=
\Big\{\one+\phi \ \big|\ \phi\in End^{(i)}_\k(R),\  \phi(f_1f_2)=
\phi(f_1)f_2+f_1\phi(f_2)+\phi(f_1)\phi(f_2),\ \forall\ f_1,f_2\in R\Big\}.
\eeq
Indeed, $\one+\phi$ is $\k$-linear, while the almost Leibniz rule ensures the multiplicativity.
And any topologically unipotent endomorphism is invertible, as the equation $x=y+\phi(y)$ is solvable over $R$, by $IFT_\one$.
\medskip

\item Let $R$ be a $\k$-subalgebra of $\k[[x]]/J$ and assume that $R$ contains the images of all  $x_i\in \k[[x]]$.
 Consider elements of $R$ as power series, $f(x)\in R$.
 Then any topologically unipotent  automorphism  $\phi\in Aut_\k(R)$  acts as a coordinate change.
 Namely, fix its action on the generators $x_i\in \k[[x]]$,
\beq
x_i \mapsto \phi(x_i)=x_i+h_i(x),\quad
h_i(x)\in\cm^2\sset \k[[x]],
\eeq
 and consider $\{x_i\}$ as coordinates on $Spec(R)$. Then   $\phi(f(x)) = f(\phi(x)) $   \text{and}  $ \phi(J)= J$.

(See e.g., \cite[2.11]{Abhyankar}, \cite[Lemma I.1.23]{Gr.Lo.Sh} or \cite[Lemma 3.1]{B.K.motor}.
In that lemma $\k$ was assumed  to be a field of zero characteristic, but the proof is characteristic-free.)

Because of this, the group $Aut_\k(R)$ (resp. $Aut^{(1)}_\k(R)$) is also called the   {\em right group} (resp. {\em unipotent right group}) or the group of right equivalences of function germs.

For more general rings there exist ``non-geometric" automorphisms, not arriving from coordinate changes as above.
Yet, if the subring $R\sseteq C^\infty(\R^p,0)/J$ admits
the Taylor expansion up to any order, then the coordinate changes are dense inside all the automorphisms.
Therefore, for  $R \sseteq C^\infty(\R^p,0)/J$,     (following \cite{Tougeron1968}) we  denote by $Aut_\k(R)$ the subgroup of all automorphisms $\phi$ of $R$
 arriving from coordinate changes, i.e., $\phi (f(x)) = f(\phi(x))$, with $\phi =(\phi_1,...,\phi_n)$ a $C^\infty$-coordinate change of $(\R^p,0)$
 satisfying $\phi(J)=J$. As before,  $Aut^{(1)}_\k(R)$ denotes the subgroup of coordinate changes satisfying $\phi(x) = x+h(x)$ with $h(x)\in\cm^2$.

\medskip

\item
Let $M$ be a (not necessarily free) $R$-module and $G$ a subgroup of  $GL_R(M)\rtimes Aut_\k(R)$.
 We assume  some prescribed {\em $R$-multiplicative action} $G\circlearrowright M$. Namely, for $g = (W,\phi) \in G\subseteq GL_R(M)\rtimes Aut_\k(R)$ one has
 \beq
 g (a\cdot z) = \phi(a)\cdot g(z), \quad a\in R, \ z \in M,
 \eeq
where  $\phi(a)$ is the standard action $Aut_\k(R)\circlearrowright R$.

For example, for a free module we fix a basis, i.e., identify $M=R^n$. Then we have the action  $GL_R(R^n)\rtimes Aut_\k(R)\circlearrowright R^n$,
 where $Aut_\k(R)$ acts component-wise. This recovers the classical right and contact equivalence of maps.

If $M$ is not a free module then it is a quotient of a free module and the action must of course respect the relations. In particular,
the possible actions $Aut_k(R)\ni \phi\circlearrowright M$ are
restricted. For example, the support and the annihilator of $M$ must be
preserved by $\phi$.
 Indeed, if $f\cdot M=0$ then
 $ \phi(f)\cdot M = \phi(f)\cdot \phi(M) = \phi(f \cdot M) =0$.

 \item
We consider the following group actions, with $M= \Mat$ the free $R$-module of rectangular matrices $A=\{a_{ij}\}$:
\bei
\item $GL(m,R)\circlearrowright\Mat$, by $(U,A)\to UA$;
\item  $GL(n,R)\circlearrowright\Mat$, by $(V,A)\to AV^{-1}$;
\item $Aut_\k(R)\circlearrowright\Mat$ by $(\phi,A)\to \phi(A)=\{\phi(a_{ij})\}$;
\item products of these groups, e.g. $GL(m,R)\times GL(n,R)\rtimes Aut_\k(R)$;
\item (for $m=n$) the  conjugation $A\to UAU^{-1}$, the congruence $A\to UAU^t$, etc.
\eei
 The filtration $R=I_0\supset I_1\supset\cdots$ induces the filtration $\{Mat_{m\times n}(I_j)\}$, accordingly one gets the subgroups $G^{(j)}\le G$.

We remark: if the action $G\circlearrowright\Mat$ is $R$-linear and preserves the subset of degenerate matrices then $G$ is contained in the group of left-right multiplications,
   $G_{lr}:=GL(m,R)\times GL(n,R)$.
  (See \cite[\S3.6]{B.K.fin.det.1} for the precise statement.)

\eee

\subsection{The condition $z+w\in Gz$ as implicit function equation}\label{Sec.Condition.z+w=gz.As.implicit.function.eq}
\noindent
For the initial determinacy question (whether $Gz\supseteq\{z\}+M_{N+1}$) we should resolve the equation (for the unknown  $g\in G$):
\beq
z+w=gz,\quad \quad \text{ for any given } w\in M_{N+1}.
\eeq

Let us explain in detail how these conditions are presentable  as a finite system of implicit function equations,  $z+w-gz =F(g) = 0$,
for some power series equations $F( y)\in R[[ y]]^{s}$, if $M$ is finitely presented as $R$-module.

Fix any finite presentation   $R^q\stackrel{A}{\to}R^p\stackrel{\al}{\to}M\to0$ and
consider the $R$-multiplicative action of some group  element $g$ on $M$,
\beq
GL_R(M)\rtimes Aut_\k(R)\ni g=(W,\phi)\circlearrowright M.
\eeq
The easiest way to see that the action of $g=(W,\phi)$ lifts to the presentation (and is determined by it) is
by describing it by an $R$-linear map as follows.
For $\phi \in  Aut_\k(R)$ denote by $M_\phi$ the
$R$-module, which is equal to $M$ as $\k$-module but with the new scalar multiplication
     $$a*z := a*_{\phi}z := \phi(a)\cdot z, \ a\in R, \ z \in M,$$
 also called an $R$-linear homomorphism over $\phi$. Then the map
$$W_{\phi}: M \to M_\phi \text{ with }  W_{\phi}(z) := (W,\phi) \cdot z$$
 is $R$-linear since it is $\k$-linear and satisfies
   $W_{\phi}(az) = a* W_{\phi}(z).$
The filtration $R=I_0\supset I_1\supset\cdots$ induces the filtrations $I_iM$ and  $I_i*M_\phi$ and $W_\phi$ is a map of filtered $R$-modules,  assuming that $\phi$ preserves $\{I_j\}$.

If  $A = [a_{i,j}]$,
then
$R^q  \xrightarrow{\text{$A_\phi$ }}R^p \xrightarrow{\text{$\alpha_\phi$ }} M_\phi \to 0$ is a presentation of $M_\phi$,
with $A_\phi :=\phi^{-1}(A)= [\phi^{-1} (a_{i,j})]$ and
$\alpha_\phi(\sum_i a_i  e_i)=\sum_i a_i * \alpha(e_i)$, $\{e_i\}$
the canonical basis of $R^p$.
As is well known, the linear map $W_\phi : M \to M_\phi$  lifts to a morphism of the presentations of $M$ and $M_\phi$
 (unique up to homotopy of complexes)
 and we get a commutative diagram
\beq \label{ComDiag}
\begin{xy}
  \xymatrix{
R^q \ar[r]_A \ar[d]_{V}   &   R^p \ar[r]_\alpha \ar[d]_{U}   &  M   \ar[r]  \ar[d]_{W_{\phi}} & 0\\
R^q \ar[r]^{A_\phi} &  R^p \ar[r]^{\alpha_\phi} &  \   M_\phi \ar[r]   & 0
}
\end{xy}
\eeq
of filtered modules ($R^r$ being filtered by $I_iR^r$). Since $W_\phi$ is an isomorphism, it follows that the linear maps $U= [u_{i,j}]$ and $V= [v_{i,j}]$ (depending on the action of $g=(W,\phi) $) are isomorphisms.
 To see this we may
assume that $R$ is local (just localize diagram (\ref{ComDiag}) w.r.t. the maximal ideals of $R$).
 For local rings the result follows from (the proof of)  \cite[Theorem 20.2]{Eisenbud}), since in our situation the free modules of the presentation have the same rank.
The action on $M \cong R^p/im(A)$ is hence induced by an action on $R^p$ preserving $im(A)$.
Conversely, for a given $\phi$ the left-hand square of (\ref{ComDiag}) determines
a commutative diagram as above with an isomorphism
$W_\phi :M:= coker(A) \to coker(A_\phi) =: M_\phi$ and hence a $\k$-linear action of $(W,\phi)$ on $M$ (note that $M$ and $M_\phi$ coincide as $\k$-modules).

\

We will mainly work with the action of $(U,V,\phi) \in GL_R(R^p)\rtimes GL_R(R^q)\rtimes Aut_\k(R)$ on a presentation of $M$.
  If $z= \sum z_i\alpha(e_i)$ and  $w= \sum w_i\alpha(e_i)$ are  elements in $M$
set $\tilde z := \sum z_ie_i$ and $\tilde w := \sum w_ie_i$
for the $\alpha$-preimages of $z$ and $w$.
Then  $\phi^{-1}(\tilde z) = \sum \phi^{-1}(z_i)e_i$ and $\phi^{-1}(\tilde w) = \sum \phi^{-1}(w_i)e_i $ are $\alpha_\phi$-preimages  in $R^p$ of $z$ and $w$.
By diagram (\ref{ComDiag}) the condition
$z+w = gz$ reads then:  $\phi^{-1}(\tilde z) +\phi^{-1}(\tilde w) = U\tilde z+\phi^{-1}(A) v $ for some $ v \in R^q$ and $UA=\phi^{-1}(A)V$. After applying  $\phi$ and renaming we get the system of equations
\beq \label{Cond}
\tilde z +\tilde w = U\cdot \phi(\tilde z)+ A\cdot v, \ U\cdot \phi(A)=A\cdot V,
\eeq
with  $\tilde z, \tilde w, A$ given and (the entries of) $U,V,\phi, v$ as unknowns.

%
%

 To this explicit form \eqref{Cond}  we add the condition  ``$g(J) = J$"  (if $R$ is a quotient ring mod $J$)
 and the condition  ``$g \in G$",  if $G$ is a subgroup of  $ GL_R(R^p)\rtimes GL_R(R^q)\rtimes Aut_\k(R)$, as follows.
\bei
\item[i.]
Assume $R\sseteq S/J$, where $S=\k[[x]]$ or $C^\infty(\R^n,0)$, thus
$\phi\in Aut_\k(R)$ is induced by a coordinate change $\phi \in Aut_\k(S)$
  that satisfies $\phi(J) = J$. Fix some generators $\{f_j\}$ of $J$, then the condition $\phi(J) = J$ becomes the system of (formal)
   power series  (or $C^\infty$-germs) equations
    $f_j(\phi(x)) = \sum_i t_{j,i}f_i$, in the unknowns  $\phi(x)=(\phi_1,...,\phi_n)$ and  $\{t_{j,i}\}$.
\item[ii.]
Assume that $G\sseteq GL_R(M)\rtimes Aut_\k(R)$ is defined by (formal) power series. That is, for a presentation of $M$ (as above) and the lift of $G$   to a subgroup of $ GL_R(R^p)\rtimes GL_R(R^q)\rtimes Aut_\k(R)$,
 the conditions on
   $(U,V,\phi) \in GL_R(R^p)\rtimes GL_R(R^q)\rtimes Aut_\k(R)$ can be written as some (formal)
power series equations in the unknowns $U, V, \phi$.
\eei

Altogether, the equations i. for $J$, ii. for $G$ and the equations \eqref{Cond}
give a finite system of implicit equations
\beq \label{CondEqu2}
F(y) =0, \ \ F \in R[[y]]^s,
\eeq
in the unknowns  $\{u_{i,j}\}$,   $\{v_{i,j}\}$,  $\{\phi_i\}$, $\{v_i\}$ (the entries of $U, V, \phi, v$) and $\{t_{i,j}\}$  denoted all together by $y=(y_1,...,y_q)$.
(In the $C^\infty$-case the equations are of $C^\infty$-type.)

\beR
i. Suppose the ideal $J\sset \k[[x]]$ admits polynomial generators, the group
 $G\sseteq GL_R(M)\rtimes Aut_\k(R)$ is defined by polynomial equations, and $\{z_i\}$, $\{w_i\}$ in the expansion $z=\sum z_i\cdot \al(e_i)$, $w=\sum w_i\cdot \al(e_i)$ as well as the entries $\{a{_{i,j}}\}$ of the presentation matrix $A$ are polynomials.
  Then all the equations in \eqref{CondEqu2} are polynomial, i.e., $F(y)\in R[y]^s$.

ii. Similarly, if the data of $A,J,G,z, w$ are algebraic resp. analytic power series, then $F(y)\in R\langle y\rangle^s$ resp. $F(y)\in R\{y\}^s$.
\eeR

\subsection{The relevant approximation theorems}\label{Sec.Approximation.Theorems}
The determinacy criteria of sections  \ref{Sec.Finite.Determinacy.Poinwise.Lie.type}, \ref{Sec.Finite.Determinacy.Poinwise.weak.Lie.type}
 provide only an order-by-order solution to the  condition $gz=z+w$, $g\in G$.
 However, in many cases these are implicit function equations, see section \ref{Sec.Condition.z+w=gz.As.implicit.function.eq}.
 Then we can use  fundamental approximation theorems, to achieve ordinary solutions. We explain this in detail.

\subsubsection{The passage from an order-by-order solution to a formal solution}
When the coefficient ring is complete we use the following important result
(due to Pfister and Popescu).
\bthe\label{Thm.Popescu}
Let $(R,\cm)$ be a complete  Noetherian local ring (of arbitrary characteristic, not necessarily over a field).
Fix any $F( y)\in R[[ y]]^{s}, y=(y_1, \dots, y_q)$,
and  suppose that the system of equations $F( y)=0$ has an order-by-order solution.
(That is, there exists a sequence $\{ y_n \in \cm\cdot R^q\}_{n\ge1}$ such
that $F( y_n)\equiv 0 \ (mod\ \cm^n)$.)
Then there exists an ordinary solution, i.e., $ y\in R^{q}$ such that $F( y)=0$.
\ethe

A stronger version of this theorem  first appeared as \cite[Theorem 2.5]{Pfister-Popescu} for the particular case
$R=\k[[ x]]$, where $\k$ is either a field or a complete DVR, of zero characteristic.
For the case of arbitrary characteristic see \cite[Theorem 7.1]{Denef-Lipshitz}.
The generalization to the case where $(R,\cm)$ is a complete Noetherian local  ring is done as follows
(we thank D. Popescu for this explanation):
\bee[i.]
\item First one observes that if this property holds for $R$ then
it holds for any finite algebra over $R$ \cite[page 82, Satz 1.2]{Ku.Pf.Po.Ro.Mo}.
\item Then one uses the Cohen Structure Theorem for complete  Noetherian local rings,  \cite[Theorem 12]{Cohen}.
\eee


\subsubsection{The passage from a formal solution to an ordinary solution for polynomial or analytic equations}
If the equations $F( y)=0$ are polynomial, i.e., $F( y)\in R[ y]^{s}$, then {\em Artin's approximation theorem} \cite{Artin} can be used whenever $R$ has the Artin approximation propery.
\bed\label{Def.AP}
The Artin approximation property, AP, holds for a  ring $R$ with filtration $\{I_j\}$ if  for every finite system of polynomial
equations over $R$, a solution in the completion $\widehat{R}^{(I_j)}$ implies a solution in $R$, which moreover can be chosen arbitrarily close to the formal solution in the filtration topology.
\eed

The famous characterization of such rings reads:
\bthe\label{Thm.Artin.Popescu}\cite[Remark 2.15]{Popescu}
A commutative local Noetherian filtered ring with the filtration $\{\cm^j\}$  has the Artin approximation property if and only if it
 is excellent and Henselian.
\ethe
This result implies AP for rather general filtrations due to the following observation:

\bel\label{Thm.Approximation.general.filtration}\cite{B.B.K.19a}
Suppose a commutative, unital ring $R$ (not necessarily local or Noetherian) has AP   for   a filtration $\{I_j\}$.
  Then $R$ has AP for any   filtration by finitely generated ideals, $\{\ca_j\}$,
     such that $\ca_j\sseteq I_{d_j}$, for some sequence satisfying $\liml_{j\to \infty}d_j=\infty$.
\eel

\bex
The most important rings that have AP
  (assuming $\{I_j\}$ satisfy  $I_j\sseteq \cm^{d_j}$, with $\liml_{j\to\infty}d_j=\infty$):
\bei
\item $\k[[ x]]/J$ (trivially);
\item   $\k\langle x\rangle/J$, with $\k$ a field or a discrete valuation ring or a multi-variable formal power series ring;
\item  $\k\{ x\}/J$, with $\k$ a complete normed field
\eei
\eex
\beR\label{Thm.Approximation.W.systems}
 Sometimes the equations are not polynomial, e.g. for $\k\langle x\rangle/J$, $\k\{ x\}/J$  one has algebraic/$\k$-analytic equations.
  In this case  the approximation property of definition \ref{Def.AP} still holds:

 \begin{itemize}
\item  If $\k$ is a valued field then the analytic Artin Approximation theorem holds over $\k\{x\}/J$, if and only if the completion of $\k$ with respect to the absolute value is separable over $\k$, see \cite{Schemmel}.
%
%
  This condition is satisfied if either $\k$ is a field of characteristic zero, or a perfect field of prime characteristic.

\item For $\k$-algebraic approximation (and more generally, for W-systems),
we refer to \cite[Theorem 1.1]{Denef-Lipshitz}. The proofs in these papers are given for the  filtration $\{\cm^j\}$,    but the generalization to  $\{I_j\}$
    is done as in lemma \ref{Thm.Approximation.general.filtration}.
   \end{itemize}
\eeR

\subsubsection{The ${ C^\infty}$-case} \label{Sec.C.infty}
 The ring $C^\infty(\R^p,0)$ is not Noetherian and has no Artin approximation property, due to the ideal $\cm^\infty$ of flat functions
  whose Taylor series are identically zero.
For example, consider the equation $x_1\cdot y=\tau(x)$, here $\tau(x)\in \cm^\infty$, but  $\tau(x)\not\in (x_1)$.
The completion of this equation has the solution, $y=0$,  but the equation has no continuous solutions.

Yet, for $\R$-analytic equations there exists a stronger approximation property by
$C^\infty$-functions (the Taylor series of the smooth solution does not only approximate the formal solution but coincides with it):
\bthe\label{Thm.Approximation.Tougeron} \cite[Theorem 1.2]{Tougeron1976})
Let $F(x,y)\in \R\{x,y\}^s$ and suppose the equation $F(x,y)=0$ has a formal solution, $\hy(x)\in \R[[x]]^q$. Then there exists a smooth solution
 $y(x)\in (C^\infty(\R^p,0))^q$, whose Taylor series  is equal to $\hy$.
\ethe

For $C^\infty$-equations there is another approximation result:

\bthe\cite[Theorem 5.3]{B.K.IFT}\label{Thm.Approximation.Cinfty.equations}
Let $F(x,y)\in C^\infty(\R^p\times\R^n,0)$ and
 suppose the $\cm$-completion $\widehat{F(x,y)}=0$ has a formal solution $\hy_0(x)$.
Take a germ ${y}_0\in \big(C^{\infty}(\R^p,0)\big)^n$ whose Taylor series is $\hy_0(x)$.
Suppose  $\det\bbm F'_y(x,{y}_0)(F'_y(x,{y}_0))^t\ebm\cdot \cm^\infty=\cm^\infty$.
 Then there exists a $C^\infty$-solution, $y(x)\in C^\infty(\R^p,0)^n$, $F(x,y(x))\equiv0$,
  whose Taylor series at the origin is precisely $\hy_0(x)$.
\ethe

For the case of a general filtration, $\{I_j\}$, and for more statements of Artin-Tougeron type, see [B.B.K.19a].

\subsubsection{Transition from order-by-order determinacy to $G$-determinacy}\label{Sec.From.obo.determinacy.to.G.determinacy} \
 We want to pass from the stage ``$\overline{Gz} \supseteq\{z\}+M_{N+1}$" to the stage ``$Gz  \supseteq\{z\}+M_{N+1}$".
Equations   \eqref{CondEqu2} are often non-polynomial, thus one cannot use Artin approximation directly.
 However, with some natural weak assumptions one can ensure ordinary solutions.

Below  $\k$ is a ring, and $R\sseteq \k[[x]]/J$ is   a $\k$-subalgebra that contains the images of $\{x_i\}$.
 Let $M$ be a finitely presented $R$-module, with a filtration $\{M_j\}$ satisfying: $M_{N_j}\sseteq (x)^j\cdot M$ for any $j$ and some  $N_j<\infty$.    Then a subgroup  $G \sseteq GL_R(M)\rtimes Aut_\k(R)$  acts on $M$.
\bprop \label{Thm.Use.of.Approximation}
Suppose $G \sseteq GL_R(M)\rtimes Aut_\k(R)$ is defined by   power series, in the sense of condition ii.  of the explicit form \eqref{Cond}  at the end of section \ref{Sec.Condition.z+w=gz.As.implicit.function.eq}.
Suppose $\overline{Gz} \supseteq\{z\}+M_{N+1}$, i.e., $z \in M$ is  order-by-order $N$-determined.
\bee[1.]
 \item The $\cm$-adic completions satisfy:
 $\widehat{G} z  \supseteq\{z\}+\widehat M_{N+1}$, i.e., $z$ is formally $N$-determined.

 \item  Suppose $R$ has the Artin approximation property, $J\sset\k[[x]]$ has polynomial generators,
and  $G$ is defined by polynomial equations. Suppose also that the $G$-action preserves $M_{N+1}$
   and $z\in M$ is polynomially presented modulo $M_{N+1}$ (i.e., there exists a presentation $R^q\stackrel{A}{\to}R^p\stackrel{\al}{\to}M\to0$ with $A\in Mat_{p\times q}(\k[x]/J)$,
and an approximation $\tz_{pol}\in (\k[x]/J)^p$ such that
$\al(\tz_{pol})-z\in M_{N+1}$).
Then $Gz  \supseteq\{z\}+M_{N+1}$, i.e., $z$ is $N$-determined.

 \item   Let $R=\k\{x\}/J$, with $\k$ a complete normed ring. Suppose $G \sseteq GL_R(M)\rtimes Aut_\k(R)$ is  defined
  by analytic equations.
Then $Gz  \supseteq\{z\}+M_{N+1}$, i.e., $z$ is  analytically  $N$-determined.
 \item   Let $R=\k\langle x\rangle /J$, with $\k$ a field or a DVR. Suppose $G \sseteq GL_R(M)\rtimes Aut_\k(R)$ is  defined
  by algebraic power series equations.
 Then $Gz  \supseteq\{z\}+M_{N+1}$, i.e., $z$ is  algebraically  $N$-determined.
\eee
 \eprop
\bpr
We use notations from section \ref{Sec.Condition.z+w=gz.As.implicit.function.eq} and use the equations in \eqref{CondEqu2}  to present the condition
$z+w = gz$, $w \in M_{N+1}$  as a system of implicit function equations, $F(y) =0$,  for $F(y) \in R[[y]]^s$. Let $z = \alpha(\tilde z)$ and $w = \alpha(\tilde w)$.

\bee[ \it 1.]
\item  As the equations $F(y) =0$ are (formal) power series equations,
the existence of an order-by-order solution implies a formal solution
$y \in \cm\cdot \widehat{R}^q$ by Popescu's Theorem \ref{Thm.Popescu}.

\item  Let us assume first that $z\in M$ admits a polynomial presentation
with $A$
and $\Tilde z$ having polynomials components.
Then $F(y)$ is  polynomial in $y$ and we can use the  Artin Approximation (for the filtration $\{\cm^j\}$) to get an ordinary solution
$y \in \cm\cdot R^q$.
This implies $Gz  \supseteq\{z\}+M_{N+1}$.

The general case, when the components of $\tilde z$ are not polynomials, is reduced to the polynomial case as follows.
By the assumption $z$ is polynomial modulo $M_{N+1}$, so we fix some $\tz_{pol}\in (\k[x]/J)^p$ satisfying $\al(\tz_{pol})-z\in M_{N+1}$.
 The conditions $z=g\cdot \al(\tz_{pol})$, $g\in G$,
  are polynomial equations in the entries of $g$.  By part (1) these equations are formally solvable, so there exists $\hg\in \hG$ such that
  $z=\hg\cdot \al(\tz_{pol})$. Thus, by Artin approximation,
 there exists an ordinary solution, $z=g_0 \cdot \al(\tz_{pol})$, $g_0\in G$, i.e.,
$ z \in G\al(\tz_{pol}) $. Now the
 $G$-orbits of $z$ and $\al(\tz_{pol})$ coincide and we
 can continue with $\al(\tz_{pol})$ (using that $G$ preserves $M_{N+1}$).

\item In this case the equations   \eqref{CondEqu2} are analytic, $F(y) \in R\{y\}^s$.
 (Note that the entries of $A$ and of $\tz $ are in $\k\{x\}/J$.)
 Thus we apply the analytic approximation (remark \ref{Thm.Approximation.W.systems})
  to get a  solution $y \in  R^q$  analytic in $x$.

\item In this case the equations   \eqref{CondEqu2} are algebraic, $F(y) \in R\langle y\rangle ^s$.
 Apply the approximation of remark \ref{Thm.Approximation.W.systems}.
\eee

\beR
\bee[\bf i.]
\item The assumptions of part 2 are rather weak:
\bei
\item If $G$ is one of the groups  $GL_R(M)$,  $Aut_\k(R)$, $GL_R(M)\rtimes Aut_\k(R)$, or of the groups of section \ref{Sec.Examples.Actions.with.Aut} iv., then its defining equations are polynomial.
\item
Suppose the submodule $M_{N+1}\sset M$ satisfies $M_{N+1}\supseteq \cm^n\cdot M$, for $n\gg1$.
 Then  any element $z\in M$ is polynomial mod $M_{N+1}$, for any presentation.
 More generally, if $z\in I\cdot M$ then a sufficient condition (for polynomiality mod $M_{N+1}$) is: $M_{N+1}\supseteq \cm^n\cdot I\cdot M$, for $n\gg1$.
\eei
 \item
 The ring $C^\infty(\R^p,0)$ has no Artin approximation property, see \S\ref{Sec.C.infty}.
 Yet, we obtain the relevant approximation result in section \ref{Sec.Approximation.Cinfty.case}, using the structure of Lie pairs.


\item
In section \ref{determinacy for maps that preserve a subscheme}  we study determinacy relative to a germ. Besides
 equations \eqref{CondEqu2} this involves the condition  $\phi(I)=I$, for an automorphism $\phi\in Aut_\k(R)$ and some ideal $I\sset R$.
 This condition is treated as in the case of automorphisms of $\k[[x]]/J$, see (i) after equation \eqref{Cond}.
The approximation problem in this case is addressed as before, in particular it does not involve any nested Artin approximation.
  (See also \cite[Remark 2.28 ]{Rond}.)
\eee
\eeR

\section{Pairs of Lie type}  \label{Sec.Pairs.of.Lie.type}

Fix a $\k$-module $M$  with a filtration $\{M_i\}$. This induces   filtrations on $\{End^{(i)}_\k(M)\}$,  $\{GL^{(i)}_\k(M)\}$, see section
  \ref{Sec.Preparations.Induced.Filtration}. As in section \ref{Sec.Intro.MainConstruction} we consider pairs  $(T_{(G^{(1)},M)},G^{(1)}),$
where  $G^{(1)}\sseteq  GL^{(1)}_\k(M)$ is a (topologically unipotent) subgroup and  $T_{(G^{(1)},M)}\sseteq End^{(1)}_\k(M)$
 is a (topologically nilpotent) submodule. The submodule $T_{(G^{(1)},M)}$
generalizes the {\em tangent image} introduced in \cite{Greuel-Pham.2018}, i.e.,
for any $z\in M$ one may consider  $T_{(G^{(1)},M)}z$ as a substitute of the tangent space at $z$ of the orbit $G^{(1)}z\sseteq M$.

\subsection{Definitions and basic examples}
 The typical example for $(T_{(G^{(1)},M)},G^{(1)})$ is the pair $(Der^{(1)}_\k(R),Aut^{(1)}_\k(R))$, see section \ref{Sec.Examples.Actions.with.Aut}. As in the case of Lie groups/algebras we would like to ``integrate the vector fields", i.e. to assign to each derivation  $\xi\in Der^{(1)}_\k(R)$ a local coordinate change, i.e.  an automorphism $\phi_\xi\in Aut^{(1)}_\k(R)$.
 However, working with rings other than $C^\infty(\R^p)$, $\C\{x\}$, and without any topology/norm on $\k$,  we cannot use the standard results on differential equations.
As we have only the filtration topology, the naive definition is by using the exponential $exp(\xi)=\sum_j\frac{\xi^j}{j!}$.
 However each term $\frac{\xi^j}{j!}$ is not well defined, unless $\k\supseteq \Q$. Even if each $\frac{\xi^j}{j!}$ is well defined, there is no convergence notions for the  infinite sum, unless $R$ is complete. Therefore we are looking for  weak replacements of the exponential map $exp$ and the logarithmic map $ln$, preserving the passage $Der^{(1)}_\k(R)\rightleftarrows Aut^{(1)}_\k(R)$.  The would-be-exponential map should begin as $\xi\to \one+\xi+\dots$, and then one should specify the condition on the higher order terms ``$\dots$". We give four versions of this condition, they all mimic the vector field integration, and each is important in particular examples.

\

Given a pair  $(T_{(G^{(1)},M)},G^{(1)})$  we define the pairs $\big\{(T_{(G^{(i)},M)},G^{(i)})\big\}_i$:
\beq
 G^{(i)}:=G^{(1)}\cap GL^{(i)}_\k(M)\quad\quad \text{and}\quad\quad
   T_{(G^{(i)},M)}:=T_{(G^{(1)},M)}\cap  End^{(i)}_\k(M), \quad i\ge2.
\eeq
As in section \ref{Sec.Preparations.Induced.Filtration} the first intersection means that $g$ and $g^{-1}$ are of the form $\one + \phi$  with $\phi\in End_\k^{(i)}(M)$.

\bed
\label{Def.Pairs.of.Lie.Type}
\bee[1.]
\item  The pair $(T_{(G^{(1)},M)},G^{(1)})$ is called  {\bf of Lie type}  if  the following holds:
\bee[i.]
\item
For any $\xi\in T_{(G^{(1)},M)}$ there exists  $g\in G^{(1)}$, such
that for any $z\in M$,
\[\text{if $ord(\xi\cdot z)<\infty$ then $ord((g-\one-\xi)\cdot z) > ord(\xi\cdot z)$.}
\]
\item
For any $g\in G^{(1)}$  there exists  $\xi\in T_{(G^{(1)},M)}$, such
that for any $z\in M$,
\[\text{if $ord((g-\one)z)<\infty$ then $ord((g-\one-\xi)\cdot z) > ord((g-\one)z)$.}
\]
\eee

\item  The pair $(T_{(G^{(1)},M)},G^{(1)})$ is called  {\bf   of weak Lie type}    if the following holds:
\bee[i.]
%
\item
For any $\xi\in T_{(G^{(1)},M)}$ there exists  $g\in G^{(1)}$, such
that $ord((g-\one-\xi)\ge 2  ord(\xi)$.
\item
For any $g\in G^{(1)}$  there exists  $\xi\in T_{(G^{(1)},M)}$, such
that $ord((g-\one-\xi)\ge 2  ord(g-\one)$.
\eee


 \item
The pair $(T_{(G^{(1)},M)},G^{(1)})$ is called {\bf pointwise of Lie type} if the following holds:
\bee[i.]
\item For any $\xi\in T_{(G^{(1)},M)}$ and $z\in M$    there exists $g\in G^{(1)}$ satisfying:
\[
\text{
if $ord(\xi\cdot z)<\infty$ then
$ord((g-\one-\xi)z)>ord(\xi\cdot z)$}.\]
\item For any $g\in G^{(1)}$ and $z\in M$  there exists $\xi\in T_{(G^{(1)},M)}$ satisfying:
\[\text{if   $ord(g-\one)z)<\infty$ then $ord((g-\one-\xi)z)>ord(g-\one)z)$}.\]
\eee

\item
The pair $(T_{(G^{(1)},M)},G^{(1)})$ is called {\bf pointwise of weak Lie type} if the following holds  for any $i\ge1$:
\bee[i.]
\item For any $\xi\in T_{(G^{(i)},M)}$ and $z\in M$  there exists $g\in G^{(i)}$ satisfying:
\[\text{
 if $ord(\xi\cdot z)<\infty$ then $ord((g-\one-\xi)z)\ge 2ord(\xi)+ord(z)$.}\]
\item For any $g\in G^{(i)}$ and $z\in M$ there exists $\xi\in T_{(G^{(i)},M)}$ satisfying:
\[\text{
 if $ord(g-\one)z)<\infty$ then
 $ord((g-\one-\xi)z)\ge 2ord(g-\one)+ord(z)$.}\]
\eee
\eee
\eed

\beR
\bee[i.]

\item
Being of Lie type  implies pointwise Lie type. In most cases ``Lie type" implies also ``weak Lie type".

Note that ``$ord(\zeta\cdot z)\geq ord(\xi\cdot z)$ for any $z\in M$" implies
``$ord(\zeta)\geq ord(\xi)$". The converse does not hold,  as for some $z$ there may be cancelation among the
 lowest terms in $\xi\cdot z$ causing $ord(\xi\cdot z)>ord(\xi)+ord(z)$. Thus in many case ``weak Lie type" is weaker than ``Lie type".

\item In the definition of pointwise weak Lie type we ask: ``for any $\xi\in T_{(G^{(i)},M)}$  exists $g\in G^{(i)}$" rather than
 ``for any $\xi\in T_{(G^{(1)},M)}$  exists $g\in G^{(1)}$", and similarly from $g$ to $\xi$.
 This  condition for arbitrary $i$ is needed
  for the determinacy theorems in section 
 \ref{Sec.Finite.Determinacy.Poinwise.weak.Lie.type}.

\item Suppose $(T_{(G^{(1)},M)},G^{(1)})$ is a pair of pointwise weak Lie type and the following two conditions hold:
\bei
\item for any $\xi\in T_{(G^{(i)},M)}$, $z\in M$, with $ord(\xi\cdot z)<\infty$, exists $\txi\in T_{(G^{(i)},M)}$ satisfying:
\[
ord(\txi\cdot z-\xi\cdot z)>ord(\xi \cdot z) \quad\quad \text{and} \quad\quad ord(\txi\cdot z)=ord(\txi)+ord(z);
\]
\item
for any $g\in G^{(i)}$, $z\in M$, with $ord((g-\one)z)<\infty$, exists $\tg\in G^{(i)}$ satisfying:
\[
ord((\tg-g)\cdot z)>ord((g-\one)\cdot z)
 \quad\quad \text{and} \quad\quad
ord((\tg-\one)\cdot z)=ord(\tg-\one)+ord(z).
\]
\eei
Then (by direct check), $(T_{(G^{(1)},M)},G^{(1)})$ is a pair of  pointwise Lie type.

\item
The classical Lie groups (over $\R$, $\C$) admit the exponential and logarithmic maps
\beq
exp,ln:\ T_{(G^{(1)},M)}\rightleftarrows G^{(1)},\quad\quad exp(\xi)=\sum_{j=0}^\infty\frac{\xi^j}{j!},\quad\quad\quad ln(g)=\sum_{j=0}^\infty(-1)^j\frac{(g-\one)^j}{j}.
\eeq
In our general context these maps do not exist, both because of non-convergence and because of the denominators, e.g. when $char(\k)>0$.
 The natural weaker versions could be of the form
\beq\label{expln}
\Psi^{(exp)},\Psi^{(ln)}:\ T_{(G^{(1)},M)}\rightleftarrows G^{(1)},\quad\quad
\Psi^{(exp)}(\xi)=\one+\xi+F^{(exp)}(\xi), \quad\quad \Psi^{(ln)}(\xi)=(g-\one)+F^{(ln)}(g-\one).
\eeq
Here $F^{(exp)}(\xi)$, $F^{(ln)}(g-\one)$ should represent the higher order terms in the following sense:
 for any $\xi\in T_{(G^{(1)},M)}$, $g\in  G^{(1)}$ and $z\in M$ holds
\beq
ord(F^{(exp)}(\xi)\cdot z)\ge 2 ord(\xi\cdot z),\quad\quad\quad\quad  ord(F^{(ln)}(g-\one)\cdot z)\ge 2 ord((g-\one)z).
\eeq
Even this condition is too restrictive and does not hold e.g. for the group $Aut_\k(R)$.
The conditions of definition \ref{Def.Pairs.of.Lie.Type} are further weakening of \eqref{expln}
  that yet ensure the needed tight relation $T_{(G^{(1)},M)}\rightleftarrows G^{(1)}$.
  See \cite{B.B.K.19b} for some rings/derivations that admit the full exponential.

Recall that in differential geometry (over $\R,\C$) the map $\Psi^{(exp)}$ integrates a vector field to a flow, and is transcendental.
 Therefore, even if one assumes weaker conditions on $\Psi^{(exp)}$, in most cases one cannot expect  $\Psi^{(exp)}$ to be algebraic.
\eee\eeR

\bex\label{Ex.Lie.type.pair.over.Q} (The classical characteristic zero case).
Suppose $\k\supseteq\Q$ and assume that $G^{(1)}\sseteq GL^{(1)}_\k(M)$ and
let $T_{(G^{(1)},M)}\sseteq End^{(1)}_\k(M)$ be the tangent space. They often admit   ``order-by-order" exponential and logarithmic maps,
i.e., for any $\xi\in T_{(G^{(1)},M)}$ and $g\in G^{(1)}$ we have for $n\gg1$:
\beq\ber
\Psi^{(exp)}_n: \xi\to \one+\xi+\frac{\xi^2}{2!}+\cdots+\frac{\xi^n}{n!}\in G^{(1)}\cdot GL^{(n+1)}_\k(M), \\
\Psi^{(ln)}_n: g\to (g-\one)-\frac{(g-\one)^2}{2}+\cdots+(-1)^n\frac{(g-\one)^n}{n}\in T_{(G^{(1)},M)}+End^{(n+1)}_\k(M).
\eer\eeq
Then $(T_{(G^{(1)},M)},G^{(1)})$ is a pair of Lie type.
(The verification of  condition {\em (1)} in Definition \ref{Def.Pairs.of.Lie.Type} is immediate.)
The  simplest case when this happens is when $G^{(1)}$, $T_{(G^{(1)},M)}$ are complete \wrt their filtrations.
 Then, instead of the order-by-order maps, we take just the standard exponential and logarithmic maps,
$exp, ln: T_{(G^{(1)},M)} {\rightleftarrows}G^{(1)}$.
\eex

\bex \label{Ex.Lie.type.pair.GLM}
\bee[i.]
\item Suppose  \mbox{ $GL^{(1)}_\k(M)=\{\one\}+End^{(1)}_\k(M)$}, see section  \ref{Sec.Examples.of.groups}.
Here  one can take just the simplest maps: $\Psi^{(exp)}(\xi)=\one+\xi$ and $\Psi^{(ln)}(g)=g-\one$. The  pair $(End^{(1)}_\k(M), GL^{(1)}_\k(M))$ is then of Lie type.
\item If $M$ happens to be a filtered module over a larger filtered ring $R$, and one  has  $GL^{(1)}_R(M)=\{\one\}+End^{(1)}_R(M)$, then
the pair  $(End^{(1)}_R(M),GL^{(1)}_R(M))$ is of Lie type with maps $\Psi^{(exp)}$ and $\Psi^{(ln)}$ the same as in i.
\item More generally, suppose $G^{(1)}\sseteq GL^{(1)}_\k(M)$ is of the form $\{\one\}+\Lambda$, for some $\k$-submodule $\Lambda\sseteq End^{(1)}_\k(M)$.
Then the same maps, $\xi\mapsto \one+\xi$, $g\mapsto g-\one$, provide the pair of Lie type $(\Lambda,G^{(1)})$.
\eee\eex
For this example the maps $\Psi^{(exp)}$ and $\Psi^{(ln)}$ are non-unique, e.g., for $\k\supseteq\Q$ we could take the exponential and the logarithmic  maps or some of their approximations.

\bex\label{Ex.Lie.type.pairs.filtration}
Suppose a pair  $(T_{(G^{(1)},M)},G^{(1)})$ is of (pointwise) (weak) Lie type.
It follows directly from the conditions of Definition \ref{Def.Pairs.of.Lie.Type} that the pairs $(T_{(G^{(i)},M)}, G^{(i)})$ are of (pointwise) (weak) Lie type as well, i.e. the properties in Definition \ref{Def.Pairs.of.Lie.Type} hold also for $G^{(i)}$ instead of $G^{(1)}$.
\\If moreover the pair is given with some maps, $\{\Psi^{(exp)}_n,  \Psi^{(ln)}_n\}$, satisfying the relevant conditions of definition \ref{Def.Pairs.of.Lie.Type},
  then the maps restrict to the filtration,
\[
T_{(G^{(i)},M)}\stackrel{\Psi^{(exp)}_n}{\to}G^{(i)}\cdot GL^{(n+1)}_\k(M),
\quad\quad\quad
G^{(i)} \stackrel{\Psi^{(ln)}_n} {\rightarrow} T_{(G^{(i)},M)}+End^{(n+1)}_\k(M).
\]
\eex

\subsection{The pair $\boldsymbol{(Der^{(1)}_\k(R),Aut^{(1)}_\k(R))}$}\label{Sec.Pairs.of.Lie.Type.Aut(R)}
 We continue part i. and ii. of section \ref{Sec.Examples.Actions.with.Aut}.
Let $\k$ be a ring, $R$ be a $\k$-subalgebra of  $\k[[x]]/J$ or of $C^\infty(\R^p,0)/J$, $x=(x_1,\dots,x_p)$, that contains the images of $\{x_i\}$.
 Accordingly we consider the elements of $R$ as power series or function-germs and denote them $f(x)$.
 Take the maximal ideal $\cm=\langle x \rangle$. Fix a filtration
$\{I_j\}$  satisfying  $I_j\cdot I_k\sseteq I_{j+k}$, $I_1\sseteq\cm$.
This induces the filtration on the module of $\k$-linear derivations,
$$\begin{array}{lllc}
Der_\k^{(j)}(R)&: \ =&Der_\k(R)\cap End^{(j)}_\k(R) \\
&\ \, =& \{\xi\in End_\k(R)|\ \xi(ab)=\xi(a)b+a\xi(b),\ \text{and}\ \forall\ i:\ \xi(I_i)\sseteq I_{j+i}\}.
\end{array}$$

We take the module $Der_\k^{(1)}(R)$ as the tangent space for the (topologically unipotent) group $Aut^{(1)}_\k(R)$.

\subsubsection{The case $\boldsymbol{\k\supseteq\Q}$} \label{sec.3.2.1}
\bthe\label{Thm.Aut(R).pointwise.Lie} Let $\k$ be a commutative ring, $\k\supseteq\Q$.
\bee[1.]
\item If $R$ is complete \wrt  filtration $\{I_j\}$ then the pair $(Der_\k^{(1)}(R),Aut^{(1)}_\k(R))$ is of Lie type.
\item
 Suppose   $R$ is one of $\k\langle x \rangle/J$, $\k\{x\}/J$, or $C^\infty(\R^p,0)/J$.
 Suppose the filtration $\{I_j\}$ satisfies:
 $Der^{(1)}_\k(R)(x)\sseteq (x)^2$ and $I_{N+1}\cdot Der_\k(R)\sseteq Der^{(1)}_\k(R)$ for $N\gg1$.
  For $R=C^\infty(\R^p,0)$ we assume also:
\bei
\item $J$ is analytically generated;
\item $\cap^\infty_{j=1} I_j\supseteq\cm^\infty$ and  $\{I_j\}$ are analytically generated $mod\ \cm^\infty$.
\eei
  Then the pair $(Der_\k^{(1)}(R),Aut^{(1)}_\k(R))$ is pointwise of Lie type.
\eee
\ethe
\bpr
\bee[\em 1.]
\item
 As $R$ is complete \wrt $\{I_j\}$, and $\k\supseteq\Q$, we can use the classical exponential and logarithmic maps,
\beq\label{Eq.exp.ln.map}
\xi\to exp(\xi):=\suml^\infty_{j=0}\frac{\xi^j}{j!},\quad\quad
g\to ln(g):=\suml^\infty_{j=0}\frac{(-1)^{j+1}(g-\one)^j}{j}.
\eeq
Note that $exp(\xi)$, $ln(g)$ are well defined self-maps of $R$, because $R$ is complete and $\xi,(g-\one)$ are topologically nilpotent.

The map $exp(\xi)\circlearrowright R$ is  $\k$-linear and topologically unipotent.
 Moreover, it is multiplicative (by repeatedly applying the Leibniz rule for $\xi$).
 Finally, this map is invertible, as $exp(-\xi)\cdot exp(\xi)=Id=exp(\xi)\cdot  exp(-\xi)$.
Hence $exp(\xi)\in Aut^{(1)}_\k(R)$.

Similarly, $ln(g)\circlearrowright R$ is topologically nilpotent, $\k$-linear. Moreover, it satisfies the Leibniz rule (by repeatedly
 applying the multiplicativity of  $g$). Therefore $ln(g)\in Der^{(1)}_\k(R)$.

Finally we check the conditions of Definition \ref{Def.Pairs.of.Lie.Type}.
We get immediately:
\bee[i.]
\item  If $ord(\xi\cdot f)<\infty$   then $ord(\xi\cdot f)< ord(exp(\xi)-\one-\xi)f)$;
\item  If $ord((g-\one)\cdot f)<\infty$   then $ord((g-\one)\cdot f)< ord( (g-\one-ln(g))f)$.
\eee
Therefore the maps of equation \eqref{Eq.exp.ln.map} equip the pair $(Der^{(1)}_\k(R),Aut^{(1)}_\k(R))$ with the Lie-type structure.

\

\item First we record the useful form of the maps $exp(\xi)$, $ln(g)$ of equation \eqref{Eq.exp.ln.map} in the particular case, $R=\k[[ x]]/J$.
 For any power series $f( x)=\sum a_\um x^\um$ and any $\xi\in Der^{(1)}_\k(R)$, $g\in Aut^{(1)}_\k(R)$, we have (by part 1):
\beq\label{Eq.action.of.exp.ln}
\ber
exp(\xi)(f( x))=\sum exp(\xi)(a_\um x^\um)=\sum a_\um\Big(exp(\xi)( x)\Big)^\um=f(exp(\xi)( x)),
\\
ln(g)(f( x))=\sum_i ln(g)(x_i)\di_i f(x).
\eer
\eeq

\bee[i.]
\item We check condition {\em 3.i.} of definition \ref{Def.Pairs.of.Lie.Type}.

Let $f\in R$ and $\xi\in Der_\k^{(i)}(R)$ such that $ord(\xi\cdot f)<\infty$. Fix $N\ge1$ satisfying:
\beq
I_{N+1}\cdot Der_\k(R)\sseteq Der^{(i)}_\k(R),\quad \quad\quad
N>ord(\xi\cdot f),\quad\quad\quad
I_{N+1}\cdot Der_\k(R)(f)\sseteq I_{ord(\xi\cdot f)+1}.
\eeq
\bei
\item
First consider the case $J=0$.
We define $g\in Aut^{(i)}_\k(R)$ by $q(x)\to q(\suml^N_{j=0}\frac{\xi^j(x)}{j!})\in R$.

In more detail, write  $h(x)=\suml^N_{j=1}\frac{\xi^j(x)}{j!}\in R^{p}$.   The map $q(x)\to q(x+h(x))$ is $\k$-linear and multiplicative.
 It is invertible, as the equation $y+h(y)=x$ is resolvable over $R$. (Note that $h(y)\in (y)^2$ and apply $IFT_\one$.)
 Therefore $g\in Aut_\k(R)$.

 We claim:   $g\in Aut^{(i)}_\k(R)$. For any $q(x)\in I_l$ we should verify: $q(\suml^N_{j=0}\frac{\xi^j(x)}{j!})-q(x)\in I_{l+i}$, i.e.
 this difference goes to $0\in R/{I_{l+i}}$.
 For $R$ one of $\k\langle x \rangle$, $\k\{x\}$, or $C^\infty(\R^p,0)$  this can be checked at the level of Taylor power series of $q(x)$, i.e. we can take
  the completion and consider $q(x)\in \k[[x]]$. (For the case $R=C^\infty(\R^p,0)$ we use: $\capl I_j\supseteq\cm^\infty$.)

 Now, for  $q(x)\in \k[[x]]$ we can use the formal exponent:
 \beq\label{Eq.check.inside.proof.Aut}
q(\suml^N_{j=0}\frac{\xi^j(x)}{j!})=q\Big(exp(\xi)(x)-\suml^\infty_{j=N+1}\frac{\xi^j(x)}{j!}\Big)\in \suml^N_{j=0}\frac{\xi^j(q)}{j!}+I_{N+1}Der_\k(R)(q).
 \eeq
By the assumption, the right hand side here is contained in  $\suml^N_{j=0}\frac{\xi^j(q)}{j!}+Der^{(i)}_\k(R)(q)$.
Therefore the map $q(x)\to q(x+h(x))$ defines an element $g\in Aut^{(i)}_\k(R)$.

Finally we compare the orders, as in the  condition {\em 3.i.} of definition \ref{Def.Pairs.of.Lie.Type}.
We note, by (\ref{Eq.action.of.exp.ln}):
\beq
(g-\one-\xi)f=f\Big(\suml^N_{j=0}\frac{\xi^j}{j!}(x)\Big)-f(x)-\xi(f)\in \{\suml^N_{j=2}\frac{\xi^j}{j!}(f)\}+I_{N+1}.
\eeq
 This follows by checking the Taylor series, as in the equation \eqref{Eq.check.inside.proof.Aut}.

    Therefore, for the chosen $\xi,f$ and the constructed $g$, holds: $ord((g-\one-\xi)f)>ord(\xi\cdot f))$.
     Note that $g$ depends on $\xi$ and $f$, via $N$.

\item For the general case, $J\neq 0$,  the map $q(x)\to q(\suml^N_{j=0}\frac{\xi^j(x)}{j!})$ is not well defined, as it
 does not necessarily preserve   $J\sset R$, see example \ref{Ex.problem.non-reg.rings}. We adjust this map by higher order terms as follows.

 Let $S$ be one of $\k\langle x\rangle$, $\k\{x\}$,  $C^\infty(\R^p,0)$, and  take the full preimages of $\{I_j\}$ to get the filtration $\{\tI_j\}$ of $S$.
 For $C^\infty(\R^p,0)$ one has   $\tI_j\supseteq J+\cm^\infty$, for any $j$, and $\tI_j$ are analytically generated $mod(\cm^\infty)$.

  For $\xi\in Der^{(i)}_\k(R)$ take its representative $\txi\in Der^{(i)}_\k(S)$ then $\txi(J)\sseteq J$.

 In the $C^\infty(\R^p,0)$ case the coefficients of  $\txi=\sum\txi_i\di_i$ are smooth functions. But there exists an analytic derivation,
  $\tilde\txi\in Der^{(i)}_\R(\R\{x\})$, satisfying: $\tilde\txi(J)\sseteq J$ and $\tilde\txi-\txi\in Der^{(2N)}_\R(C^\infty(\R^p,0))$.
 Indeed, fix some analytic generators $\{q^J_\al\}$ of $J\sset S$, and expand $\txi=\sum \txi_i\di_i$. Then the condition $\txi(J)\sseteq J$
  amounts to the $\R\{x\}$-linear equations $\sum \txi_i\di_i q^{J}_\al=\sum t_{\al\be}q^J_\be$. (Here the unknowns are $\{\txi_i\}, \{t_{\al\be}\}$.)
   These equations posses the formal solution, the image of $\txi$ under the completion.
    Therefore, by Artin approximation, there exists an analytic solution,
   $\tilde\txi$, that approximates $\txi$, see remark \ref{Thm.Approximation.W.systems}.
We replace $\txi$ by this analytic derivation $\tilde\txi$.

Now   we construct the map $\tphi\in Aut_\k^{(i)}(S)$ by
    $q(x)\to q(\suml^N_{j=0}\frac{\txi^j(x)}{j!}+\tilde{h}(x))$, where $\tilde{h}(x)\in \tI_{N+1}\cdot S^p$ and $\tphi(J)=J$.
    Fix some (finite)  sets of analytic generators, $\{q_\al^{(J)}\}$ of $J$, and $\{q_\be^{(I)}\}$ of $I_{N+1}$.
     Expand $\tilde{h}(x)$ as $\sum t_\be q_\be^{(I)}$, for the new variables $\{t_\be\}$. Then the condition $\tphi(J)=J$ amounts to:
\beq
q_\al^{(J)}\Big(\suml^N_{j=0}\frac{\txi^j(x)}{j!}+\sum t_\be q_\be^{(I)}\Big)=\sum_\be w_{\al\be}q_\be^{(J)},\quad\quad \forall\ \al.
\eeq
This is a finite system of implicit function equations on the variables $\{w_{\al\be}\}$, $\{t_\be\}$. The equations are
  algebraic/analytic power series. (Also in the $C^\infty$ case!) Thus, by the relevant approximation property,
  it suffices to demonstrate a formal solution. Note that $\suml^\infty_{j=N+1}\frac{\txi^j(x)}{j!}\in I_{N+1}$, thus we can take $\{t_\be\}$ satisfying
   $\suml^\infty_{j=N+1}\frac{\txi^j(x)}{j!}-\sum_\be t_\be q_\be^{(I)}=0$. This gives the formal  solution,
\beq
q_\al^{(J)}(\suml^N_{j=0}\frac{\txi^j(x)}{j!}+\sum_\be t_\be q_\be^{(I)})=q^{(J)}_\al(exp(\txi)(x))=
exp(\txi)(q^{(J)}_\al)=\sum_\be w_{\al\be}q_\be^{(J)}\in J\cdot \hat{S},
 \quad\quad \forall\ \al.
\eeq
Now, by approximation, we get the needed analytic/algebraic solution.

Altogether, we have constructed a coordinate change $\tphi_\txi:\ q(x)\to q(\suml^N_{j=0}\frac{\txi^j(x)}{j!}+\tilde{h}(x))$ that preserves $J\sset S$.
 It descends to the element $\phi_\txi\in Aut_\k(R)$. As in the case $J=0$, one checks that $\phi_\txi$ is topologically unipotent
  and  $ord((\phi_\txi-\one-\xi)f)>ord(\xi\cdot f) $.
\eei

\item
(We check condition {\em 3.ii.} of definition \ref{Def.Pairs.of.Lie.Type}.)  Let $f\in R$ and  $g\in Aut^{(i)}_\k(R)$ such that $ord((g-\one)f)<\infty$.
Fix $N\gg1$ satisfying:
\beq
I_{N+1}\cdot Der_\k(R)\sseteq Der^{(i)}_\k(R),\quad \quad
N>ord((g-\one)f),\quad
I_{N+1}\cdot Der_\k(R)(f)\sseteq I_{ord((g-\one)\cdot f)+1}.
\eeq
\bei
\item
First assume $J=0$.  We define   $\xi\in Der_\k^{(i)}(R)$
 by
\beq
 \xi(q(x))=\suml^N_{j=1}\frac{(-1)^{j+1}(g-1)^j(x)}{j}\cdot q'(x)=\suml^N_{j=1}\suml_i\frac{(-1)^{j+1}(g-1)^j(x_i)}{j}\frac{\di q(x)}{\di x_i}.
\eeq
  This is a derivation by construction. (Note that $\k\langle x\rangle$ is closed under differentiation.

  To check the topological nilpotence we should verify: if $q(x)\in I_l$ then $\xi(q(x))\in I_{l+i}$, i.e.
   $\xi(q(x))$ goes to $0\in R/I_{l+i}$. In the same way as in the check of topological unipotence,  we can pass to completion and
    use the formal relation $ln(g)q(x)=ln(g)(x)\cdot q'(x)$,     of equation \eqref{Eq.action.of.exp.ln}.
(For $C^\infty$ case  we use $\cap I_j\supseteq\cm^\infty$, as before.)
   Therefore
\beq
\suml^N_{j=1}\frac{(-1)^{j+1}(g-1)^j(x)}{j}\cdot q'(x)\in End^{(i)}_\k(R)(q)+I_{N+1}\cdot Der_\k(R)(q).
\eeq

Finally, for the condition {\em 3.ii.} of definition \ref{Def.Pairs.of.Lie.Type} we note, using equation \eqref{Eq.exp.ln.map}:
\beq
(g-\one-\xi)f=\suml^N_{j=2}\frac{(-1)^{j+1}(g-1)^j(x)}{j}f'(x)
=ln(g)f'-(g-\one)f-\suml^\infty_{j=N+1}\frac{(-1)^{j+1}(g-1)^j(f)}{j}.
\eeq
And the later belongs to
$\{ln(g)f'-(g-\one)f\}+I_{N+1}$ i.e. to $\{\suml^N_{j=2}\frac{(-1)^{j+1}(g-1)^j}{j}(f)\}+I_{N+1} $.

 Thus we have constructed $\xi\in Der^{(i)}_\k(R)$ satisfying $ord\big((g-\one-\xi)f\big)> ord((g-\one)f)$. Note that $\xi$ depends on $g,f$ (via $N$).
\item For the general case, $J\neq\{0\}$, we cannot use $\Big(\suml^N_{j=1}\frac{(-1)^{j+1}(g-1)^j(x)}{j}\Big) q'(x)$ as this does not necessarily preserve $J$.
 Thus we adjust this derivation,
 \beq
 \xi(q(x)):=\Big(\suml^N_{j=1}\frac{(-1)^{j+1}(g-1)^j(x)}{j}+h(x)\Big) q'(x).
 \eeq
 Here $h\in I_{N+1}$   is chosen to ensure $\xi(J)\sseteq J$.

  First we provide a formal solution for this condition,  $\hat{h}(x)=\suml^\infty_{j=N+1}\frac{(-1)^{j+1}(g-1)^j(x)}{j}$.
  Now as in the step {\em 2.i.}, for the case $R=\k\langle x\rangle/J$, $\k\{x\}/J$ we use the relevant approximation directly.
    For $C^\infty$-case  we first modify $\xi$
    by high enough terms to ensure that it is an analytic/polynomial derivation.

 Thus the needed $h(x)\in R$ exists in all the cases and we get
   $\xi\in Der^{(i)}_\k(R)$. The condition $ord\big((g-\one-\xi)f\big)> ord((g-\one) f)$ is checked as before.
\epr
\eei
\eee
\eee
\bex As a special case of Theorem \ref{Thm.Aut(R).pointwise.Lie} we get:
Let  $\k\supseteq \Q$ (e.g $\k$ a field of characteristic 0) and $R$ one of $\k[[x]]/J$, $\k\{x\}/J$, $\k\langle x\rangle/J$ and  $C^\infty(\R^p,0)/J$ ($J$ analytically generated), with filtration $\{\mathfrak{m}^j\}_j$.
 Then the pair $(Der_\k^{(1)}(R),Aut^{(1)}_\k(R))$ is of pointwise Lie type.
\eex
\beR The constructions in part 2 of this proof are somewhat involved, for various reasons.
\bee[\bf i.]
\item
 We could not use the full power series  $exp(\xi)$, $ln(g)$, as they do not always act on rings like  $\k\langle x\rangle$  and  $C^\infty(\R^p,0)$.
For example, for any derivation $D$  there exists a smooth function $f$ such that $(\sum \frac{D^j}{j!})f$ diverges at each point off the origin, \cite{B.B.K.19b}. Thus we cannot define the action  $exp(D)\circlearrowright C^\infty(\R^p,0)$ via
  $exp(D)=\sum \frac{D^j}{j!}$.

 Similarly, for $D=(x^p+x^{p+1})\frac{d}{d x}$, $p\ge 2$, the power series $exp(D)(x)$ is non-algebraic (though it is analytic for $\k=\R,\C$).
Therefore $exp(D)$ does not act on $\k\langle x\rangle$.
\item
 We could replace $exp(\xi)$, $ln(g)$ by their n'th jet approximation,
\beq
\Psi^{(exp)}_n: \xi\to \one+\xi+\frac{\xi^2}{2!}+\cdots+\frac{\xi^n}{n!},\quad   \quad
\Psi^{(ln)}_n: g\to (g-\one)-\frac{(g-\one)^2}{2}+\cdots+(-1)^n\frac{(g-\one)^n}{n}.
 \eeq
This brings another problem: the map $f\to \Psi^{(exp)}_n(\xi)(f)$ is not necessarily multiplicative, and for small $n$ not necesarily
invertible.
 For example, let $R=\C\{x\}$ and $\xi=x^2\di_x$.
Then the solution to $(\one+\xi)y=x$ is $y=\sum n!(-1)^nx^{n+1}$, which is diverging off the origin.
 Therefore  instead of this map, we had to use the map $f\to f(\Psi^{(exp)}_n(\xi)(x))$, for $\Psi^{(exp)}_n$ as above,
and then  to adjust it in the case of $R=S/J$.
\eee
\eeR

\beR\label{Ex.problem.non-reg.rings} Another problem is due to non-regular rings.
  Let $\k\supseteq\Q$ and  $R=\k[[x,y]]/{(x^p-y^q)}$, $p<q$, filtered by $\{\cm^j\}$.
\bee[i.]
\item
Let $\xi=x(q\cdot x\di_x+p\cdot y\di_y)$,
 so that $\xi\cdot (x^p-y^q)\in \cm\cdot (x^p-y^q)$. Thus $\xi\in Der^{(1)}_\k(R)$. However, for $f(x,y)=x^p-y^q$,
\[
f\big((x,y)+\xi(x,y)\big)=f(x+qx^2,y+pxy)=x^p(1+qx)^p-y^q(1+px)^q\not\in (x^p-y^q).
\]
Thus $f(x,y)\to f\big((x,y)+\xi(x,y)\big)$ is not an automorphism of $R$.
\item Let $g:(x,y)\to(x+y^q,uy)$, where $u\in \k[[x,y]]$ is defined by the condition $u^q=1+\frac{(x+y^q)^p-x^p}{y^q}$, $u(0,0)=1$.
  Then $g(f)\in (f)$ thus $g\in Aut^{(1)}_\k(R)$. But $\xi:=\sum (g-\one)(x)\di_i$ is not a derivation of $R$, as $\xi(f)=y^q(-px^{p-1}+(u-1)q)\not\in(f)\sset\k[[x,y]]$.
\eee
\eeR

\beR
As one sees from the proof, the statement of part 2  of Theorem 3.6 holds for more general rings.
 We only use in the proof the following: $R$ is a ring over a field of zero characteristic, $R$ is closed under differentiation,
 and the approximation property in a ring $R$ holds for the particular class of equations over $R$  that we need to solve.
 For instance, our proof works when $R=W_n/J$, where $\{W_n\}_{n\geq 1}$ is either a Weierstrass system over $\k$,  see
 \cite[Definition of page 2]{Denef-Lipshitz},
or   a convergent Weierstrass system  over $\R$, \cite[page 798]{vandenDries}. (The filtration is $\{\cm^j \}$.)
\eeR

\subsubsection{ The pair $\boldsymbol{ Der^{(1)}_\k(R),Aut_\k(R)}$ for the case of arbitrary $\boldsymbol{\k}$}
In the previous section  \ref{sec.3.2.1} we considered the case $\k \supset \Q$ and proved conditions for pointwise Lie type in Theorem \ref{Thm.Aut(R).pointwise.Lie}. Now we
do not assume that the ring $\k$ contains $\Q$ and establish  conditions to guarantee weak Lie type under certain Assumptions \ref{Assumptions.on.the.ring}.

Let $R$ be a filtered subring of $\k[[x]]/J$, $x=(x_1,\dots,x_p)$.
  We assume that the filtration satisfies  $I_1\sseteq\cm=(x)$ and $Der^{(1)}_\k(R)(x)\sseteq (x)^2$.
 Moreover, in this section we make the following assumptions about  $R$ and the filtration:

\begin{Assumptions} \label{Assumptions.on.the.ring}
\bee[i.]
\item Elements of  $Der^{(1)}_\k(R)$ induce coordinate changes: for any $\xi\in Der^{(1)}_\k(R)$ there exists
 $h^{(\xi)}(x)=(h_1(x),\dots,h_p(x))$, with $ord(h_i(x))\ge 2ord(\xi(x_i))$, such that $f(x)\to f(x+\xi(x)+h(x))$ is a well defined morphism of $R$, i.e. maps $J$ to $J$.
\item
 Elements of $Aut^{(1)}_\k(R)$ induce  derivations: for any $g\in Aut^{(1)}_\k(R)$ there exists $h^{(g)}(x)=(h_1(x),\dots,h_p(x))$,
  with $ord(h^{(g)}_i(x))\ge 2ord((g-\one)(x_i))$,
 such that $\sum_i\big( (g-\one)x_i+h^{(g)}_i\big)\di_i \in Der^{(1)}_\k(R)$.
\item
 $R$ admits Taylor expansion  up to second order for coordinate changes, i.e., for any $\xi\in Der^{(1)}_\k(R)$ and a corresponding $h^{(\xi)}(x)$
  and any  $f(x)\in R$ holds
  $f(x+\xi(x)+h^{(\xi)}(x))=f(x)+ \xi\cdot f(x)+F(\xi,f,x)$, where   $F$ satisfies $F(Der^{(j)}_\k(R),I_i,x)\sseteq I_{2j+i}$ for any $i\ge1$, $j\ge1$.
\item $IFT_\one$ holds over $R$, i.e., any equation $ y+F( y)= x$, with $F( y)\in ( y)^2$ is solvable for $y$ over $R$.
 \eee
\end{Assumptions}

\bex
These assumptions  hold for the rings  of Theorem \ref{Thm.Aut(R).pointwise.Lie}.
(Assumption i. is verified after equation (37) and assumption ii.  is verified after equation (41), and $IFT_\one$ holds for those rings.)

We verify the assumptions also for  $R$   one of $\k[[x]]/J$,  $\k\langle x\rangle/J$, $\k\{x\}/J$,  where $J\sset\cm^2$ and $\k$   is a (complete normed) ring.
\bee[i.]
\item Condition  $i.$   holds trivially when $J=0$ (i.e. regular rings), with any filtration.
 One just takes  $h^{(\xi)}=0$.
 If $J\neq0$ then condition  $i.$  is non-trivial  and the non-zero correction  $h^{(\xi)}$ is needed, see remark \ref{Ex.problem.non-reg.rings}.
\item
 Condition  $ii.$   holds trivially when $J=0$ (one takes $h^{(g)}=0$) for the filtrations $\{\cm^j\}$ or $I_{j+1}=\ca^j\cdot I_1$, $I_1\sseteq\ca^2\sseteq\cm^2$.
 More generally, a sufficient condition on the filtration is: $(I_1\cap (x)^2) \di I_j\sseteq I_{j+1}$.

 If $J\neq0$ then condition  $ii.$  is non-trivial  and the non-zero correction  $h^{(g)}$ is needed, see remark \ref{Ex.problem.non-reg.rings}.

\item Condition  $iii.$ holds for our rings $\k[[x]]/J$,  $\k\langle x\rangle/J$, $\k\{x\}/J$, see example \ref{Ex.rings.with.IFT1}

\item
Condition $iv.$ holds  for filtrations $\{I_j=\cm^j\}$ and
   $\{I_{j+1}=\ca^j\cdot I_1\}$, with $I_1\sseteq\ca^2\sset\cm$. More generally, it holds
 when the filtration $\{I_j\}$ satisfies $I_j\cdot I_j\cdot Der_\k(R)\Big(Der_\k(R)(I_i)\Big)\sseteq I_{2j+i}$.

\eee
\medskip
In particular, Assumptions \ref{Assumptions.on.the.ring} hold for $\k[[x]]$,  $\k\langle x\rangle$, $\k\{x\}$ and for the filtrations $\{\cm^j\}$
  or $I_{j+1}=\ca^j\cdot I_1$, $I_1\sseteq\ca^2\sseteq\cm^2$.
\eex

\bex\label{more examples when conditions hold} Here are  more examples of rings satisfying the assumptions,  for the filtration $\{\cm^j\}$.
Let $\k$ be either a field or a discrete valuation ring, and $\{W_n\}_{n\geq 1}$
 a Weierstrass system over $\k$. Then, for any integer $n\geq 1,$ the ring
$W_n$ satisfies \ref{Assumptions.on.the.ring}.

Indeed, $R$ is closed under compositions, see  \cite[page 3]{Denef-Lipshitz}, and admits Taylor expansions up to second order, see \cite[page 4]{Denef-Lipshitz}.
The $IFT_\one$ holds over $R$, see Example \ref{Ex.rings.with.IFT1}.
Finally, $R$ is a $\k$--subalgebra of $\k[[x]]$, therefore the
topologically unipotent $\k$--automorphisms of $R$ come as coordinate changes.
 Thus condition $iv.$ holds for $R$.
\medskip

%
\eex

\bprop \label{Thm.Aut(R).of.weak.Lie.type}
Suppose $R$ satisfies Assumptions  \ref{Assumptions.on.the.ring}. Then the pair  $(Der_\k^{(1)}(R),Aut^{(1)}_\k(R))$ is of weak Lie type.
\eprop
\bpr
\bee[\bf i.]
\item (Verifying condition {\em 2.i.}  of definition \ref{Def.Pairs.of.Lie.Type})

For any $\xi\in Der_\k^{(1)}(R)$ we take the corresponding coordinate change, $f(x)\to f(x+\xi(x)+h^{(\xi)}(x))$.
 This is a well defined (set-theoretic) self-map of $R$, by assumption $i.$ It is $\k$-linear and multiplicative by construction.
 It is invertible by $IFT_\one$.  Thus we get an element $g_\xi\in Aut_\k(R)$.

By Taylor-expansion (assumption \ref{Assumptions.on.the.ring} $iii.$), we have: $f( x+\xi( x)+h^{(\xi)}(x))=f( x)+\xi\cdot f( x)+F(\xi,f, x)$. Here, for $f\in I_i$, $\xi\in Der_\k^{(j)}(R)$
we have: $F(\xi,f, x)\in I_{i+2j}$, thus $g-\one-\xi\in End^{(2j)}_\k(R)$.
 In particular,  $g_\xi$ is topologically unipotent and  $ord(g-\one-\xi)\ge 2ord(\xi)$.

\item (Verifying condition {\em 2.ii.}of definition \ref{Def.Pairs.of.Lie.Type})
To any $g\in Aut^{(1)}_\k(R)$ we associate $\xi_g\in Der^{(1)}_\k(R)$ by
\beq
\xi_g \cdot f(x):=\sum \big((g-\one)(x_i)+h^{(g)}_i(x)\big)\di_i f(x).
\eeq
 This is a derivation, $\xi\in Der^{(1)}_\k(R)$,  by assumption \ref{Assumptions.on.the.ring} $ii$.

Now we check:
\beq
(g-\one)f(x)-\xi_g\cdot f(x)=f(x+(g-\one)(x))-f(x)-\xi_g\cdot f(x)=F((g-\one)(x),f,x)-\sum h^{(g)}_i(x)\di f.
\eeq
 Thus, by assumption \ref{Assumptions.on.the.ring}.ii,  $ord((g-\one-\xi_g)f)\ge 2ord((g-\one)f)$.
As this holds for any $f$ we get:  $ord(g-\one-\xi_g)\ge 2ord(g-\one)$.  \epr
\eee

\bex \label{ex3.10}
For $\k$ a field, the rings $\k[[ x]]$, $\k\{ x\}$, $\k\langle x\rangle$ as well as Weierstrass system and functions of the Denjoy Carleman class,
with filtration $\{\cm^j\}_j$ or $\{\ca^j I_1\}$, $I_1\sseteq\ca^2 \sseteq \cm^2 $,
satisfy the assumptions \ref{Assumptions.on.the.ring}. Therefore for these rings the pair $(Der_\k^{(1)}(R),Aut^{(1)}_\k(R))$ is of weak Lie type.
\eex

\subsubsection{ The subgroup $\boldsymbol{Aut_{\k,\ca}(R)\sset Aut_\k(R)}$}\label{Sec.(w)Lie.type.of.relative.group}
Sometimes one considers only those automorphisms of the ring that preserve a given ideal $\ca$,
\beq
Aut_{\k,\ca}(R):=\{\phi\in Aut_{\k}(R)|\ \phi(\ca)=\ca\}.
\eeq
Accordingly, one considers the module of $\ca$-logarithmic derivations (also known as the module of $\ca$-preserving
derivations \cite[Definition 2.1.1]{Brumatti-Simis95}),
$Der_{\k,\ca}(R):=\{\xi\in Der_\k(R)|\ \xi(\ca)\sseteq \ca\}$. This is an $R$-submodule of $Der_\k(R)$.

As before, we get the induced filtrations of the group $\{Aut^{(i)}_{\k,\ca}(R)\}$ and of the module $\{Der^{(i)}_{\k,\ca}(R)\}$.
\bprop Suppose $\k\supseteq\Q$. \label{Prop.3.16}
\bee[1.]
\item
If $R=\k[[ x]]/J$ then $(Der^{(1)}_{\k,\ca}(R),Aut^{(1)}_{\k,\ca}(R))$ is a pair of Lie type for any filtration.
\item
 Suppose   $R$ is one of $\k\langle x \rangle/J$, $\k\{x\}/J$, or $C^\infty(\R^p,0)/J$.
 Suppose the filtration $\{I_j\}$ satisfies:
 $Der^{(1)}_\k(R)(x)\sseteq (x)^2$ and $I_{N+1}\cdot Der_\k(R)\sseteq Der^{(1)}_\k(R)$ for $N\gg1$.
 For $R=C^\infty(\R^p,0)$ we assume also:
\bei
\item $J,\ca$ are analytically generated;
\item $\cap I_j\supseteq\cm^\infty$ and  $\{I_j\}$ are analytically generated $mod\ \cm^\infty$.
\eei
  Then the pair $(Der_{\k,\ca}^{(1)}(R),Aut^{(1)}_{\k,\ca}(R))$ is pointwise of Lie type.
\eee
\eprop
The proof is the same as for Theorem \ref{Thm.Aut(R).pointwise.Lie}.   The correspondences  $Der^{(1)}_{\k,\ca}(R)\ni \xi \rightleftarrows g\in  Aut^{(1)}_{\k,\ca}(R)$,
 are defined in the same way, one just imposes the condition that the images preserve $\ca$. (This results in algebraic resp. analytic equations).

\beR\bee[i.]
\item
Note that the natural homomorphism of groups $Aut_{\k,\ca}(R)\to Aut_{\k}(R/\ca)$ is not necessarily injective or surjective.
Neither is the natural map $Der_{\k,\ca}(R)\to Der_{\k}(R/\ca)$, though it is known to be surjective in some particular cases,
e.g., for $R$ one of $ \k[ x],\k[[ x]]$, \cite[Lemma 2.1.2]{Brumatti-Simis95}, \cite[p. 2717]{Narvaez.Macarro} or when $R$ is a regular local ring of zero
characteristic with some technical conditions (see \cite[Theorems 30.6 and 30.8]{Matsumura} or
\cite[Theorem 2.1]{KallstromTadesse2015}).
Therefore the determinacy problem for the pair $(Der_{\k,\ca}(R),Aut_{\k,\ca}(R))$
is in general not reduced to the pair $(Der_{\k}(R/\ca),Aut_{\k}(R/\ca))$.
 Thus we need this special proposition.
 \item
  The analogue of Proposition \ref{Prop.3.16} for weak Lie type is more technical. We need to assume  the conditions \ref{Assumptions.on.the.ring} for $R$ and for $R/\ca$,
   and moreover some lifting properties. Therefore we omit it.
\eee
\eeR

\subsection{Constructing new Lie pairs from old ones}
\subsubsection{Diagonal action}\label{Sec.Diagonal.Action}
Suppose a pair of Lie type  (resp. pointwise and/or weak Lie type) acts on several modules. Namely, we have a collection of (not necessarily injective) homomorphisms
\beq
\Big\{(T_{(G^{(1)},M)},G^{(1)})\stackrel{\psi_\al}{\to}(End^{(1)}_\k(M_\al),GL^{(1)}_\k(M_\al))\Big\}_\al,
\eeq
such that the conditions of Definition \ref{Def.Pairs.of.Lie.Type} are satisfied for the pairs $\{(\psi_\al(T_{(G^{(1)},M)}),\psi_\al(G^{(1)})\}_\al$.

Take the diagonal action, $(T_{(G^{(1)},M)},G^{(1)})\circlearrowright \oplus_\al M_\al$, by $(\xi,g)\cdot \{z_\al\}:=(\{\psi_\al(\xi) z_\al\},\{\psi_\al(g) z_\al\})$.
 This corresponds to the diagonal homomorphism, $(T_{(G^{(1)},M)},G^{(1)})\stackrel{\De}{\to}\oplus_\al (T_{(G^{(1)},M)},G^{(1)})$.
\bel
Suppose each pair $(\psi_\al(T_{(G^{(1)},M)}),\psi_\al(G^{(1)})$is of  Lie type (resp. pointwise and/or weak Lie type).
Then the same holds for the pair $(\De(T_{(G^{(1)},M)}),\De(G^{(1)}))$.
\eel
\bpr
 By the assumption to every $\oplus\psi_\al(\xi)$ is assigned $\oplus \psi_\al(g)$ and vice versa.
 To verify the (relevant) conditions of definition \ref{Def.Pairs.of.Lie.Type} we should compare the orders of
 $\oplus \psi_\al(g-\one-\xi)$ to the orders of $\oplus \psi_\al(g-\one)$,
 $\oplus \psi_\al(\xi)$.
(For the pointwise case we should compare the order  $ \oplus \psi_\al(g-\one-\xi)(z)$ to the orders of $\oplus \psi_\al(g-\one)(z)$, $\oplus \psi_\al(\xi)(z)$.)

 For each case it is enough to observe: $ord(\oplus w_\al)=min\{ord(w_\al)\}$.
\epr
\bex
Suppose the pair $(Der^{(1)}_\k(R),Aut^{(1)}_\k(R))$ is of  Lie type (resp. pointwise and/or weak Lie type) for the action on $R$. Then the same holds  for the diagonal actions  $Aut_\k^{(1)}(R)\circlearrowright R^{\oplus n}$, $Aut_\k^{(1)}(R)\circlearrowright \Mat$.
\eex

\subsubsection{A group generated by two groups.}
Fix some subgroups $G^{(1)},H^{(1)}\sseteq  GL^{(1)}_\k(M)$ and denote the subgroup they generate by $\langle G^{(1)},H^{(1)}\rangle$.
 The Lie-type structures of $G^{(1)},H^{(1)}$ often induce the one on $\langle G^{(1)},H^{(1)}\rangle$ as we show now.

In this subsection we assume: $\langle G^{(1)},H^{(1)}\rangle=G^{(1)}\cdot H^{(1)}$, i.e., any
element of  $\langle G^{(1)},H^{(1)}\rangle$ is presentable in the form $g\cdot h$, for some  $g\in G^{(1)}$, $h\in H^{(1)}$.
(This holds e.g., for semi-direct products.)

\bel\label{Thm.Lie.type.of.(H.G)}
Suppose  $\langle G^{(1)},H^{(1)}\rangle=G^{(1)}\cdot H^{(1)}$. Suppose moreover that for any $i\ge1$ the following holds:
\[(T_{(G^{(1)},M)}+T_{(H^{(1)},M)})\cap End^{(i)}_\k(M)=T_{(G^{(i)},M)}+T_{(H^{(i)},M)},\quad
 (G^{(1)}\cdot H^{(1)})\cap GL^{(i)}_\k(M)=G^{(i)}\cdot H^{(i)}.
\]
\bee[1.]
\item  If the pairs $(T_{(G^{(1)},M)},G^{(1)})$ and $(T_{(H^{(1)},M)},H^{(1)})$ are  of (pointwise) weak Lie type, then
$(T_{(G^{(1)},M)}+T_{(H^{(1)},M)},G^{(1)}\cdot H^{(1)})$ is a pair of (pointwise) weak Lie type.
\item
Let $(T_{(G^{(1)},M)},G^{(1)})$ and $(T_{(H^{(1)},M)},H^{(1)})$ be pairs of  pointwise  Lie type. Suppose that the following holds:
for any $z\in M$ and any $u\in T_{(G^{(i)},M)}+T_{(H^{(i)},M)}$ there exist $\xi\in T_{(G^{(i)},M)}$, $\eta\in T_{(H^{(i)},M)}$,
and  $g\in G^{(i)}$, $h\in H^{(i)}$, such that $u=\xi+\eta$ and
\[ord(u\cdot z)<\Big\{ord((g-\one-\xi)z),\ ord((h-\one-\eta)z),\ ord((g-\one-\xi)\eta\cdot z,\ ord(\xi\cdot\eta\cdot z))\Big\}.
\]
Then $(T_{(G^{(1)},M)}+T_{(H^{(1)},M)},G^{(1)}\cdot H^{(1)})$ is a pair of  pointwise  Lie type.
\item
Let $(T_{(G^{(1)},M)},G^{(1)})$ and $(T_{(H^{(1)},M)},H^{(1)})$ be  pairs of   Lie type. Suppose the following holds:
for  any $u\in T_{(G^{(i)},M)}+T_{(H^{(i)},M)}$ there exist $\xi\in T_{(G^{(i)},M)}$, $\eta\in T_{(H^{(i)},M)}$,
and  $g\in G^{(i)}$, $h\in H^{(i)}$, such that $u=\xi+\eta$ and for any $z\in M$ holds:
\[ord(u\cdot z)<\Big\{ord((g-\one-\xi)z),\ ord((h-\one-\eta)z),\ ord((g-\one-\xi)\eta\cdot z,\ ord(\xi\cdot\eta\cdot z))\Big\}.
\]
Then $(T_{(G^{(1)},M)}+T_{(H^{(1)},M)},G^{(1)}\cdot H^{(1)})$ is a pair of    Lie type.
\eee
\eel
Note that the long condition in part {\it 2.} is a weakening of the conditions
\beq
\Big(T_{(G^{(1)},M)}z+T_{(H^{(1)},M)}z\Big)\cap M_j=\Big(T_{(G^{(1)},M)}z\cap M_j\Big)+\Big(T_{(H^{(1)},M)}z\cap M_j\Big),
\eeq
\[\Big(G^{(1)}\cdot H^{(1)}z-\{z\}\Big)\cap M_j=\Big((G^{(1)}z-\{z\})\cap M_j\Big)+\Big(( H^{(1)}z-\{z\}) \cap M_j\Big).
\]
\bpr
\bee[\bf 1.]
\item
 For the pairs of weak Lie type, take $u\in T_{(G^{(1)},M)}+T_{(H^{(1)},M)}$ and choose (any) presentation $u=\xi+\eta$ such that $ord(u)\le ord(\xi), ord(\eta)$.

 By the assumption of weak Lie pairs we fix $g\in G^{(1)}$, $h\in H^{(1)}$ satisfying: $ord(g-\one-\xi)\ge 2\xi$, $ord(h-\one-\eta)\ge 2\eta$.
 Finally we present
 \beq\label{Eq.product.of.groups}
gh-\one-\xi-\eta=(g-\one-\xi)(h-\one-\eta)+(g-\one-\xi)+(h-\one-\eta)+(g-\one-\xi)\eta+\xi(h-\one-\eta)+\xi\cdot \eta.
 \eeq
 Thus we get $ord(gh-\one-\xi-\eta)\ge 2ord(u)$, i.e., the bound of part {\em 3.i.} of definition \ref{Def.Pairs.of.Lie.Type}.

For the bound of part ii. we take any element $\ca\in G^{(1)}\cdot H^{(1)}$ and present it as $\ca=gh$, with $ord(g-\one)\ge ord(\ca-\one)\le ord(h-\one)$.
 By the assumption of weak Lie pairs we fix $\xi,\eta$ satisfying: $ord(g-\one-\xi)\ge 2(g-\one)$, $ord(h-\one-\eta)\ge 2(h-\one)$.
 Now, the bound $ord(gh-\one-\xi-\eta)\ge 2ord(\tg-\one)$ follows from equation \eqref{Eq.product.of.groups}. This verifies condition
  {\em 3.ii.} of definition  \ref{Def.Pairs.of.Lie.Type}.

\smallskip

For the pairs of pointwise weak Lie type the proof is the same, just apply all the formulas to some $z\in M$, and note:
 if $\xi\in T_{(G^{(i)},M)}$, $\eta\in  \xi\in T_{(G^{(i)},M)}$ then $g\cdot h\in G^{(i)}\cdot H^{(i)}$ (and vice versa).
\item Take $z\in M$ and $u\in T_{(G^{(1)},M)}+T_{(H^{(1)},M)}$ and choose (any) presentation $u=\xi+\eta$ satisfying the assumptions.
 Then equation  \eqref{Eq.product.of.groups} ensures $ord((gh-\one-\xi-\eta)z)>ord(uz)$.

Vice versa, take  $\tg\in G^{(1)}\cdot H^{(1)}$ and choose a presentation $\tg=gh$ satisfying the assumptions.
 Then equation  \eqref{Eq.product.of.groups} ensures $ord((gh-\one-\xi-\eta)z)>ord((\tg-\one)z)$.

\item The proof for pairs of Lie type is similar. \epr
\eee

\bex\label{Ex.GL.times.GL.times.Aut.of.weak.Lie.type}
(Continuing part iii.  of section \ref{Sec.Examples.Actions.with.Aut}.) Let $R$ be a filtered ring over a ring $\k$, and $M$ a filtered $R$-module.
 Fix a system of generators  $\{e_i\}$ of $M$, and some
    $R$-multiplicative action $Aut^{(1)}_\k(R)\circlearrowright M$, by $\phi(\sum a_je_j)=\sum\phi(a_j)e_j$.
 Thus  $Aut^{(1)}_\k(R)$ preserves $\{e_j\}$ and we put  $Der^{(1)}_\k(R)(e_j)=0$.
Take some subgroups
 $G\sseteq GL_R(M)$,   $H\sseteq Aut^{(1)}_\k(R)$ such that  the pairs
 $( T_{(G^{(1)},M)},G^{(1)})$, $( T_{(H^{(1)},M)},H^{(1)})$ are  of  weak  Lie type.
 We claim: the pair $(T_{(G^{(1)},M)}+T_{(H^{(1)},M)},G^{(1)}\cdot H^{(1)})$  is of weak Lie type.
\bpr
 The assumptions of part 1 of the last lemma are satisfied. Indeed, suppose $(\xi+\eta)\in (T_{(G^{(1)},M)}+T_{(H^{(1)},M)})\cap End^{(i)}_\k(M)$. Apply this to the generators $\{e_j\}$
  to get: $\xi(e_j)\in M_i$. Thus $\xi\in T_{(G^{(i)},M)}$ and hence $\eta\in T_{(H^{(i)},M)})$.
   Similarly, suppose $g\cdot h\in (G\cdot H)\cap GL^{(i)}_\k(M)$. As $h(e_j)=e_j$, we get: $g(e_j)-e_j\in M_i$, hence $g\in G^{(i)}$. But then $h\in H^{(i)}$.
     \epr

As a particular case take the action of the group  $G=GL(m,R)\times GL(n,R)\rtimes Aut_\k(R)$ on matrices     $M=\Mat$.
 We get: the pair
\beq
\Big(End^{(1)}_R(m)\oplus End^{(1)}_R(n)\oplus Der^{(1)}_\k(R),\ GL^{(1)}(m,R)\times GL^{(1)}(n,R)\rtimes Aut^{(1)}_\k(R)\Big)
\eeq
is of weak Lie type.
\eex

\subsubsection{ The pair $\boldsymbol{(T_{(G^{(1)},M)}+T_{(H^{(1)},M)},G^{(1)}\cdot H^{(1)})}$ in the case $\boldsymbol{\k\supseteq\Q}$}
 While lemma \ref{Thm.Lie.type.of.(H.G)} is rather general, the assumptions of its parts 2 and 3 are often difficult to check.
 However, in the case $\k\supseteq\Q$ we have a much simpler statement
(see Lemma \ref{Thm.(GH).of.pointwise.Lie.type}).

Suppose  $(T_{(G^{(1)},M)},G^{(1)})$  $(T_{(H^{(1)},M)}, H^{(1)})$ are pairs of pointwise Lie type.
 Suppose for any $z\in M$ the correspondences of definition \ref{Def.Pairs.of.Lie.Type},
  $T_{(G^{(1)},M)} \rightleftarrows G^{(1)}$, $T_{(H^{(1)},M)}\rightleftarrows H^{(1)}$ are
 realized by maps
 \beq\label{Eq.exp.log.of.special.type}
 \Psi^{(exp)}(\xi,z)=\sum^N_{j=0}\frac{\xi^j}{j!}+\phi^{(exp)}_{N,z},\quad\quad\quad
  \Psi^{(ln)}(g,z)=\sum^N_{j=1}\frac{(-1)^{j+1}(g-\one)^j}{j}+\phi^{(ln)}_{N,z}.
 \eeq
Here  $\phi^{(exp)}_{N,z},\phi^{(ln)}_{N,z}\in End^{(N+1)}_\k(M)$ are the higher order terms, with $N\ge max\big(ord(\xi z),ord((g-\one)z)\big)$,
 when this $max$ is finite.  (These $\phi^{(exp)}_{N,z},\phi^{(ln)}_{N,z}$ are not necessarily the same for $G,H$.)

For example, this holds for $GL_R(M)$ and in many cases for $Aut_\k(R)$, see the proof of theorem \ref{Thm.Aut(R).pointwise.Lie}.

\bel\label{Thm.(GH).of.pointwise.Lie.type}
Assume $\k\supseteq\Q$ and suppose
$\langle G^{(1)},H^{(1)}\rangle=G^{(1)}\cdot H^{(1)}$. Moreover, assume for the commutator $[T_{(G^{(1)},M)},T_{(H^{(1)},M)}]\sseteq T_{(H^{(1)},M)}\sseteq End^{(1)}_\k(M)$. Then
 the pair $(T_{(G^{(1)},M)}+T_{(H^{(1)},M)},G^{(1)}\cdot H^{(1)})$ is of pointwise Lie type.
\eel
\bpr
\bee[i.]
\item (Verifying condition {\em 3.i.} of definition \ref{Def.Pairs.of.Lie.Type})
 Take an element   $\xi+\eta\in T_{(G^{(1)},M)}+T_{(H^{(1)},M)}$.
 Assume $N:=ord((\xi+\eta)z)<\infty$.
 To define the corresponding element in $G\cdot H$ we recall the Baker-Campbell-Hausdorff formula,
  \beq
  exp(\xi)exp(\eta)=exp\Big(\xi+\eta+\frac{[\xi,\eta]}{2}+\frac{[\xi,[\xi,\eta]]+[\eta,[\eta,\xi]]}{12}+\cdots\Big).
  \eeq
As we do not assume completeness of $M$, we truncate this formula at order $N$ to get:
\beq
\Psi^{(exp)}(\xi,z)\cdot \Psi^{(exp)}(\eta,z)\in  \Psi^{(exp)}(\xi+\eta+\cdots,z)\cdot GL^{(N+1)}_\k(M).
\eeq
 Fix $\tilde\eta\in T_{(H^{(1)},M)}$ by the condition
  $\tilde\eta+\frac{[\xi,\tilde\eta]}{2}+\frac{[\xi,[\xi,\tilde\eta]]+[\tilde\eta,[\tilde\eta,\xi]]}{12}+\cdots\in \{\eta\}+End^{(N+1)}_\k(M)$.
  Such $\tilde\eta$ is obtained iteratively, one puts $\tilde\eta=\eta-\frac{[\xi,\eta]}{2}-\cdots$, then adjusts the higher order terms, and so on.
  Each summand here belongs to $T_{(H^{(1)},M)}$, because $[T_{(G^{(1)},M)},T_{(H^{(1)},M)}]\sseteq T_{(H^{(1)},M)}$.
   And we stop at order $N$.

Finally, associate to $\xi+\eta\in T_{(G^{(1)},M)}+T_{(H^{(1)},M)}$ the element $g\cdot h:=\Psi^{(exp)}(\xi,z)\cdot \Psi^{(exp)}(\tilde\eta,z)\in G^{(1)}\cdot H^{(1)}$.
 By construction we have:
 \beq
 \Big(g\cdot h-\one-(\xi+\eta)\Big)z\in \Big(\{\sum^N_{j=2}\frac{(\xi+\eta)^j}{j!}\}+End^{(N+1)}_\k(M)\Big)z.
 \eeq
Therefore $ord \Big(\big(g\cdot h-\one-(\xi+\eta)\big) z\Big)>ord((\xi+\eta)z)$.
\item (Verifying condition {\em 3.ii.} of definition \ref{Def.Pairs.of.Lie.Type}) For any element of $G^{(1)}\cdot H^{(1)}$ choose a presentation $g\cdot h$.
  Assume $N:=ord((g\cdot h-\one)z)<\infty$.
 Let $\xi:=\Psi^{(ln)}(g,z)$ and $\eta:=\Psi^{(ln)}(h,z)$.
  Define $\tilde\eta=\eta+\frac{[\xi,\tilde\eta]}{2}+\frac{[\xi,[\xi,\tilde\eta]]+[\tilde\eta,[\tilde\eta,\xi]]}{12}+\cdots \in T_{(H^{(1)},M)}$,
   the summation up to order $N$. Then $g\cdot h\in  \Psi^{(exp)}(\xi+\tilde\eta,z)\cdot GL^{(N+1)}_\k(M)$. Associate to $g\cdot h$
   the element $\xi+\tilde\eta\in T_{(G^{(1)},M)}+T_{(H^{(1)},M)}$. By construction holds:
$ord  \big((g\cdot h-\one-\xi-\tilde\eta)z\big)>ord((g\cdot h-\one)z)$.\epr
\eee
\bex\label{Ex.Gl(m,R)GL(n,R)Aut(R).char.0.pointwise.Lie.type}
\bee[i.]
\item Let $\k\supseteq\Q$ and $R$ be one of $\k[[ x ]]/J$, $\k\langle x \rangle/J$, $\k\{x\}/J$, or $C^\infty(\R^p,0)/J$.
 Suppose the filtration satisfies the assumptions of theorem \ref{Thm.Aut(R).pointwise.Lie}.
 Take $H=GL_R(M)$ and $G=Aut_\k(R)$ and suppose $GL_R(M)\rtimes Aut_k(R)$ acts on $M$, see \S\ref{Sec.Examples.Actions.with.Aut}.iii.
We claim: the pair $(End^{(1)}_R(M)+Der^{(1)}_\k(R),GL^{(1)}_R(M)\rtimes Aut^{(1)}_\k(R))$ is  of pointwise Lie type.
\bpr
The pairs $(End^{(1)}_R(M),GL^{(1)}_R(M))$ and $(Der^{(1)}_\k(R),Aut^{(1)}_\k(R))$ are  of pointwise Lie type, with the maps
  $ \Psi^{(exp)}(\xi,z)$, $ \Psi^{(ln)}(\xi,z)$ as in equation \eqref{Eq.exp.log.of.special.type}, see the proof of theorem \ref{Thm.Aut(R).pointwise.Lie}.
   In addition there holds $[Der^{(1)}_\k(R),End^{(1)}_R(M)]\sseteq End^{(1)}_R(M)$.  Indeed, for any $\xi\in Der^{(1)}_\k(R)$, $\eta\in End^{(1)}_R(M)$
    and some generators $\{e_i\}$ of $M$ one has: $[\xi,\eta]\sum_i c_i e_i=\sum_{ij} c_i\xi(\eta_{ij})e_j$, where $\eta_{ij}$ is
     the presentation matrix of $\eta$.
      Now apply the last lemma.
\epr
\item In particular, let $M=R^n$, then the group of contact equivalences, $\cK:=GL(n,R)\rtimes Aut_\k(R)$, gives a pair of pointwise Lie type.
\item In the same way one gets, for $M=\Mat$,  the pair  of pointwise Lie type:
\beq
(End^{(1)}_R(R^m)+End^{(1)}_R(R^n)+Der^{(1)}_\k(R)\ ,\ GL^{(1)}(m,R)\times GL^{(1)}(n,R)\rtimes Aut^{(1)}_\k(R))
\eeq

For square matrices, $m=n$, we get pointwise Lie pairs for the groups $G_{congr}\rtimes Aut_\k(R)$, $G_{conj}\rtimes Aut_\k(R)$, etc,
  see \S\ref{Sec.Examples.Actions.with.Aut}.iv.
 (The verification of $[T_{(G^{(1)},M)},T_{(H^{(1)},M)}]\sseteq T_{(H^{(1)},M)}$ is immediate.)
\eee
\eex
\bex\label{Ex.Gl(m,R)GL(n,R)Aut(R).char0.Lie.type}
 In some cases one can take the full Taylor expansions in equation \eqref{Eq.exp.log.of.special.type}. This holds, e.g. when $M$
 is complete for its filtration. (For example, $R$ is complete wrt. $\{I_j\}$ and $M_j=I_j\cdot M$.)
  Then the maps $ \Psi^{(exp)}(\xi,z)$, $ \Psi^{(ln)}(\xi,z)$ do not depend on $z$, and one gets: the pair
 $(T_{(G^{(1)},M)}+T_{(H^{(1)},M)},G^{(1)}\cdot H^{(1)})$ is of Lie type.
 (In the proof of lemma \ref{Thm.(GH).of.pointwise.Lie.type} one puts $N=\infty$, and all formulas become exact, not just up to $GL^{(N+1)}_\k(M)$, $End^{(N+1)}_\k(M)$.)
\eex
\medskip

\subsubsection{ The case of  direct product, $\boldsymbol{G\times H}$}
We treat this special case separately since it requires different assumptions.

\bel
Fix two pairs of Lie type,  $T_{(G^{(1)},M)}\stackrel{\Psi^{(exp)}_G}{\to }G^{(1)}$, $T_{(H^{(1)},M)}\stackrel{\Psi^{(exp)}_H}{\to }H^{(1)}$.
 Suppose that the maps $\Psi^{(exp)}_G,\Psi^{(exp)}_H$
are power series, with {\em the same coefficients,} i.e.
\[
\Psi^{(exp)}_G(\xi)=\one+\xi+\suml_{j\ge2}a_j\xi^j,\quad\quad\quad  \Psi^{(exp)}_H(\eta)=\one+\eta+\suml_{j\ge2}a_j\eta^j.
\]
Suppose that the same holds also for the maps $\Psi^{(ln)}_G,\Psi^{(ln)}_H$.
Then the pair $\Big(T_{(G^{(1)},M)}\oplus T_{(H^{(1)},M)},G^{(1)}\times H^{(1)}\Big)$ is of Lie type for
the maps
\[
\Psi^{(exp)}_{G\times H}(\xi,\eta)=\Psi^{(exp)}_{G}(\xi)\cdot\Psi^{(exp)}_H(-\eta)^{-1},\quad
\Psi^{(ln)}_{G\times H}(g,h)=\Psi^{(ln)}_{G}(g)-\Psi^{(ln)}_H(h^{-1}).
\]
\eel
(Here $\Psi^{(exp)}_H(-\eta)^{-1}$ is the inverse of $\Psi^{(exp)}_H(-\eta)$ inside $H$.)

\bpr
The map $\Psi^{(exp)}_{G\times H}$ is well defined and one has: $\Psi^{(exp)}_{G\times H}(\xi,\eta)=\one+\xi+\eta+F(\xi,\eta)$,
with $F(\xi,\eta)\in (\xi^2,\eta^2,\xi\cdot\eta)$.

Moreover, $\Psi^{(exp)}_{G\times H}(\xi,-\xi)=\one$, therefore
$F(\xi,\eta)=(\xi+\eta)\tilde{F}(\xi,\eta)$, where $\tilde{F}(\xi,\eta)\in (\xi,\eta)$. Therefore, for any $z\in M$ holds:
$ord\left(F(\xi,\eta)\cdot z\right)>ord(\xi\cdot z+\eta\cdot z)$.

Similarly, for the logarithmic map we have:
\beq
\Psi^{(ln)}_{G\times H}(g,h)=(gh-\one)+F_G(g-\one)-F_H(h^{-1}-\one)-(g-\one)(h-\one)-h^{-1}(h-\one)^2, \quad \Psi^{(ln)}_{G\times H}(g,g^{-1})=\zero.
\eeq
Thus the total higher order expression, $F_G(g-\one)-F_H(h^{-1}-\one)-(g-\one)(h-\one)-h^{-1}(h-\one)^2$, is presentable in the form $(gh-\one)\cdot\tilde{F}(g,h)$, where
$\tilde{F}(g,h)\in (g-\one)+(h-\one)$. In this way we get the needed bounds of Definition \ref{Def.Pairs.of.Lie.Type}.
\epr
\bex
We continue Example \ref{Ex.Lie.type.pair.GLM}.
Let $(R,\cm)$ be a local ring and $GL(n,R)\circlearrowright R^n$.
The pair $(Mat_{n\times n}(\cm),GL^{(1)}(n,R))$ is of Lie type for the map $A\to\one+A$. Therefore the pair
\[
Mat_{m\times m}(\cm)\oplus Mat_{n\times n}(\cm)\to GL^{(1)}(m,R)\times GL^{(1)}(n,R)
\]
is of Lie type for $\Psi^{(exp)}(A,B)=(\one+A)(\one-B)^{-1}$. (And not just of weak Lie type, as was proved in Example \ref{Ex.GL.times.GL.times.Aut.of.weak.Lie.type}).
\eex

\section{The general criteria of determinacy} \label{Sec.Finite.Determinacy.Criteria}
Fix a filtered action $G\circlearrowright M$ and a pair of  pointwise  (weak) Lie
type, $(T_{(G^{(1)},M)},G^{(1)})$.
By definition, $T_{(G^{(1)},M)}\sseteq End^{(1)}_\k(M)$, this defines the action  $T_{(G^{(1)},M)}\circlearrowright M$.
For any $z\in M$ we take  the corresponding orbits   $T_{(G^{(1)},M)}(z)$ resp. $G^{(1)}z\sseteq M$ and their closures $\overline{T_{(G^{(1)},M)}(z)}$ resp.
 $\overline{G^{(1)}z}$, in the filtration topology.
  We prove the statements of the type ``$\overline{T_{(G^{(1)},M)}(z)}$ is large iff  $\overline{G^{(1)}z}$ is large".

Recall that  $z\in M$ is order-by-order $N$-determined iff $\{z\}+M_{N+1}\sseteq \overline{G^{(1)}z}$. Thus the statements read roughly
 ``finite determinacy is equivalent to large tangent space".

\subsection{$\overline{T_{(G^{(1)},M)}(z)}$ vs  $\overline{G^{(1)}z}$ for pairs of pointwise Lie type}\label{Sec.Finite.Determinacy.Poinwise.Lie.type}
\bthe\label{Thm.Finite.Determinacy.pointwise.Lie.type}
 Suppose $(T_{(G^{(i)},M)},G^{(i)})$ is pointwise of Lie type, for some $i\ge1$
 (i.e., the properties in Definition \ref{Def.Pairs.of.Lie.Type} hold  for $G^{(i)}$ instead of $G^{(1)}$).
 Then for any $z\in M$ holds:
\[
\quad\quad\quad\overline{T_{(G^{(i)},M)}z}\supseteq M_{N+1}\quad\quad \text{ iff } \quad\quad
\overline{G^{(i)}z}\supseteq \{z\}+M_{N+1}.
\]
\ethe
Thus, $z$ is order-by-order  $N$-determined iff it is ``infinitesimally" order-by-order  $N$-determined.

\bpr
 $"\Rightarrow":$
 Take any  $w_{N+1}\in M_{N+1}$. We construct a sequence $\{g_n\in G^{(i)}\}$ satisfying: $g_n^{-1}\cdots g_1^{-1}(z+w_{N+1})\in \{z\}+M_{N+1+n}$.

  By the assumption $w_{N+1}\in \xi\cdot z+M_{N+2}$ for some $\xi\in T_{(G^{(i)},M)}$.
  Since  $(T_{(G^{(i)},M)},G^{(i)})$ is pointwise of Lie type there is   $g_1\in G^{(i)}$,
 satisfying $ord((g_1-\one-\xi)z)>ord(\xi\cdot z)$. Then $g_1z-z-w_{N+1}\in M_{N+2}$, i.e., $g^{-1}_1(z+w_{N+1})=z+w_{N+2}$, for some $w_{N+2}\in M_{N+2} \subseteq M_{N+1} $.
 By  assumption $w_{N+2}\in \xi\cdot z+M_{N+3}$ for some $\xi\in T_{(G^{(i)},M)}$ and there is a $g_2\in G^{(i)}$,
 satisfying $g_2z-z-w_{N+2}\in M_{N+3}$,
  i.e., $g^{-1}_2(z+w_{N+2})=z+w_{N+3}$, $w_{N+3}\in M_{N+3}.$
Proceed inductively.

 $"\Leftarrow":$ Let $w_{N+1}\in M_{N+1}$, we construct a sequence $\{\xi_n\in T_{(G^{(i)},M)}\}$ satisfying: $w_{N+1}-(\sum^n_1 \xi_j) z\in M_{N+1+n}$.
 By the assumption,   $gz-z-w_{N+1}\in M_{N+2}$ for some  $g\in G^{(i)}$.
Agian,  since  $(T_{(G^{(i)},M)},G^{(i)})$ is pointwise of Lie type there is
   $\xi_1\in T_{(G^{(i)},M)}$, satisfying: $ord((g-\one-\xi_1)z)>ord((g-\one)z)$. Thus $w_{N+1}\in \xi_1\cdot z+M_{N+2}$, i.e.
  $w_{N+1}- \xi_1\cdot z=w_{N+2}\in M_{N+2}$ and we can proceed inductively.
  \epr

\beR
One is naturally tempted to a stronger statement, ``Suppose $(T_{(G^{(1)},M)},G^{(1)})$ is of Lie  type, then
$\{z\}+ T_{(G^{(k)},M)}(z) \sseteq \overline{G^{(k)}z}$ for some $k\ge1$."
This does not hold, see \cite[Remark 2.3]{B.K.motor}.
\eeR

\subsection{$\overline{T_{(G^{(1)},M)}(z)}$ vs  $\overline{G^{(1)}z}$ for pairs of pointwise weak  Lie type}\label{Sec.Finite.Determinacy.Poinwise.weak.Lie.type}
\bthe\label{Thm.Finite.Determinacy.weak.Lie.pairs}
  Suppose $(T_{(G^{(1)},M)},G^{(1)})$ is a pair of pointwise weak Lie type. Then for any $z\in M$ holds:
\bee[1.]
\item
If $\overline{T_{(G^{(k)},M)}z}\supseteq M_{N+k+ord(z)}$ for any $k\ge1$ then $\overline{G^{(k)}z}\supseteq \{z\}+M_{N+k+ord(z)}$ for any $k>N$.
\item
If $\overline{G^{(k)}z}\supseteq \{z\}+M_{N+k +ord(z)}$  for any $k\ge1$ then $\overline{T_{(G^{(k)},M)}z}\supseteq M_{N+k+ord(z)}$ for any $k>N$.
\eee
\ethe
In the first part we get:  $z$ is  order-by-order $(2N+ord(z))$-determined.

\bpr
\bee[\bf 1.]
\item
For any $w\in M_{N+k+ord(z)}$, with $k>N$, we construct a sequence $\{g_n\in G^{(k)}\}$ satisfying: $g_n^{-1}\cdots g_1^{-1}(z+w)\in \{z\}+M_{N+k+n+ord(z)}$.

By assumption we have $w\in \xi\cdot z+M_{N+k+1+ord(z)}$, for some $\xi\in T_{(G^{(k)},M)}$. By the pointwise weak Lie type assumption there exists
 $g_1\in G^{(k)}$ satisfying:  $(g_1-\one-\xi)z\in M_{2ord(\xi)+ord(z)},  \
 ord(\xi)\geq k \geq N+1$.
 Thus $g_1 z-z-w\in  M_{N+k+1+ord(z)}$.
 Hence $g_1^{-1}(z+w) = z + w_1$ with $w_1 \in M_{N+k+1+ord(z)}$. We continue now with $w_1$
and with $k$ replaced by $k+1$. Then we have $w' \in \xi\cdot z+M_{N+k+2+ord(z)}$, for some $\xi\in T_{(G^{(k+1)},M)}$, i.e. $ord(\xi) \geq k+1$.
The pointwise weak Lie type assumption implies that there exists
 $g_2\in G^{(k+1)} \subseteq G^{(k)}$ satisfying:  $(g_2-\one-\xi)z\in M_{2(k+1)+ord(z)}$.
 Thus $g_2 z-z-w'\in  M_{N+k+2+ord(z)}$ and $g_2^{-1}(z+w') = z + w_2$ with $w_2 \in M_{N+k+2+ord(z)}$.
 Proceed inductively.

 \item
 For any $w\in M_{N+k+ord(z)}$, with $k>N$, we construct a sequence $\{\xi_n\in T_{(G^{(k)},M)}\}$ satisfying: $\sum^n_1 \xi_j(z)-w\in M_{N+k+n+ord(z)}$.

 By the assumption,   $gz-z-w\in M_{N+k+1+ord(z)}$ for some  $g\in G^{(k)}$. Take the corresponding
  $\xi\in T_{(G^{(k)},M)}$, satisfying: $ord((g-\one-\xi)z)\ge 2k+ord(z)$. Thus $w-\xi\cdot z\in M_{N+k+1+ord(z)}$. Similar as in 1. we proceed inductively. \epr
\eee

\bcor\label{corollary of finite determinacy} Suppose $(T_{(G^{(1)},M)},G^{(1)})$ is a pair of pointwise weak Lie type.
\bee[1.]
\item   Assume that $M_{N+1+ord(z)}\sseteq \overline{T_{(G^{(1)},M)}(z)}$  implies
$M_{N+k+ord(z)}\sseteq \overline{T_{(G^{(k)},M)}(z)}$ for all $k\ge1$. If  $M_{N+1}\sseteq \overline{T_{(G^{(1)},M)}(z)}$
then $\{z\}+M_{2N+1+ord(z)}\sseteq \overline{G^{(N+1+ord(z))}z}.$

\item Assume that $\{z\}+M_{N+1+ord(z)}\sseteq \overline{G^{(1)}z}$ implies $\{z\}+M_{N+k+ord(z)}\sseteq \overline{G^{(k)}z}$ for all $k\ge1$.
 If $\{z\}+M_{N+k}\sseteq \overline{G^{(k)}z}$ then
$\overline{T_{(G^{(N+1)},M)}(z)}\supseteq M_{2N+ord(z)+1}$.
\eee
\ecor

The assumption of part 1. holds for many pairs of pointwise weak Lie type, see  the following Example \ref{Ex.for.cond.50} i.. The assumption of part 2. is more subtle.

\bex\label{Ex.for.cond.50}
\bee[i.]
\item Suppose a ring $R$ is filtered by ideals $\{I_j\}$, satisfying the condition $(I_{i+k}:I_i)\cdot I_j\supseteq I_{k+j}$ for any $i,j,k\ge1$.
 (The filtrations $\{I^j\}$, $\{\ca^j\cdot I_1\}$ satisfy this condition.)
 Let $M$ be $R$-module with the induced filtration $M_j=I_j\cdot M$.
 Suppose $(T_{(G^{(1)},M)},G^{(1)})$ is a pointwise weak Lie pair and $M_{N+1+ord(z)}\sseteq \overline{T_{(G^{(i+1)},M)}(z)}$,  for some $i\ge0$. Then
  $ M_{N+k+ord(z)}\sseteq \overline{T_{(G^{(i+k)},M)}(z)}$ for any $k\ge1$. Indeed, we observe:
\beq
T_{(G^{(i+k)},M)}=T_{(G^{(i+1)},M)}\cap End^{(i+k)}_\k(M)\supseteq (I_{i+k}:I_{i+1})\cdot T_{(G^{(i+1)},M)}.
\eeq
 Here the embedding holds because $(I_{i+k}:I_{i+1})\cdot End^{(i+1)}_\k(M)\sseteq End^{(i+k)}_\k(M)$. Thus
\beq
\overline{T_{(G^{(i+k)},M)}(z)}\supseteq \overline{(I_{i+k}:I_{i+1})\cdot   T_{(G^{(i+1)},M)}(z)}\supseteq
 (I_{i+k}:I_{i+1})\cdot  I_{N+1+ord(z)}\cdot M\supseteq I_{N+k+ord(z)}\cdot M.
\eeq
\item If the filtration of $M$ does not satisfy some condition similar to $(I_{i+k}:I_i)\cdot I_j\supseteq I_{k+j}$, as above, then the assumptions of
 the last corollary do not hold. For example, let $M=R=\k[[x]]$  filtered by
 \beq
 I_1=\cm^n\supset I_2=\cm^{2n}\supset I_3=\cm^{2n+1}\supset I_4=\cm^{2n+2}\supset\cdots.
 \eeq
Let $f(x)=x^{n+1}$ and assume $char(\k)\neq (n+1)$. Then $T_{(Aut^{(1)}_\k(R),M)}(f)=(x)^{2n}=I_1$,  but  we have
$T_{(Aut^{(2)}_\k(R),M)}(f)=(x)^{3n}=I_{n+2}$.
 Similarly $Aut^{(1)}_\k(R)(f)\supseteq \{f\}+(x)^{2n}$ but $Aut^{(2)}_\k(R)(f)\not\supseteq \{f\}+(x)^{2n+1}$.
\eee
\eex

\subsection{Sharpness of results}\label{Sec.Finite.Determinacy.Sharpness.of.bounds}
In Theorems \ref{Thm.Finite.Determinacy.pointwise.Lie.type}, \ref{Thm.Finite.Determinacy.weak.Lie.pairs} we see essential difference between
the Lie type case (typically when $\k$ is a field of characteristic zero) and the weak Lie  type case (e.g., when $\k$ is of positive characteristic).
The natural question is whether the bounds in the weak Lie case
can be improved, brought closer to
the bounds for the Lie-type case.
The following example shows that the bound of Theorem \ref{Thm.Finite.Determinacy.weak.Lie.pairs} in the weak Lie-type case is close to being sharp.
\
Let $\k$ a field of characteristic $p>0$.  Take $R=\k[[x]]$, filtered by $\{(x)^j\}$,   and $G=Aut_\k(R)$.
  Then $T_{(G^{(1)},M)}=Der^{(1)}_\k(R)=(x^2\di_x)$.
Let $f=x^p+x^{p+N}$, with $N>2p$, $gcd(p,N)=1$. Then $T_{(G^{(1)},M)}(f)=(x^{p+N+1})$, so one would naively expect the order of determinacy
to be close to $(p+N)$, which is however not the case.

In fact, for  any univariate $f \in \k[[x]]$ the order of determinacy is computed in \cite [Proposition 2.8]{Nguyen}, assuming $\k=\bar\k$, and for this $f$ the exact order of determinacy is $p+N+\lceil\frac{N}{(p-1)}\rceil -1$.
%
In particular, for $p=2$ we get: the order of determinacy equals $p+2N=2(p+N)-ord(f)$, which has the same order as the bound
of part 1. of Theorem \ref{Thm.Finite.Determinacy.weak.Lie.pairs}.

For a multivariate example (and contact equivalence) see Remark 2.3 (b) of \cite {Boubakri.Gre.Mark}.

%
%

\medskip

The following examples illustrate the non-triviality of the prime characteristic case.
\bex
\bee[i.]
\item
Let $\k$ be an algebraically closed field, $char(\k)=p$. Let $R=\k[x]/(x^{pN+1})$, with $N>1$. Then $f(x)=x^p\in R$ is clearly $(pN-1)$-determined.
But $T_{(Aut^{(1)}_\k(R),M)}(f)=\{0\}$.
\item Over a field of zero characteristic we often have: if $\{z\}+(M_N\smin M_{N+1})\sseteq\overline{G^{(1)}z}$ then $\{z\}+M_N\sseteq\overline{G^{(1)}z}$.
This does not hold in prime characteristic. For example, let $\k$ be algebraically closed, $char(\k)=p$, $R=\k[[x]]$ and $f(x)=x^p\in R$.
Then $\{f\}+(x^{pn})\smin(x^{pn+1})\sseteq Aut^{(1)}_\k(R)\cdot f$. But $(f)+(x^{pn+1})\not\sseteq Aut^{(1)}_\k(R)\cdot f$.
\eee
\eex
\smallskip

\subsection{Finite determinacy in terms of infinitesimal stability}\label{Sec.Finite.Determinacy.vs.Infinitesimal.Stability}
In this subsection we assume $R$ to be a   local Noetherian ring, filtered by ideals $\{I_j\}$, and $M$ a  finitely generated $R$-module,
 filtered by $\{I_j\cdot M \}$.
Theorem \ref{Thm.Finite.Determinacy.pointwise.Lie.type} and Nakayama lemma  give the following corollary.
\bcor
Let   $(T_{(G^{(1)},M)},G^{(1)})$ be a pair of pointwise  Lie type.
Then  $z\in M$ is $G^{(1)}$-finitely order-by-order determined if and only if the quotient module $M/(T_{(G^{(1)},M)}(z))$ is annihilated by some $I_n$, for $n\gg1$.
\ecor

Geometrically this means: the module  $M/(T_{(G^{(1)},M)}(z))$ is not supported outside of  $V(I_n)\sset Spec(R)$.

\medskip

In many cases the elements of the module $M$ may be considered as functions on the scheme $Spec(R)$, e.g., this happens for
$M=Maps\big((\k^n,0),(\k^m,0)\big)$ or $M=\Mat$. Then we can
evaluate $T_{(G^{(1)},M)}(z)$ at points of $Spec(R)$ and compare it with
the ambient module.  More precisely, for any prime ideal $\cp\in Spec(R)$
  we take the generic point of  the corresponding (irreducible) locus, $V(\cp)\sset Spec(R)$,
   i.e., pass to the field of fractions, $Frac(R/\cp)$. Accordingly, we pass from modules over $R$ to  vector spaces over $Frac(R/\cp)$,
 \beq
T_{(G^{(1)},M)}(z)\otimes_R Frac(R/\cp)\sseteq M\otimes_R Frac(R/\cp).
 \eeq
Then the condition ``the module  $M/(T_{(G^{(1)},M)}(z))$ is not supported outside of
$V(I_n)$'' can be stated as:
\beq
{\rm if} \ \cp\not\supseteq I_n\ {\rm then} \ T_{(G^{(1)},M)}(z)\otimes Frac(R/\cp)= M\otimes Frac(R/\cp).
\eeq
Geometrically this means: the condition  $T_{(G^{(1)},M)}(z)|_{pt}=M|_{pt}$ holds  for any point $pt\in Spec(R)\smin V(I_n)$.
In the classical terminology this equality of vector spaces is called ``infinitesimal stability at a given point''.
 Therefore in the classical language we get:
\bcor
With the assumptions as before,  $z\in M$ is $G^{(1)}$-finitely order-by-order
determined if and only if $z$ is infinitesimally
stable  at any closed point of $Spec(R)\smin V(I_n)$.
\ecor
For the rings like $\C\{ x\}$ or $\R\{ x\}$, and the classical groups like right or contact equivalence,
this recovers the classically known criteria, see e.g., \cite[Theorem
2.1]{Wall-1981}.

\subsection{The passage from $\overline{Gz}$ to $Gz$ in the $C^\infty$-case.} \label{Sec.Approximation.Cinfty.case}
The results of sections \ref{Sec.Finite.Determinacy.Poinwise.Lie.type}-\ref{Sec.Finite.Determinacy.vs.Infinitesimal.Stability} are of order-by-order type,  $\overline{G^{(1)}z} \supseteq\{z\}+M_{N+1}$, i.e., they address the orbit closures.
To get the ordinary determinacy criteria, i.e., of type  $G^{(1)}z\supseteq\{z\}+M_{N+1}$,
 we use the approximation results of section \ref{Sec.From.obo.determinacy.to.G.determinacy}.

 In  Proposition \ref{Thm.Use.of.Approximation} we did not consider $C^\infty$-case, as this ring has no Artin approximation.
 Now, with the notion of Lie type pairs, we can extend the approximation results to the $C^\infty$-case, and we obtain the
  usual (not just order-by-order) determinacy statements.

Let   $R=C^\infty(\R^p,0)/J$ with the maximal ideal $\cm = \langle x_1,...,x_n \rangle$.
  Let $M$ be a finitely presented filtered $R$-module.
\bprop  \label{Thm.Approximation.for.C.infty.case}
 Suppose the ideal $J\sset C^\infty(\R^p,0)$ is generated by analytic power series. Suppose the module
    $M$ admits a presentation matrix with analytic entries and satisfies $M_{N+1}\supseteq \cm^\tN\cdot M$, for some $\tN\gg N$.
    Suppose the group $G \sseteq GL_R(M) \rtimes Aut_\R (R)$  is defined by analytic equations and
      $(T_{(G^{(1)},M)},G^{(1)})$ is a pair of pointwise Lie type.
      If $\overline{  G^{(1)}z}\supseteq \{z\}+M_{N+1}$ or $\overline{T_{(G^{(1)},M)}z}\supseteq M_{N+1}$  then $ G^{(1)}z\supseteq \{z\}+M_{N+1}$.
\eprop
\bpr
 (The proof below is based on the initial proof of G. Belitskii \cite{Belickii}.)

Fix some $w\in M_{N+1}$.
 By theorem \ref{Thm.Finite.Determinacy.pointwise.Lie.type}
 $\overline{T_{(G^{(1)},M)}z}\supseteq M_{N+1}$ implies
 $\overline{  G^{(1)}z}\supseteq \{z\}+M_{N+1}$, meaning that
the condition $gz=z+w$, $g\in G$ is resolvable order-by-order. We prove that an order-by-order solution implies an ordinary solution.
 First we show a formal solution of the completion, $\hg\hz=\hz+\hw$. Then we check the surjectivity of completion map , $G\twoheadrightarrow\hG$.
 Finally, we establish a solution $g\in G$.

 We remark that  $\overline{T_{(G^{(1)},M)}z}\supseteq M_{N+1}$ implies $T_{(G^{(1)},M)}z\supseteq M_{N+1}$.
 (Apply Nakayama and use $M_{N+1}\supseteq \cm^\tN\cdot M$.)

\bee[\bf Step 1.]
\item
Take the $\cm$-adic completion $R\to \hR=\R[[x]]/J$. 
 In this case $R$ is not Noetherian and $R \to \hR $ is surjective,   by Borel's lemma.
The kernel of $R \twoheadrightarrow \widehat R$ is $\cm^\infty$, the ideal of flat functions.
 Accordingly we have the maps
\beq
R^n\to \hR^n, \quad\quad\quad
M\to \hM \cong M\otimes_R \hR, \quad\quad\quad
G\to \hG\sseteq GL_\hR(\hM)\rtimes Aut_\k(\hR).
\eeq
As in the proof of proposition \ref{Thm.Use.of.Approximation} we lift the condition $gz=z+w$, $g\in G$ to the presentation of $M$,
to get implicit function equations, \eqref{CondEqu2}. By our assumption they admit order-by-order solution, thus their completion admits a formal solution.
 Therefore we have:
$\hG\hz\supseteq \{\hz\}+\widehat{M_{N+1}}$. (Here $\widehat{M_{N+1}}$ is the image of $M_{N+1}$ in $\hM$.)

\item  We verify the surjectivity of the completion map $G\twoheadrightarrow \hG$.
Using the diagram \eqref{ComDiag} any  $\hg\in \hG$ is presentable as   $(\hU,\hV,\hat\phi)\in GL_\hR(\hR^p)\rtimes GL_\hR(\hR^q)\rtimes Aut_\R(\hR)$.
 To find a corresponding preimage $g\in G$ we should find $(U,V,\phi) \in  GL_R(R^p)\rtimes GL_R(R^q)\rtimes Aut_\k(R)$
  that satisfy   $UA=\phi(A)V$, $\phi(A) = [a_{i,j}(\phi)]$ (and further equations of $G$) and whose Taylor expansions are $(\hU,\hV,\hat\phi)$.

As before, we get the system of implicit function equations, \eqref{CondEqu2}.
By the assumption $J$ is generated by analytic series (hence finitely generated), the entries of $A$ are analytic
  and the equations of $G\sseteq GL_R(M)\rtimes Aut_\R(R)$ are analytic. And we are given a formal solution $(\hU,\hV,\hat\phi)$.

  Therefore,   Tougeron  approximation (theorem \ref{Thm.Approximation.Tougeron})
   implies the existence of $C^{\infty}$-solution, $(U,V,\phi)$, whose Taylor expansion is precisely $(\hU,\hV,\hat\phi)$.
 This is the needed preimage $g\in G$.

By restricting to the unipotent part we get the surjection $G^{(1)}\twoheadrightarrow \hG^{(1)}=\widehat{G^{(1)}}$.

\item
For $z\in M$ and any $w\in M_{N+1}$ we want to resolve the condition $gz=z+w$, $g\in G$.
By Step 1 this has a formal solution, $\hg \hz=\hz+\hw$, with $\hg\in \hG$.
By Step 2, there is a preimage $g\in G$ of $\hg$,
  and it satisfies: $gz-z-w\in \cap \al(\cm^j\cdot R^p)=\al(\cm^\infty\cdot R^p)= \cm^\infty\cdot M$, for $\al$  as in \S\ref{Sec.Condition.z+w=gz.As.implicit.function.eq}.
  Therefore it is enough to prove: $Gz\supseteq \{z\}+\cm^\infty\cdot M$.

\

Fix $z$ and vary $g,w$ then we get the (lifted) orbit $\{(g,w)|\ gz=z+w\}\sset G\times M$. We recall its defining equations,
 as in \S\ref{Sec.Condition.z+w=gz.As.implicit.function.eq}.
Fix the canonical basis $\{e_i\}$ of $R^p$, thus $z=\sum z_i\al(e_i)$, $w=\sum w_i\al(e_i)$, here $z_i\in R$. A group element is $g=(U,V,\phi)\in G$
 and the explicit definition  of this orbit   is:
 \beq
 \{(U,V,\phi,\{w_i\},\{t_j\})|\ \sum_j U_{ij}z_j(\phi(x))=z_i+w_i+\sum_j V_{ij}t_j,\ \forall i\}.
 \eeq
 These are $C^\infty$-implicit function equations, $F(y,x,w)=0$. (Here $y$ is the tuple of  the unknowns $U,V,\phi,\{t_{ij}\}$.)
  In our case $w_i\in \cm^\infty$ thus the completion of these equations has the trivial solution, $y_0=(Id,Id,x,0)$. To ensure an ordinary $C^\infty$
   solution we use theorem \ref{Thm.Approximation.Cinfty.equations}.

In our case $\di_y F|_{y_0}$ is exactly the tangent space map, $T_{(G^{(1)},M)}\to T_{(G^{(1)},M)}z\sseteq M$, $\xi\to \xi(z)$. By the assumption, the
 image of this map contains $\cm^\tN\cdot M$. Thus the ideal of maximal minors of  $\di_y F|_{y_0}$ contains $\cm^{\tilde{N}}$ for some $N<\tilde{N}<\infty$.

Finally we estimate $\det\Big(F'_y(x,y_0)(F'_y(x,y_0))^t\Big)$.
\bee[i.]
\item
Recall the Cauchy-Binet formula: if $L\in \Mat$, $m\le n$, then $det(LL^t)=\sum det(L_\Box)^2$, where $L_\Box$ is a maximal ($m\times m$) minor and the sum goes
 over all the maximal minors.
\item
Recall the statement:  if the ideal generated by $\{f_i\}$ contains $\cm^N\sset C^\infty(\R^p,0)$, for some $N<\infty$,
  then $(\sum f^2_i)\cdot \cm^\infty=\cm^\infty$. Indeed, for any $g\in \cm^\infty$ the vanishing order of $\frac{g}{\sum f^2_i}$ at
  $0\in \R^p$ is $\infty$.
   Therefore this fraction extends to a $C^\infty$ function on $(\R^p,0)$.
\eee

Altogether we get  $\det\Big(F'_y(x,y_0)(F'_y(x,y_0))^t\Big)\cm^\infty=\cm^\infty$, thus theorem \ref{Thm.Approximation.Cinfty.equations} ensures the $C^\infty$ solution.
 This resolves the condition $gz=z+w_\infty$, $g\in G$.\epr
\eee

\beR
Proposition \ref{Thm.Approximation.for.C.infty.case} does not hold if the ideal $J\sset C^\infty(\R^p,0)$ is not generated by analytic functions.
For example, let $R=C^\infty(\R^2,0)/(x_2)^\infty$, then $\hR=\R[[x_1,x_2]]$ and the projection $Aut(R)\to Aut(\hR)$  is non surjective.
 The order of $Aut_\R(\hR)$-determinacy of the element $\hx_2\in \R[[x_1,x_2]]$ is 1, but the element
  $x_2\in C^\infty(\R^2,0)/(x_2)^\infty$ is not finitely $Aut_\R(R)$-determined. For example, $x_2\not\sim x_2+x^N_1$ for any $N$.

In addition, if $J$ is not generated by analytic functions then it might be not finitely generated, as $C^\infty(\R^p,0)$ is not Noetherian.
Then condition $\phi(J) = J$ brings  infinitely many  equations.
\eeR

\subsection{}
Finally we state the criterion for the orbits, not just their closures.
  Suppose the group $G\sseteq GL_R(M)\rtimes Aut_\k(R)$    acts on $M$.
   We assume that the filtration satisfies $Der^{(1)}_\k(R)(x)\sseteq(x)^2$ and  $I_N\cdot Der_\k(R)\sseteq Der^{(1)}_\k(R)$ for $N\gg1$.

 As a consequence of the previous results we state the following finite determinacy result which covers positive characteristic, characterisitc 0 and the $C^\infty$-case (see Remark \ref{Thm.Approximation.W.systems} and section \ref{Sec.C.infty}).

\bcor
\bee[1.]
\item
Let $\k$ be any field, and let $R$ be one of $\k[[x]]/J, \k \langle x \rangle/J$, or
$\k\{x\}/J$ (for the latter case assume that $\k$ is a valued field such that its completion with respect to the absolute value is separable over $\k$).
 Suppose the conditions \ref{Assumptions.on.the.ring} are satisfied.
   Suppose $G\sseteq GL_R(M)\rtimes Aut_\k(R)$
  is defined by formal/algebraic/analytic power series.
   If $T_{(G^{(1)},M)}z\supseteq  M_{N+1}$ then
$G^{(k)}z\supseteq \{z\}+M_{N+k}$, for any $k>max(0,N-ord(z))$.
\item Suppose $\k\supseteq\Q$ and  $R$ is one of $\k[[x]]/J$, $\k\langle x\rangle /J$, $\k\{x\}/J$.
 Suppose $G\sseteq GL_R(M)\rtimes Aut_\k(R)$  is defined by formal/algebraic/analytic equations.
 Then   $T_{(G^{(1)},M)}z\supseteq M_{N+1}$  iff $G^{(1)}z\supseteq \{z\}+ M_{N+1}$.
\item
Suppose $R=C^\infty(\R^p,0)/J$ with $J$ analytically generated, $M$ admits a presentation matrix with analytic entries,
 and $G\sseteq GL_R(M)\rtimes Aut_\k(R)$  is defined by  analytic equations.
  Suppose $M_{N+1}\supseteq \cm^\tN M$ for some $\tN\gg N$.
   Then   $T_{(G^{(1)},M)}z\supseteq M_{N+1}$  iff $G^{(1)}z\supseteq \{z\}+ M_{N+1}$.
\eee
\ecor
\bpr
 By theorem  \ref{Thm.Aut(R).pointwise.Lie}  and proposition \ref{Thm.Aut(R).of.weak.Lie.type} we have the (pointwise) (weak) Lie type pairs.

  For part 1 use theorem \ref{Thm.Finite.Determinacy.weak.Lie.pairs}
 to get $\overline{G^{(k)}z}\supseteq \{z\}+\cm^{N+k}\cdot M$.
For parts 2,3 use  theorem  \ref{Thm.Finite.Determinacy.pointwise.Lie.type}, to get:
\beq
\overline{T_{(G^{(1)},M)}z}\supseteq M_{N+1}\quad \text{ if and only if }\quad \overline{G^{(1)}z}\supseteq \{z\}+ M_{N+1}.
\eeq

 Then use propositions \ref{Thm.Use.of.Approximation} and \ref{Thm.Approximation.for.C.infty.case} to remove the closures.
\epr

\section{Applications and examples}\label{Sec.Applications.Examples}
In this section we derive several explicit determinacy statements. We work under the following assumptions:
\begin{Assumptions}\label{Assumptions.on.the.ring.for.applications}

\bee[i.]
\item Let $\k$  be any field and, when $R$ is $\k\{x\}/J$, we assume that $\k$  is a complete valued field (of any characteristic).

\item $R\sseteq \k[[x]]/J$ or $R=C^\infty(\R^p,0)/J$, with filtration by some ideals $\{I_j\}$, satisfying:
\[
Der^{(1)}_\k(R)(x)\sseteq (x)^2, \quad \text{ and}  \quad
 I_{N+1}\cdot Der_\k(R)\sseteq Der^{(1)}_\k(R)\quad \text{ for}  \quad N\gg1.
\]
\item
\bei
\item For $R=C^\infty(\R^p,0)/J$ we assume $J$ is analytically generated,  $\cap I_j\supseteq\cm^\infty$, and $\{I_j\}$ are analytically generated $mod\ \cm^\infty$.
\item If $R$ is not one of  the rings $\k[[x]]/J$, $\k\langle x\rangle /J$, $\k\{x\}/J$, $C^\infty(\R^p,0)/J$, then we also assume that
 $R$ has the Artin approximation property;
  $\{I_j\},J$ are polynomially generated and the objects for which we derive the determinacy (e.g. $f\in R$ or $f\in R^n$ or
 $f\in \Mat$) is polynomial $mod\ I_N$, for some (specified) $N$.
\item If $\k\not\supseteq\Q$ then  in addition we assume conditions \ref{Assumptions.on.the.ring}.
\eei
\eee
\end{Assumptions}
These conditions are needed to get the (pointwise resp. weak) Lie type pair and to apply the approximation results.

Recall that the filtrations $\{\cm^j\}$, $ \{I_{j+1}=I^j_1\cdot \ca\}$, for some ideals $I_1\sseteq \ca^2$  satisfy these assumptions.

Recall that $ord(f)$ is defined in \S\ref{Sec.Preparations.Induced.Filtration}.

For a group action $G\circlearrowright M$ and a filtration $\{M_j\}$ of $M$ we say that $z\in M$
  is infinitesimally $d$-determined if $T_{(G^{(1)},M)}z\supseteq M_{d+1}$.

\subsection{Right determinacy of germs of functions}

We formulate our results for the classical case of function germs and show how
they specialize to classically known results on finite right-determinacy (our results are more general with respect to the allowed rings and filtrations).

\bcor \label{cor5.1} Suppose $R$ with filtration $\{I_j\}$  satisfy the Assumptions \ref{Assumptions.on.the.ring.for.applications}. For $f\in R$ the following holds:
\bee
\item [\bf 1.]
\bee[i.]
\item   If $Der^{(1)}_\k(R)f\supseteq I_{N+1}$ then $\{f\}+I_{2N+1-ord(f)}\sseteq Aut^{(N+1-ordf)}_\k(R)\cdot f$.

(Thus  $f$ is $(2N-ord(f))$-determined.)
\item  [ii.]
If $\{f\}+I_{N+k}\sseteq Aut_{\k}^{(k)} (R)\cdot f$ for any $k\ge1$, then  $I_{2N+1-ord(f)}\sseteq Der^{(1)}_\k(R)(f).$

(Thus  $f$ is infinitesimally $(2N-ord(f))$-determined.)
\eee
\item [\bf 2.]
Suppose $\k\supseteq \Q$ and $R$ is one of $\k[[x]]/J$, $\k\{x\}/J$, $\k\langle x\rangle/J$, $C^\infty(\R^p,0)/J$.
 Then  $I_{N+1}\sseteq Der_\k^{(1)}(R)(f)$ iff $\{f\}+I_{N+1}\sseteq Aut^{(1)}_\k(R)\cdot f$.

Thus   $f$ is  right-$N$-determined iff it is infinitesimally right-$N$-determined.
\eee\ecor
\bpr
\bee[\bf 1.]
\item  By Proposition \ref{Thm.Aut(R).of.weak.Lie.type} the pair
$(Der_\k^{(1)}(R),Aut^{(1)}_\k(R))$ is of weak Lie type.
 For the chosen filtration the condition
$I_{N+1} \sseteq Der_\k^{(1)}(R)(f)$ implies $I_{N+k} \sseteq Der_\k^{(k)}(R)(f)$ for any $k \geq 1$.
 Then part {\em 1.} of Theorem \ref{Thm.Finite.Determinacy.weak.Lie.pairs} gives:
\beq
\{f\}+I_{2N+1-ord(f)}\sseteq \overline{Aut^{(N+1-ord(f))}_\k(R)\cdot f}.
 \eeq
Now apply proposition \ref{Thm.Use.of.Approximation} to remove the closure.

Part {\em ii.} is proved similarly, using part {\em 2.} of Theorem \ref{Thm.Finite.Determinacy.weak.Lie.pairs}.

\item
By Theorem \ref{Thm.Aut(R).pointwise.Lie} the  pair $(Der_\k^{(1)}(R),Aut^{(1)}_\k(R))$ is of pointwise Lie type.
 Then, by Theorem \ref{Thm.Finite.Determinacy.pointwise.Lie.type} we have:
 $I_{N+1}\sseteq \overline{Der_\k^{(1)}(R)(f)}$  iff  $\{f\}+I_{N+1}\sseteq \overline{Aut^{(1)}_\k(R)\cdot f}$.
  Finally, apply propositions \ref{Thm.Use.of.Approximation}, \ref{Thm.Approximation.for.C.infty.case} to get the result.
\epr
\eee
\bex
Suppose $\k$ is a field and $R=\k[[ x]]$, or $\k\{ x\}$ ($\k$ perfect), or $\k\langle x\rangle$, or $C^\infty(\R^p,0)$, filtered by $\{\cm^j=( x)^j\}_j$.
Then $Der_\k^{(1)}(R)(f)=\cm^2\cdot (\di_1 f,\dots,\di_p f)$,
and Corollary \ref{cor5.1} gives the classical criteria (for right equivalence).
\bee[i.]
\item For $\k\in \C,\R$ and  $R=\k[[ x]]$ or $\k\{ x\}$ see either \cite[Theorem 2.23]{Gr.Lo.Sh} or \cite[Theorem 1.2]{Wall-1981}.
If $\k$ is a field of characteristic zero this is \cite[Theorem 3.1.13]{Bou09}.
\item For $R=\k[[ x]]$ with $\k$ algebraically closed of arbitrary characteristic,  our Corollary \ref{cor5.1} 1.i. is \cite[Theorem 3, part 1]{Boubakri.Gre.Mark}.
\eee
\eex
\bex\label{complicated automorphism group example}
\bee[i.]
\item More generally, suppose  $\k$ is a field and  $R$ is one of $\k[[ x]]/J$, $\k\{ x\}/J$ ($\k$ perfect), $\k\langle x\rangle/J$,
$C^\infty(\R^p,0)/J$, filtered by $\{\cm^j\}$. Then the assumptions of \ref{Assumptions.on.the.ring} are satisfied and
we obtain the determinacy criteria for functions on singular germs. For $R=\C\{ x\}/J$ this was obtained in \cite{Damon84},
see also \cite[Theorem 2.2. i.]{Bruce-Roberts} and \cite[Proposition 1.4]{Dimca}.
\item
If the ring is non-regular then (even for $\k\supseteq\Q$) finite determinacy does not imply that the singularity is isolated. For example, let
$f(x,y,z,w)=w\in \k[[x,y,z,w]]/(xy)$. Then $f$ is obviously 1-determined
(for $w+h, h \in \cm^2$, apply the coordinate change $w \mapsto w+h$),
but the hypersurface $f^{-1}(0)$ has a non-isolated singularity.
\item
We remark that the group $Aut^{(1)}_\k(R)$ can be rather  small when the ideal $J$ is complicated, and similarly for the module $Der^{(1)}_\k(R)$.
In such cases there are no finitely right determined functions at all.
\eee
\eex
\bex
(Newton filtrations, $\k\supseteq\Q$) Let $R$ be one of $\k[[ x]]$, $\k\{ x\}$, $\k\langle x\rangle$.
 For any Newton diagram $\Ga\sset\R^p_{\ge0}$ (of some element $f\in R$) take the corresponding ideal $I(\Ga)$, generated by monomials lying on or above the diagram. Take a sequence of decreasing Newton diagrams, i.e.
$\Ga_i+\R^p_{\ge 0}\supsetneq \Ga_{i+1}+\R^p_{\ge 0}$.  Then  the ideals $\{I(\Ga_j)\}_i$ form a decreasing filtration.
 Suppose this filtration   satisfies  assumptions \ref{Assumptions.on.the.ring.for.applications}, e.g.
   \bei
   \item To ensure $Der^{(1)}_\k(R)(x)\sseteq (x)^2$ it is enough to assume that no point of $\Ga_1$ lies under the hyperplane $\{\sum x_i=2\}\sset\R^p$.
   \item To ensure $I_N\cdot Der_\k(R) \sseteq Der^{(1)}_\k(R)$, for $N\gg1$, it is enough to assume that the distances between the successive
    diagrams are bounded from above.
   \eei
Then we get the determinacy criterion for Newton filtration (see also \cite{Arnold}, \cite{Wall-1999}, \cite{Boubakri.Gre.Mark2011} for a detailed treatment of piecewise filtrations  and a refined determinacy criterion in \cite[Corollary 4.8] {Boubakri.Gre.Mark2011}).
\eex

\bex
Let   $R$ be one of $\k[[ x]]$, $\k\{ x\}$, $\k\langle x\rangle$, $C^\infty(\R^p,0)$, with filtration $\{\cm^j\}$.
Another way to control the finite determinacy is to bound it by the Milnor number, $\mu(f):=dim_\k R/Der_\k(R)(f)$.
Using the obvious bound  $\cm^{\mu(f)}\sseteq Der_\k(R)(f)$ we get
\[
\cm^{\mu(f)+2}\sset \cm^2\cdot Der_\k(R)(f)\sseteq Der^{(1)}_\k(R)(f).
\]
Therefore we get from Corollary \ref{cor5.1}:
\bee[i.]
\item If $\mu(f)<\infty$ then $f$ is $\left(2\mu(f)-ord(f)+2\right)$-right-determined.
This is \cite[Corollary 1, part 1]{Boubakri.Gre.Mark}.
\item If $\k\supseteq\Q$, (e.g., $\k$ is a field of characteristic zero) then $f$ is $\left(\mu(f)+1\right)$-right-determined.
For $\k=\C$ this is \cite[Corollary 2.24 (1)]{Gr.Lo.Sh}.
If $f$ is $N$-right determined and  $\k\supseteq\Q$ then  $\mu(f)\le\bin{N+n}{n}$.
\item
 For positive characteristic we get: if $\{f\}+\cm^{N+k}\sseteq Aut_{\k}^{(k)} (R)\cdot f$ for any $k\ge1$, then $\mu(f)\le\bin{2N-ord(f)+n}{n}$.  In particular, $f$ has at most an isolated singularity.

 Note that in \cite[Theorem 4.1]{Boubakri.Gre.Mark}  (with a correction in \cite[Theorem 4.13]{Greuel-Pham.2017}) the conclusion was proved with the weaker assumption: $\{f\}+\cm^{N+1}\sseteq Aut_{\k}^{(1)} (R)\cdot f$, and hence, finite determinacy is equivalent to isolated singularity. The proof
 in \cite{Greuel-Pham.2017} required
a special construction in positive characteristic, and it is not clear if such a construction is available in our general setting.

 \eee
\eex

\subsection{Right (in)determinacy of germs of maps}
We consider the free $R$-module $R^{n}$, which can be identified with the space of maps $Spec(R)\to (\k^n,0)$. The group  $G=Aut_\k(R)$ acts  componentwise on
$R^{n}$. In this section we assume   $\k\supseteq\Q$.
For simplicity assume that the filtration of $R$ is by powers of the maximal ideal.

\bprop\label{Thm.no.finite.determ.for.maps}
 Let $(R,\cm)$ be a local $\k$-algebra, $\k\supseteq\Q$, filtered by $\cm^j$. Suppose  that the $\cm$-adic completion $\hR$ of $R$ is a complete Noetherian ring of positive Krull dimension with
 $\hR/\cm\hR =\k$. If $n>1$ and a tuple $f = (f_1,\cdots, f_n) \in\cm\cdot R^n$ is  $Aut^{(1)}_\k(R)$-finitely determined then the following holds.
\bee[1.]
\item For any associated prime $\cm\neq \cp\in Ass(R)$ holds: $dim(R/\cp)\ge n$
 and the image of $f$ in $(R/\cp)^n$ is 1-determined.
\item The images of $\{f_i\}$ in the vector space $\cm/\cm^2$ are linearly independent.
 \eee
\eprop
In particular we get: if $f_i\in \cm^2$ for some $i$,   then $f$ is not  $Aut^{(1)}_\k(R)$-finitely determined.

(See also Proposition 4.5 and Corollary 4.6 in \cite{Greuel-Pham.2018}  for the case $R = \k[[x]]$).

\bpr
\bee[\bf Step 1.]
\item
We reduce the proof to the particular case: $R$ is a  complete local Noetherian domain (over a field of zero characteristic).
\bee[a.]
\item
   Take the $\cm$-adic completion, $R\to \hR$.
 If $ f \in R^n$ is $N$-determined for  $Aut^{(1)}_\k(R)$,
then  $\hat f \in \hR^n$ is
 N-determined  for $Aut^{(1)}_\k(\hR)$. Indeed, fix any $\hat{h}\in \cm^{N+1}\cdot \hR^n$ and take a representing sequence $\{h_j\in \cm^{N+1}\cdot R^n\}_j$.
 It satisfies: $h_{i+j}-h_j\in \cm^j\cdot R^n$ and
   the image of this sequence in $\hR^n$ converges to $\hat{h}$.
   Then $\hf+\hat{h}_j\in Aut^{(1)}_\k(\hR)(\hf)$, thus $\hf+\hat{h}\in Aut^{(1)}_\k(\hR)(\hf)$, by completeness.
 Thus it is enough to prove the statement for $R$ -  complete  local Noetherian.

\item
Let $\{\cp_i\}$ be the associated primes of $R$. Any automorphism
$\phi\in Aut^{(1)}_\k(R)$ acts on the set $\{\cp_i\}$. Moreover,
 by the unipotence, $\phi$ preserves each $\cp_i$.
 Thus, for any such $\cp$, we have the natural homomorphism of groups, $Aut^{(1)}_\k(R)\to Aut^{(1)}_\k(R/\cp)$. Therefore,
if $ f \in R^{n}$ is $N$-determined for $Aut^{(1)}_\k(R)$, then its image
in $(R/\cp)^{n}$ is $N$-determined for  $Aut^{(1)}_\k(R/\cp)$.
 As $\cp$ is a prime ideal, $R/\cp$ is a domain. It is of positive dimension iff $\cp\neq\cm$.
 Thus it is enough to prove the statement for $R$ a local complete Noetherian domain.

\eee

\item
\bee[a.]
\item
 Let $R$ be a local complete Noetherian domain, containing $\Q$. The pair
$(Der^{(1)}_\k(R),Aut^{(1)}_\k(R))$ is of  Lie type,
   by Theorem \ref{Thm.Aut(R).pointwise.Lie}. Thus (Theorem \ref{Thm.Finite.Determinacy.pointwise.Lie.type}) the finite determinacy of $ f \in R^n$ implies:
    $Der^{(1)}_\k(R)( f )\supseteq \cm^{N+1}\cdot R^n$, for some $N$. In particular,      $Der_\k(R)( f )\supseteq \cm^{N+1}\cdot R^n$,
      i.e. $\cm$ is the support of the module $R^n/Der_\k(R)( f )$ and is the only associated prime of
    $Ann(R^n/Der_\k(R)( f ))$. Hence  $Ann(R^n/Der_\k(R)( f ))$ is $\cm$-primary and  its height must be $dim(R)$.

Fix some generators $\{D_\al\}$ of $Der_\k(R)$ and consider the matrix $\{D_\al(f_i)\}_{i,\al}$, which is a presentation matrix of  $R^n/Der_\k(R)( f )$.
By \cite[Proposition 20.7]{Eisenbud} $Ann(R^n/Der_\k(R)( f ))$ and the $n$'th
determinantal ideal  $I_n[D_\al(f_i)]$ have the same radical.
We want to check
whether    $I_n[D_\al(f_i)]\supseteq\cm^{\tilde{N}}$ for some $\tilde N$.

\item
 Recall that for a local complete Noetherian domain the rank of the module of derivations equals $dim(R)$.
  Therefore for   $D_1\dots D_{dim(R)+1}\in Der_\k(R)$
  holds: $\sum a_i D_i=0$ for some $a_i\in R$, where not all $a_i$'s are zero. Therefore any $dim(R)+1$ columns of $[D_\al(f_i)]$ are $R$-linearly dependent.
   Thus it suffices to consider only some block of $dim(R)$ columns. So, we assume the matrix $[D_\al(f_i)]$ is of size $n\times dim(R)$.

If $dim(R)<n$ then $I_n[D_\al(f_i)]=0$ and no power of $\cm$ is $0$ since  $dim(R)>0$.
Thus we assume $dim(R)\ge n$. Evaluate this matrix at the origin and check its rank, over the field $R/\cm$.
 Suppose $r:=rank[D_\al(f_i)]|_0<n$, then the matrix is left-right equivalent to $\one_{r\times r}\oplus B$, where $B\in Mat_{(n-r)\times(dim(R)-r)}(\cm)$.
 Thus $I_n[D_\al(f_i)]=I_{n-r}(B)$. But the later ideal  has height at most $(dim(R)-r)-(n-r)+1<dim(R)$, see \cite[(2.1)]{Bruns-Vetter}.
  Thus $I_{n-r}(B)$ cannot contain any power of $\cm$.

Therefore the finite determinacy of $f$ implies:
$rank[D_\al(f_i)]|_0=n$. Then $[D_\al(f_i)]$ is left-right
equivalent to $[\one_{n\times n}|\ \zero_{n\times (dim(R)-n)} ]$, the unit and the zero matrices.
 But then $Der_\k(R)( f )= R^n$ and hence
 \beq
 Der^{(1)}_\k(R)( f )\supseteq  \cm\cdot Der_\k(R)( f )\supseteq \cm\cdot R^n.
 \eeq
  Thus, by theorem \ref{Thm.Finite.Determinacy.pointwise.Lie.type},
  $f\in R^n$ is 1-determined.
\item
Finally, if over the initial ring the images of $\{f_i\}$ in $\cm/\cm^2$ are not $\k$-linearly independent, then they are not linearly
 independent also after the completion and the projection to $R/\cp$, for an associated prime $\cp$. But then  $I_n[D_\al(f_i)]$ cannot contain any power of
  $\cm$. (The row operations correspond to $\k$-linear operations on $\{f_i\}$, thus we can assume $f_1\in \cm^2$ and get: $rank[D_\al(f_i)]|_0<n$.)
  \epr\eee

\eee
\bex
\bee[i.]
\item
 Suppose $(R,\cm)$ is regular Noetherian and a tuple $(f_1\dots f_n)$ is finitely $Aut_\k(R)$ determined, $n>1$. Then the $f_i$ are
  linearly independent mod $\cm^2$ and hence can be extended to  local coordinates on $Spec(R)$, i.e. a regular system of parameters of $R$ (by Nakayama).
 For $\C\{x\}$ this was proved by Mather, see proposition 2.3 of  \cite{Wall-1981}.

\item
 For non-regular rings the tuple   $\{f_i\}$ can be finitely determined without even being a regular sequence.
  For example, let $R=k[[x,y]]/ ( (x)\cdot(x,y) )$.  Here x, y are multi-variables.
Then the sequence $\{y_i\}$ is linearly independent mod $\cm^2$, hence finitely determined (in fact 1-determined).
 But $y_i$ is a zero-divisor in $R$ and  $\{y_i\}$ is not a regular sequence.
\eee
\eex
\beR
The  statement of  Proposition \ref{Thm.no.finite.determ.for.maps} is restricted to the case $char(\k)=0$,
   because only in this case we have established the implication (Theorem \ref{Thm.Finite.Determinacy.pointwise.Lie.type})
 \beq
\overline{G^{(1)}z}\supseteq \{z\}+M_{N+1} \Rrightarrow \overline{T_{(G^{(1)},M)}z}\supseteq  M_{N+k},\ \text{for some }k<\infty.
\eeq
\eeR

\bigskip

\subsection{Contact determinacy of germs of maps}
In this subsection we consider the determinacy under the action of the contact group, $\cK:=GL(n,R)\rtimes Aut_\k(R)\circlearrowright R^{n}$.
In this case the tangent space  $T_{(\cK^{(1)},M)}$ is generated by $Der_\k^{(1)}(R)$ and $End^{(1)}_R(R^n)=Mat_{n\times n}(I_1)$, where
$\{I_j\}$ is the filtration of $R$.
 Denote the ideal generated by a tuple $f=(f_1,\dots,f_n)\in R^n$ by $(f)$. Denote by $Der_\k^{(1)}(R)(f)\sseteq R^n$  the submodule obtained by
  applying the derivations to the tuple $f$.
 The pair $(T_{(\cK^{(1)},M)},\cK^{(1)})$
is of pointwise (weak) Lie type,
see Examples  \ref{Ex.GL.times.GL.times.Aut.of.weak.Lie.type}, \ref{Ex.Gl(m,R)GL(n,R)Aut(R).char.0.pointwise.Lie.type}.
As in the case of right equivalence we have:

\bcor \label{cor5.9}
Suppose $R$ with filtration $\{I_j\}$ satisfies the
assumptions \ref{Assumptions.on.the.ring.for.applications}.
\bee[1.]
\item
\bee[i.]
\item If $I_{N+1}\cdot R^{n}\sseteq Der_\k^{(1)}(R)(f)+ (f)\cdot I_1\cdot  R^n$ then $(f)\cdot R^n +I_{2N+1-ord(f)}\cdot R^{n}\sseteq \cK^{(N+1-ord(f))}\cdot f$.

(Thus $f$ is $(2N-ord(f))$-contact-determined.)
\item
Suppose
$(f)\cdot R^n+I_{N+k}\cdot R^n\subseteq \cK^{(k)}\cdot f$ for any $k\ge1$, then $I_{2N+1-ord(f)}\cdot R^n\sseteq Der^{(1)}_\k(R)(f)+(f)\cdot I_1\cdot R^n$.

(Thus $f$ is infinitesimally $(2N-ord(f))$-contact-determined.)
\eee
\item Let $\k\supseteq \Q$.
Then $I_{N+1}\cdot R^{n}\sseteq Der_\k^{(1)}(R)( f)+(f)\cdot I_1\cdot R^n$ \ iff \
$(f)\cdot R^n +I_{N+1}\cdot R^{n}\sseteq \cK^{(1)}\cdot f$.

 ($(N+1)$-contact-determinacy of $f$ vs infinitesimal $(2N+1)$-contact-determinacy of $f$)
\eee
\ecor
The proof is the same as for corollary \ref{cor5.1}.

\bex \label{ex5.10} (For function germs, $n=1$)
 For $R$ one of $\k[[ x]]$, $\k\{ x\}$, $\k\langle x\rangle$, $C^\infty(\R^p,0)$, with filtration $\{\cm^j\}$, we get the classical statements, with $Jac(f)$ the Jacobian ideal of $f$.
\bee[i.]
\item($\k\supseteq \Q$) $\cm^{N+1}\cdot  R^{n}\sseteq \cm^2\cdot Jac( f )+Mat_{n\times n}(\cm)\cdot  f $ iff
$(f)+\cm^{N+1}\cdot  R^{n}\sseteq \cK^{(1)}( f )$.

 For $\k=\C$ see \cite[Theorem 1.2]{Wall-1981}, \cite[Theorem 2.23]{Gr.Lo.Sh}.
For $\k$ a field of characteristic zero this is \cite[Theorem 3.1.13]{Bou09} and \cite[Theorem 4.2]{Boubakri.Gre.Mark}.

\item In  positive characteristic, for $R=\k[[ x]]$,
   part {\em 2.i} is \cite[Theorem 3, part 2]{Boubakri.Gre.Mark}, see
also \cite[Lemma 2.6]{GrKr.90}.
Part {\em 2.ii.} was claimed in \cite{Boubakri.Gre.Mark}, but the proof was incomplete.
 The full (quite involved) proof was given in Theorem 4.6 and 4.8 in \cite{Greuel-Pham.2017}, for the case that $\k$ is infinite,
  but assuming only: $(f)+I_{N+1}\subseteq \cK^{(1)}\cdot f$.
\eee

For maps, $n\ge1$, cf.  lemma 6.22 and theorem 6.27 in \cite{Dimca.book}.
\eex
\bex (For maps, $n\ge1$) Let $R$ be one of $\k[[ x]]/J$, $\k\{ x\}/J$, $\k\langle x\rangle/J$, $C^\infty(\R^p,0)/J$, with filtration $\{\cm^j\}$. We
 get the determinacy criteria of maps from $Spec(R)$ to $(\k^n,0)$.
 \bee[i.]
 \item For $\k\supseteq\Q$ and $J=0$ (regular rings) the condition for finite $\cK$-determinacy is in general not implied by the condition that the variety $V(f)$ has an isolated singularity, except if $(f)$ is a minimal generating set of a complete intersection ideal (in which case both conditions are equivalent, cf. \cite [Theorem 4.6]{Greuel-Pham.2017} and
 Example \ref{ICIS}). For example,   take the map $f = (xy,xz,yz)$,  then $V(f)$ has an isolated singularity but the map $f$ is not finitely determined by Corollary \ref{cor5.9}
  \item
As we have observed in Example \ref{complicated automorphism group example}, if $Spec(R)$ has
a ``complicated" singularity, then the group $Aut_{\k} (R)$ can be
very small, and there might be no finitely determined germs at all.
\eee
\eex
\bex Let  $R$ be one of $\k[[ x]]$, $\k\{ x\}$, $\k\langle x\rangle$, $C^\infty(\R^p,0)$, filtered by $\{\cm^j\}$.
As before, we can control the finite determinacy of function germs $f \in R$ in terms of  the Tjurina number $\tau(f):=dim_\k\left(R/Der_\k(R)(f)+(f)\right)$.
Use the obvious bound  $\cm^{\tau(f)}\sseteq Der_\k(R)(f) +(f)$ to get;
\[
\cm^{\tau(f)+1}\sset \cm\cdot Der_\k(R)(f)+\cm\cdot(f)\sseteq Der^{(1)}_\k(R)(f)+\cm\cdot(f).
\]
Therefore we get:
\bee[i.]
\item If $\tau(f)<\infty$ then $f$ is  $\left(2\tau(f)-ord(f)+2\right)$-contact-determined.
This is \cite[Corollary 1, part 2]{Boubakri.Gre.Mark}.
\item If $\k\supseteq\Q$, (e.g., $\k$ is a field of characteristic zero) then
$f$ is  $\left(\tau(f)+1\right)$-contact-determined. For $\k=\C$ this is \cite[Corollary 2.24]{Gr.Lo.Sh}, for $\k$ a field of
zero characteristic this is \cite[Theorem 3.1.13]{Bou09}.
\item
If $k\supseteq\Q$ and  $f\in R$ is $N$-contact determined then $\tau(f)\le \bin{ N+n}{n}$.
\eee
\eex

\bex
(Finite determinacy of complete intersection ideals)\label{ICIS}
 Let  $R$ be one of $\k[[ x]]/J$, $\k\{ x\}/J$, $\k\langle x\rangle/J$, $C^\infty(\R^p,0)/J$, filtered by $\{\cm^j\}$.
 A complete intersection ideal $I\sset R$ is called $N$-determined if for any minimal
generating sequence, $f=(f_1\dots f_n)$, and any sequence $(g_1\dots g_n)\in \cm^{N+1}$ holds: $(f)\stackrel{Aut^{(1)}_\k(R)}{\sim}(f_1+g_1,\dots,f_n+g_n)$.
 An ideal is called infinitesimally $N$-determined   if $\cm^{N+1}R^{n}\sseteq Der_\k^{(1)}(R)(f)+\cm\cdot I\cdot R^n$.
  (This condition does not depend on the choice of $f$.)
 For $\k\supseteq\Q$ Corollary \ref{cor5.9} can be stated:
$I$ is infinitesimally $N$-determined iff it is $N$-determined.
 See also \cite[Theorem 4.6 (2)]{Greuel-Pham.2017}, where this is proved for arbitrary infinite fields.
 For the case $R=\C\{x\}/J$, with $J$ a complete intersection, and $I\sset R$ a principal ideal,  we get proposition 1.3 of \cite{Dimca2}.
\eex

\subsection{Determinacy for maps relative to a germ}\label{determinacy for maps that preserve a subscheme}
Another typical scenario is when the ambient space contains a particular subscheme, $V(\ca)\sset Spec(R)$. Then one studies
the maps of $Spec(R)$ up to right or contact transformations that preserve $V(\ca)$. Thus we use the
group of relative right transformations, $\cR_\ca:=Aut_{\k,\ca}(R)$, and the group of relative contact transformations,
$\cK_\ca:=GL(n,R)\rtimes Aut_{\k,\ca}(R)$, see section \ref{Sec.(w)Lie.type.of.relative.group}.
 (As was explained in \S\ref{Sec.(w)Lie.type.of.relative.group}, this case cannot be reduced to the case $Aut_{\k}(R/\ca)$.)
The pair $(Der^{(1)}_{\k,\ca}(R),Aut^{(1)}_{\k,\ca}(R))$
is of (pointwise) Lie type, and similarly for the pair
$(Der^{(1)}_{\k,\ca}(R)\oplus \cm\cdot R^{n},\cK_\ca^{(1)})$.
Theorems \ref{Thm.Finite.Determinacy.pointwise.Lie.type},   \ref{Thm.Finite.Determinacy.weak.Lie.pairs} give the immediate corollaries:
\bcor
Suppose $\k\supseteq\Q$ and $R$ with filtration $\{I_j\}$ satisfies the assumptions of \ref{Assumptions.on.the.ring.for.applications}.
 For $R=C^\infty/J$ we assume that $\ca$ is analytically generated.
\bee[\bf 1.]
\item ($n=1$, right determinacy) $I_{N+1} \sseteq Der^{(1)}_{\k,\ca}(R)(f)$ iff  $(f)+I_{n+1}\sset \cR^{(1)}_\ca (f)$.
\item ($n\ge1$, contact determinacy) $I_{N+1}R^n\sseteq Der^{(1)}_{\k,\ca}(R)(f)+I_1\cdot(f)\cdot R^n$ iff  $(f)+I_{n+1}\cdot R^n\sset \cK^{(1)}_\ca (f)$.
\eee
\ecor
\bex
The classically studied cases are:
$R=\cO_{(\k^n,0)}$, for $\k$ a field of zero characteristic, and $I\sset R$ - a radical ideal.
 See e.g., \cite[Theorem 2.2]{Janeczko1982}, \cite[Proposition 1.4]{Dimca}, \cite[Theorem 3.5, Corollaries
3.6 and 3.7]{Grandjean}, \cite{Damon84}, \cite[Propositions
2.3, 2.4 and 2.5]{Izumiya-Matsuoka}.
Suppose $\k$ is one of  $\R,\C$ and $R$ is one of  $\k[[ x]]$, $\k\{ x\}$, $\k\langle x\rangle$, $C^\infty(\R^p,0)$.
Then we recover, e.g., \cite[Theorem 3.6]{Orefice.Tomazella} (for $\cR_I$) and \cite[Lemma 3.11]{Orefice.Tomazella} (for $\cK_I$).
The converse statement (finite determinacy implies large tangent space) is e.g., \cite[Theorem 3.5]{Orefice.Tomazella}.
\eex

\subsection{Relative determinacy for non-isolated singularities of function germs}\label{determinacy for non-isolated singularities}
Suppose an element $f\in R$ defines a non-isolated singularity, i.e., one of the ideals $Der_\k(R)(f)\sseteq R$,  $(f)+Der_\k(R)(f)\sseteq R$ has infinite colength.
Then no finite determinacy is possible for the filtration $\{\cm^j\}$. In such cases one restricts the possible deformations, taking only
those that preserve the singular locus of $f$ (with its multiplicity). This corresponds to filtration $\{\cm^j\cdot \ca\}$.
Here the ideal $\ca$ is usually non-radical, it    defines the relevant singularity scheme.
Accordingly, instead of the groups  $Aut_{\k}(R)$, $GL_R(n)\rtimes Aut_{\k}(R)$, one considers
  the subgroups $\cR_\ca=Aut_{\k,\ca}(R)$, $\cK_\ca=GL_R(n)\rtimes Aut_{\k,\ca}(R)$,
see \S  \ref{Sec.(w)Lie.type.of.relative.group}.
  As before, one has the $rel(\ca)$ notions of determinacy.
 For simplicity we restrict to the case $\k\supseteq\Q$.
\bcor
Suppose $\k\supseteq\Q$ and $R$, $\{I_{N+1}=\cm^N\cdot \ca\}$ satisfy the condition \ref{Assumptions.on.the.ring.for.applications}.

\bee[1.]
\item ($n=1$) $I_{N+1} \sseteq Der^{(1)}_{\k,\ca}(R)(f)$ iff  $(f)+I_{N+1} \sset \cR^{(1)}_\ca (f)$.
\item ($n\ge1$) $I_{N+1}\cdot R^n\sseteq Der^{(1)}_{\k,\ca}(R)(f)+\ca\cdot(f)\cdot R^n$ iff  $(f)+I_{N+1}\cdot R^n\sset \cK^{(1)}_\ca (f)$.

\eee
\ecor
The module of $\ca$-logarithmic derivations, $Der_{\k,\ca}(R)$, is in general complicated, but it often contains a simpler module,
\beq
Der_{\k,\ca}(R)\supseteq \sqrt{\ca}\cdot Der_\k(R)+Ann_{Der_\k(R)}(\ca).
\eeq
(The latter summand here denotes all the derivations that annihilate $\ca$.) This leads to a weaker statement, but with the condition easier to check.
\bex
Suppose $R=\k[[x_1,\dots,x_n]]$ and take the ideal $\ca=(x_1,\dots,x_l)^q$.
Then
\[
Der_{\k,\ca}(R)=\langle\di_{l+1},\dots,\di_n\rangle+(x_1,\dots,x_l)\cdot\langle\di_1,\dots,\di_l\rangle.
\]
 Then for  $f\in (x_1,\dots,x_l)^q\smin (x_1,\dots,x_l)^{q+1}$ we get:
\bee[i.]
\item Suppose $\k\supseteq\Q$.
If $\cm^2\langle\di_{l+1},\dots,\di_n\rangle(f)+\cm\cdot\sqrt{\ca}\cdot\langle\di_1,\dots,\di_n\rangle(f)\supseteq I\cdot\cm^{N+1}$
 then $f$ is $N$-right$^{rel(\ca)}$-determined.
(And similarly for the contact determinacy.) For the case $\ca=(x_1,\dots,x_l)^2$ this goes in the style of results
of \cite[Theorem 6.5 and Corollary 6.6]{Pellikaan88}, see also \cite[Proposition 1.5 and Corollary 1.6]{Siersma83} and \cite[Theorem 3.5 and Corollary
3.6]{Grandjean}. 
\item For an arbitrary $\k$ we have:
if $\cm^2\langle\di_{l+1},\dots,\di_n\rangle(f)+\cm\cdot\sqrt{\ca}\cdot\langle\di_1,\dots,\di_n\rangle(f)\supseteq \ca\cdot\cm^{N+1}$ then
$f$ is $(2N-ord(f))$-right$^{rel(\ca)}$-determined.   (And similarly for the contact determinacy.)
This is \cite[Theorem 3.2]{Hengxing.Jingwen}.
\eee
\eex

\subsection{Relative algebraization}
 If $f\in \k[[x]]$ is finitely ($\cR$ resp. $\cK$) determined then, in particular, it is  ($\cR$ resp. $\cK$) equivalent to a polynomial.
 Such an algebraization does not hold for non-isolated critical points resp. singularities. For example, $f(x,y,z)=xy(x+y)(x-zy)(x-e^zy)\in \C\{x,y,z\}$
 is not $\cK$-equivalent to a polynomial, see example 14.1 of \cite{Whitney}. (More references are in \cite{Kucharz}.)
 Certain non-isolated singularities can still be converted to a polynomial, \cite{Moehring.van.Straten}.

We prove now that in the non-isolated case $f$ can be converted to a  polynomial ``in the direction transversal to the critical resp. singular locus".
\bprop
Let $\k$ be a field and $f\in (x)^2\sset\k[[x]]$. Supose the ideal $I:=\sqrt{Jac(f)}$ (resp. $I:=\sqrt{(f)+Jac(f)}$)
 is of height $c$.  Suppose the projection $V(I)\to V(x_1,\dots,x_c)\sset Spec\ \k[[x]]$ by $(x_1,\dots,x_p)\to(x_{c+1},\dots,x_p)$ is a finite morphism.
Then $f$ is $\cR$ (resp. $\cK$)-equivalent to an element of
 $\k[[x_{c+1},\dots,x_p]][x_1,\dots,x_c]$.
\eprop
This result extends (in the irreducible case)  Theorem 1.1 of \cite{Kucharz} to an arbitrary field.
\bpr
 Denote $I:=\sqrt{Jac(f)}$, resp. $I:=\sqrt{(f)+Jac(f)}$. As $I$ is finitely-generated we have: $I^{N+1}\sseteq I\cdot Jac(f)$,
  resp. $I^{N+1}\sseteq I\cdot (f)+Jac(f)$, for $N\gg1$.

 Take the filtration of $\k[[x]]$ by $\{I^j\}$. Note that $\{I^j\}$ satisfies Assumptions \ref{Assumptions.on.the.ring.for.applications}, and (in the case $\k\not\supseteq \Q$) Assumptions  \ref{Assumptions.on.the.ring}.
 Thus $I^{N+1}+\{f\}\sseteq \cR^{(1)}f$, resp. $I^{N+1}+\{f\}\sseteq \cK^{(1)}f$.

\

As the scheme $V(I)$ is of co-dimension $c$ and the projection  $V(I)\to V(x_1,\dots,x_c)$ is finite, the intersection
 $V(I)\cap  V(x_{c+1},\dots,x_p)$ is a zero-dimensional scheme. Therefore $I+(x_{c+1},\dots,x_p)\supseteq (x_1,\dots,x_p)^\tN$, for some $\tN\gg1$. But then we have also  $I^N+(x_{c+1},\dots,x_p)\supseteq (x_1,\dots,x_p)^{N\cdot\tN}$.
 And this implies $I^N+\k[[x_{c+1},\dots,x_p]][x_1,\dots,x_c]=\k[[x_1,\dots,x_p]]$, for any $N$.

Combining the two steps we have:
$
\cR^{(1)}f\cap \k[[x_{c+1},\dots,x_p]][x_1,\dots,x_c]\neq\varnothing$, or  respectively
$\cK^{(1)}f\cap \k[[x_{c+1},\dots,x_p]][x_1,\dots,x_c]\neq\varnothing$.
\epr

\bex (Generalizing Corollary 1.3 of \cite{Kucharz} to an arbitrary field)
Any element $f\in \k[[x]]$ is $\cK$-equivalent to an element of $\k[[x_1,\dots,x_{p-2}]][x_{p-1},x_p]$. Indeed, take the irreducible decomposition, $f=\prod f_i^{n_i}$. The scheme $Sing(V(\prod f_i))$ is of codimension at least 2. Now apply the proposition.

In particular, any element $f\in \k[[x_1,x_2]]$ is $\cK$-equivalent to a polynomial. If $char(\k)=0$ then the critical and singular loci coincide, i.e. $\sqrt{Jac(f)}=\sqrt{(f)+Jac(f)}$, and thus  $f\in \k[[x_1,x_2]]$ is $\cR$-equivalent to a polynomial.
\eex

\subsection{Finite determinacy of matrices}
Take as $M$ the $R$-module of matrices, $\Mat$, with the filtration $\{Mat_{m\times n}(I_j)\}$. In this section the group $G$ will be one of $GL(m,R)$, $GL(n,R)$,
$Aut_\k(R)$, or their (semi-)direct products. They are of (pointwise)  (weak) Lie type,
  their tangent spaces are written down in  Examples \ref{Ex.GL.times.GL.times.Aut.of.weak.Lie.type} and \ref{Ex.Gl(m,R)GL(n,R)Aut(R).char.0.pointwise.Lie.type}.
Thus Theorems \ref{Thm.Finite.Determinacy.pointwise.Lie.type}, \ref{Thm.Finite.Determinacy.weak.Lie.pairs}  imply:

\bcor
 Suppose $R$, $\{I_j\}$ satisfy the assumption \ref{Assumptions.on.the.ring.for.applications}.
 \bee[1.] \item
\bee[i.]
\item If $Mat_{m\times n}(I_{N+1})\sseteq T_{(G^{(1)},M)}(A)$ then $\{A\}+Mat_{m\times n}(I_{2N+1-ord(A)})\sseteq G^{(N+1-ord(A))}(A)$

($(2N-ord(A))$-contact-determinacy)
\item
Suppose
$\{A\}+Mat_{m\times n}(I_{N+k})\subseteq T_{(G^{(k)},M)}(A)$, for any $k\ge1$, then
 $Mat_{m\times n}(I_{2N+1-ord(A)})\sseteq T_{(G^{(1)},M)}(A)$.

(infinitesimal $(2N-ord(A))$-contact-determinacy)
\eee
\item For $\k\supseteq \Q$:
$ Mat_{m\times n}(I_{N+1})\sseteq T_{(G^{(1)},M)}(A)$ iff
$\{A\}+Mat_{m\times n}(I_{N+1})\sseteq G^{(1)}(A)$.

 ($N+1$-contact-determinacy vs infinitesimal $N+1$-contact-determinacy)
\eee
\ecor

For $\k$ a field of characteristic zero  this statement was proved in \cite[Corollary 2.9]{B.K.motor}.
For $\k$ a field of positive characteristic and $R=\k[[ x]]$,
part (1.i) of the statement gives  \cite[Theorem 3.2]{Greuel-Pham.2018}.

For $R=\k[[ x]]$, $\k$ an arbitrary field, a more general statement for matrices was proved in \cite[Proposition 4.2]{Greuel-Pham.2018}.

\beR
In the case of matrices the submodule $T_{(G^{(1)},M)}(A)\sset\Mat$ can be rather complicated.
And a bound like $Mat_{m\times n}(I_{N+1})\sseteq T_{(G^{(1)},M)}(A)$
can be difficult to verify. Here one faces a purely commutative algebra question, to compute or bound the support of the quotients module,
$\Mat/T_{(G^{(1)},M)}(A)$, i.e., its annihilator ideal. An algorithm to compute $T_{(G^{(1)},M)}(A)$ is described in \cite{Greuel-Pham.2017b}
 for the filtration $\{\cm^j\}$.
The module $T_{(G^{(1)},M)}$ is usually close to $T_{(G,M)}$, with the simple bound
$\cm\cdot T_{(G,M)}\sseteq T_{(G^{(1)},M)}\sseteq T_{(G,M)}$. And the quotient
\beq
T^1_{(\Mat,G,A)}:=\Mat/T_{(G,M)}(A)
\eeq
is usually better behaved, thus one first studies this quotient. In \cite{B.K.fin.det.1}, \cite{B.K.fin.det.2}, \cite{Kerner}
this quotient was extensively studied for several group actions with useful bounds for the annihilator of the module $T^1_{(\Mat,G,A)}$.
These led to the simple bounds on the order of determinacy and to full control of finite determinacy.
\eeR

\subsection{Determinacy of families}
 Suppose $R$ is an $S$ algebra, $S$ an algebra over a field ${\k}$, e.g., $S={\k}[[t]]$, ${\k}\{t\}$,
 ${\k}[t]$. We consider the elements of $R$, $R^{n}$, $\Mat$ as families of some objects  over the base space $Spec(S)$.
Fixing  a section $ Spec(S)\to Spec(R) $ we have the (local) families of elements in the families of modules, $\{z_t\in M_t\}$.
   The groups  $Aut_\k(R)$,   $GL(n,R)\rtimes Aut_\k(R)$, $GL(m,R)\times GL(n,R)\rtimes Aut_\k(R)$, etc., induce  equivalence of such families.
(The equivalence acts as identity on the base $Spec(S)$ and maps fibers to fibers.)
 Then one speaks about the order of determinacy of ``families of elements inside families of modules, under the action of some group families''.
\begin{eqnarray}
\text{A family $\{z_t\in M_t\}$ is $N$-determined under the action of }
\\\nonumber
\text{the family $\{G_t\}$ if $\{\{z_t\}+M_{N+1,t}\sseteq G_t z_t\}$, for some $N$}.
\end{eqnarray}
Then, as in all the examples of this section, we get criteria for finite determinacy and for bounds of the order of determinacy for right/contact/etc. equivalence of families. Note that the determinacy is in general not semicontinuous in a flat family (cf. Example I.2.2.4 of \cite{Gr.Lo.Sh}), but it can usually be bounded by a semicontinuos invariant.
This follows from the semicontinuity theorem \cite[Proposition 3.4]{Greuel-Pham.2017} and its application in  \cite[Theorem 4.6]{Greuel-Pham.2017}.

\end{document}